\documentclass[11pt,reqno]{amsart}
\usepackage{amsfonts}
\usepackage{mathrsfs}
\usepackage{amsfonts}
\usepackage{amsmath}
\usepackage{stmaryrd}
\usepackage{amssymb}
\usepackage{amsthm}
\usepackage{mathrsfs}
\usepackage{amsfonts}
\usepackage{amscd}
\usepackage{indentfirst}
\usepackage{enumerate}
\usepackage{amsmath,amsfonts,amssymb,amsthm}
\usepackage{amsmath,amssymb,amsthm,amscd}
\usepackage{graphicx,mathrsfs}
\usepackage{a4wide}
\usepackage{amsmath}
\allowdisplaybreaks[4]
\usepackage[colorlinks,linkcolor=blue,anchorcolor=blue,linktocpage=true,urlcolor=blue,citecolor=blue]{hyperref}

\usepackage{color}

\usepackage{empheq}
\allowdisplaybreaks
\usepackage[titletoc]{appendix}

 \numberwithin{equation}{section}
\newtheorem{theorem}{Theorem}[section]
\newtheorem{proposition}[theorem]{Proposition}

\newtheorem{lemma}[theorem]{Lemma}

\newtheorem{remark}[theorem]{Remark}

\begin{document}

\title[Qualitative properties of peak solutions]
{Qualitative properties of peak solutions for Nonlinear Schr\"{o}dinger equations with nearly critical Sobolev exponents }
\author[Z. Liu, S. Tian, H. Xie and P. Yang]{Zhongyuan Liu, Shuying Tian, Huafei Xie, Pingping Yang}

 \address[Zhongyuan Liu]{School of Mathematics and Statistics, Henan University, Kaifeng, 475004,  P. R. China}
 \email{liuzy@henu.edu.cn}

  \address[Shuying Tian]{School of Mathematics and Statistics, Wuhan University of Technology, Wuhan, 430070, P. R. China}
\email{sytian@whut.edu.cn}

\address[Huafei Xie]{Department of Mathematics, Shantou University, Shantou, 515063, P. R. China}
  \email{huafeixie@ccnu.edu.cn}

 \address[Pingping Yang]{School of Mathematics and Statistics, Key Laboratory of Nonlinear Analysis
$\&$ Applications (Ministry of Education),
 Hubei Key Laboratory of Mathematical Sciences, Central China
Normal University, Wuhan, 430079, P. R. China \& Dipartimento SCFMN, $Universit\grave{a}$ degli Studi di Sassari, Via Vienna 2, 07100 Sassari, Italy}
\email{ypp15623175603@mails.ccnu.edu.cn}

%\thanks{Luo was supported  by NSFC grants (No.11701204, 11831009).}
\begin{abstract}
In this paper, we are concerned with  qualitative properties of peak solutions to the following nonlinear Schr\"{o}dinger equations
\begin{equation*}
-\Delta u+V(x)u= u^{p-\varepsilon},\,\,\,u>0\,\,\,\text{in}\,\,\,\mathbb{R}^N,
\end{equation*}
where  $V(x)\in C^1(\mathbb{R}^N) $ is a nonnegative bounded function,  $\varepsilon>0$, $p=\frac{N+2}{N-2}$,  $N\geq6$.
Precisely, we use the blow-up analysis  based on  Pohozaev  identities   to establish the local uniqueness and  Morse index
of the peak solutions obtained by Cao et al. (Calc. Var. Partial Differential Equations, 64: 139, 2025),  provided that $V(x)$ possesses  non-degenerate critical points.
\end{abstract}
%\date{\today}
\maketitle
{\small

\section{Introduction and Main Results}\label{s1}

In this paper, we consider the following nonlinear Schr\"{o}dinger equations
\begin{equation}\label{a1.1}
 -\Delta u +V(x)u=f(u),\quad u>0\quad\text{in}\quad\mathbb{R}^N,
  \end{equation}
  where   $V(x)\in C^1(\mathbb{R}^N) $ is a nonnegative  function, $f(u)$ is a continuous function of $u$.

Problem \eqref{a1.1} arises naturally in search for   standing waves for the nonlinear Schr\"{o}dinger equations
\begin{equation}\label{1.2}
i\hbar\frac{\partial\psi}{\partial t}+\hbar^2\Delta\psi-\widetilde{V}(x)\psi+f(\psi)=0,\quad(t,x)\in\mathbb{R}\times\mathbb{R}^N, \end{equation}
where $\hbar$ is the Plank constant, $i$ is the imaginary unit, $f(e^{i\theta}u)=e^{i\theta}f(u)$ for all $u\in\mathbb{R}$. A standing wave is a solution of the form  $\psi(t,x)=\exp(-iE t/\hbar)u(x)$.
Then $\psi(t,x)$ is a solution of \eqref{1.2} if and only if $u$ satisfies the nonlinear elliptic equation
 \begin{equation}\label{01.1}
 -\hbar^2\Delta u +V(x)u=f(u)\quad\text{in}\quad\mathbb{R}^N, \end{equation}
where $V(x)=\widetilde{V}(x)-E, E$ is a positive constant.
This elliptic  equation has its roots in the study of semiclassical states (i.e. $\hbar\rightarrow0$)  of \eqref{01.1} under various assumptions on the potential $V(x)$ or
 the nonlinearity $f(u)$. For more information concerning the existence, multiplicity and  asymptotics  of semiclassical states of \eqref{01.1}, we refer the interested readers to \cite{ABC,AmMS,BJ1,BJ2,BW1,BW2,Cdny,CN,CNY,CP0,CY,DDMW,DF1,DF2,FW,Gui,Li,Oh1,Oh2,Wang,WY0,ZCZ,ZZ2} and the references therein.

 % \vspace{7cm}

%Problem \eqref{01.1} originated from searching the standing wave of the following nonlinear Schr\"{o}dinger equation
%\begin{equation}\label{1.2}
%i\hbar\frac{\partial\psi}{\partial t}+\hbar^2\Delta\psi-\widetilde{V}(x)\psi+f(\psi)=0,\quad(t,x)\in\mathbb{R}\times\mathbb{R}^N, \end{equation}
%where $\hbar$ is the Plank constant, $i$ is the imaginary unit, $f(e^{i\theta}u)=e^{i\theta}f(u)$ for all $u\in\mathbb{R}$.   The standing wave is the solution to equation \eqref{1.2} with the form  $\psi(t,x)=\exp(-iE t/\hbar)u(x)$.
%Then $\psi(t,x)$ satisfies equation \eqref{1.2} if and only if $u$  solves  nonlinear elliptic equation
% \begin{equation}\label{1.3}
% -\hbar^2\Delta u +V(x)u=f(u)\quad\text{in}\quad\mathbb{R}^N, \end{equation}
%where $V(x)=\widetilde{V}(x)-E, E>0$ is a constant.
%
%In recent years, Schr\"{o}dinger equations  and its variants have been widely studied
%which mainly focus on the existence, multiplicity and qualitative properties of semiclassical states of \eqref{1.3}  with some suitable conditions on the potential $V(x)$ or the nonlinearity $f(u)$. For example, Cao and Noussair \cite{CN} used the variational methods and penalization techniques to study the existence and qualitative property of standing wave solutions to \eqref{1.3} as $\hbar\rightarrow0$. For more results related to the problem, see \cite{Cdny,CNY,WY0} and the references therein.

Beyond semiclassical states, there has been a large number of works devoted to
the existence of ground states and higher energy solutions for \eqref{a1.1} (i.e., $\hbar=1$ in \eqref{01.1}) in the last years.
%If $\hbar=1$, problem \eqref{01.1} becomes
%\begin{equation}\label{1.4}
% -\Delta u +V(x)u=f(u)\quad\text{in}\quad\mathbb{R}^N.
%  \end{equation}
% For the  subcritical or critical problems, the methods are usually based on variational or perturbation methods.
%Next, we will introduce some results related to equation \eqref{1.4}.
   Rabinowitz \cite{Ra} showed  that \eqref{a1.1} admits a ground state solutions for compact potential $V(x)$ and
    subcritical  nonlinearity $f(u)$.
  Jeanjean and Tanaka \cite{JT} obtained the existence of a positive solution to \eqref{a1.1} without  Ambrosetti-Rabinowitz superlinear condition on
  $f(u)$.
  Li, Wang and Zeng \cite{LWZ1} utilized   a ``Nehari type" argument to establish the existence of ground state solutions of \eqref{a1.1} under a super-quadratic condition on  $f(u)$.
  Wei and Yan \cite{WY2}  employed  the finite reduction method  to get
   infinitely many non-radial solutions for \eqref{a1.1} with $f(u)=u^q$, $1<q<\frac{N+2}{N-2}$,  whose energy can be made arbitrarily large.
 Recently,   Guo, Musso, Peng and Yan \cite{GMPY22} showed that the solutions obtained in \cite{WY2}
are non-degenerate in a suitable symmetric space and then  construct new peak solutions of \eqref{a1.1} by this non-degenerate result.

Pan and  Wang \cite{PW} investigated
\begin{equation}\label{1.1}
 -\Delta u +V(x)u=u^{p-\varepsilon},\quad u>0\quad\text{in}\quad\mathbb{R}^N,
  \end{equation}
where $p=\frac{N+2}{N-2}$ and $\varepsilon>0$. They obtained
asymptotic behavior of ground states to   \eqref{1.1} and
 proved that ground state solutions blow up and concentrate at a single point as $\varepsilon\rightarrow0$.
Furthermore,  Wang \cite{W} showed that   the location of   blow-up point must be  a global minimum point of $V(x)$.
 Micheletti and Pistoia \cite{MP} established the existence of  peak solutions to  \eqref{1.1} for $N\geq7$,
  which blow up at the  stable critical points of $V(x)$, under the certain assumptions on  $V(x)$  and  that $\|V\|_{L^{\frac{N}{2}}(\mathbb{R}^N)}$ is small enough.
Hirano, Micheletti and Pistoia \cite{HMP} constructed positive and sign-changing peak solutions for  \eqref{1.1}.
Subsequently, Cao, Liu and Luo \cite{CLl1} combined the local Pohozaev-type identities   with the reduction argument in suitable
weighted-$L^\infty$ spaces to get the existence of multi-peak solutions to \eqref{1.1} with $ N\geq5$.
Unlike \cite{HMP,MP}, their  approach does not require that
 $\|V\|_{L^{\frac{N}{2}}(\mathbb{R}^N)}$ is small.
For more results of \eqref{1.1} with $\epsilon\leq0$, we refer the interested readers to the references \cite{BC,CWY,CP,DDMW,MMP,MPV,PWY,VW}.

\vskip 0.1cm
On the other hand, qualitative characteristics of peak solutions have long been a focus of extensive research.
%Recently, Liu, Luo  and Peng \cite{LLmo} studied the non-degeneracy and uniqueness of the bubbling solutions to the almost critical H\'{e}non problem in $\mathbb{R}^N$.
Grossi \cite{G02} studied \eqref{01.1} for power type nonlinearity, that is,
\begin{equation}\label{a1.4}
 -\hbar^2\Delta u +V(x)u=u^q, \;\;u>0\quad\text{in}\quad\mathbb{R}^N,
\end{equation}
where $1<q<\frac{N+2}{N-2}$.  He established  the uniqueness and multiplicity results of single peak solutions of \eqref{a1.4} for  a suitable class of potential $V(x)$.
Cao and Heinz \cite{CH} used the topological degree method to obtain the uniqueness of multi-peak solutions of \eqref{a1.4} concentrating near nondegenerate critical points of $V(x)$  for $\hbar$ small.
Grossi and Servadei \cite{GS} computed the Morse index of single-peak solutions of \eqref{a1.4} under some assumptions on the (possibly degenerate) potential $V(x)$.
Later on, Cao, Li and Luo \cite{CLiL} used   local Pohozaev-type identities to prove the local uniqueness of  multi-peak solutions of \eqref{a1.4} with weaker conditions on $V(x)$. Moreover,
Luo, Tian and Zhou \cite{LTZ} studied the local
uniqueness of peak solutions of \eqref{a1.4} for a special non-admissible potential $V(x)$, which has non-isolated critical points.
Luo, Pan, Peng and Zhou \cite{LPPZ} addressed \eqref{a1.4} with very degenerate potentials and gave the existence and uniqueness of multi-peak solutions for \eqref{a1.4}.
 In very recent work, Luo, Pan and Peng \cite{LPP}  extended Grossi and Servadei's work and computed the Morse index of  multi-peak solutions to \eqref{a1.4}  provided that the critical points
 of $V(x)$ are non-isolated and degenerate. For recent
works in this direction, we refer the reader to \cite{CPY,GLY,GLW,LPW20,LPWY,LPZ23}.

Motivated by the aforementioned works, it is natural to investigate the qualitative properties of multi-peak solutions for \eqref{1.1}.
%The main purpose of  the present paper is to investigate the uniqueness and Morse index of multi-peak solutions of \eqref{1.1} obtained in \cite{CLl1}.
The main purpose of  the present paper is to establish the uniqueness and  the Morse index of peak solutions to \eqref{1.1} constructed in \cite{CLl1}.
To state the main results, we need to give some notations and well-known facts. Define
\begin{equation*}
U_{\mu,\xi}(y)=\alpha_{N}
\frac{\mu^\frac{N-2}{2}}{(1+\mu^2|y-\xi|^2)^{\frac{N-2}{2}}}, \;\;\mu>0,\;\;
\xi\in\mathbb{R}^N,
\end{equation*}
where $\alpha_{N}=\left(N(N-2)\right)^{\frac{N-2}{4}}$.
Then it is known from \cite{CGS} that the functions $U_{\mu,\xi}$
are the only solutions to the following elliptic equation with critical exponent
\begin{equation*}%\label{U11}
-\Delta u=u^{p}, \;\;\;u>0 \;\;\;\text{in}\;\; \mathbb{R}^N. \end{equation*}
Letting $\xi_i$, $i=1,\cdots,k$ be $k$ different  points in $\mathbb{R}^N$,  then we set $d=\min\limits_{i\neq j,i,j=1,\cdots,k}|\xi_i-\xi_j|$ and $\delta=\frac{d}{10}$.
Choose $\eta(x)$ is a  smooth cut-off function satisfying
$\eta(x)=1$ if $|x|\leq\delta$, $\eta=0$ if $|x|\geq2\delta$ and $0\leq\eta\leq1$. Letting $\eta_i(x)=\eta(x-\xi_i)$, $i=1,\cdots,k$, then we define
\[
W_{\mu_i,\xi_i}=\eta_iU_{\mu_i,\xi_i},\quad \boldsymbol{W_{\mu,\xi}}=\sum_{i=1}^kW_{\mu_i,\xi_i},
\]
where $\boldsymbol{\mu}=(\mu_1,\cdots,\mu_k)\in\mathbb{R}^k$, $\boldsymbol{\xi}=(\xi_1,\cdots,\xi_k)\in\mathbb{R}^{kN}$.

%%Now we  present the existence result for multi-peak solutions to \eqref{1.1} established by Cao, Liu and Luo in \cite{CLl1}.\vskip 0.1cm
%%Now let us recall the existence result for multi-peak solutions to \eqref{1.1} established by Cao, Liu and Luo in \cite{CLl1}.\vskip 0.1cm
We first recall the existence result for multi-peak solutions to \eqref{1.1}, as established by Cao, Liu and Luo in \cite{CLl1}.\vskip 0.1cm

\noindent\textbf{Theorem  A}
 \emph{
Let $N\geq5$. Assume that $V(x)$ is a bounded nonnegative $C^1$ function and $\xi^*_i$, $i=1,2,\cdots,k$ are the $k$  stable critical points of $V(x)$ such that $V(\xi_{i}^*)>0$. Then  there exists a positive constant
$\varepsilon_0$ such that for $\varepsilon\in(0,\varepsilon_0)$,  problem \eqref{1.1}  admits a $k$-peak solution $u_{\varepsilon}$ with the form
\begin{equation}\label{1a}
u_{\varepsilon}=\boldsymbol{W}_{{\boldsymbol{\mu}_\varepsilon},
{\boldsymbol{\xi}_\varepsilon}}+\phi_\varepsilon
=\sum_{i=1}^k W_{\mu_{\varepsilon,i},\xi_{\varepsilon,i}}+\phi_\varepsilon,
\end{equation}
where $\boldsymbol{\mu}_\varepsilon=
(\mu_{\varepsilon,1},\cdots,\mu_{\varepsilon,k})$ and  $\boldsymbol{\xi}_\varepsilon=
(\xi_{\varepsilon,1},\cdots,
\xi_{\varepsilon,k})$ satisfy, as $\varepsilon\rightarrow0$, $\mu_{\varepsilon,i}\sim\varepsilon^{-\frac{1}{2}}$,  $\xi_{\varepsilon,i}\rightarrow\xi_{i}^*$, $\|\phi_\varepsilon\|_{*}=O(\varepsilon^{\frac{1+\sigma}{2}})$, $\|\cdot\|_*$ is defined in \eqref{aa2.12}}.

%
%Grossi and Pacella \cite{Gp} studied the Morse index of single peak solution for the following equation
%\begin{equation*}
% -\Delta u =N(N-2)u^{p-\epsilon},\quad u>0\quad\text{in}\quad\Omega,
%  \end{equation*} where $\Omega$ is a smooth bounded domain.
%More recently,
%Luo, Pan  and Peng \cite{LPP} compute the Morse index of multi-peak solutions for the Schr\"{o}dinger equation. The qualitative property of bubbling solutions for critical elliptic equations  has been widely studied. One can see  \cite{CPY,DLY,Guo2,LLX,LTX} and references therein. Inspired by them, we also devote to analyze the uniqueness and the Morse index of the multi-peak solution obtained in \cite{CLl1}
%via two kinds of local Pohozaev identities and blow-up analysis.

%However, the theory for (1.1)and (1.4)in general is still in its infancy.
\vskip 0.1cm
Throughout this paper, we always assume that $V(x)$ is a bounded nonnegative $C^1$ function.
Now we   state  the local uniqueness of multi-peak solutions in Theorem A  as follows.
\begin{theorem}\label{athm1.1}
Let $N\geq6$  and  $\xi_{i}^{*}$, $i=1,\cdots,k$,  be the non-degenerate critical points of $V(x)$.
Assume that  $V(x)\in C^2\big(B_{2\delta}(\xi^*_{i})\big)$ and $V(\xi_i^*)>0$. Let $u_\varepsilon^{(1)}$ and $u_\varepsilon^{(2)}$ be the $k$-peak solutions of \eqref{1.1} with the form \eqref{1a}.
Then there exists $\varepsilon_0>0$ such that for
$\varepsilon\in(0,\varepsilon_0)$, $u_\varepsilon^{(1)}\equiv u_\varepsilon^{(2)}$.
\end{theorem}
We next analyze the Morse index $m(u_\varepsilon)$ of the multi-peak solutions $u_\varepsilon$ to \eqref{1.1} of the form \eqref{1a}.
%Let us recall that the Morse index $m(u_\varepsilon)$ of the solution $u_\varepsilon$ of \eqref{1.1}
%is defined as the number of the negative eigenvalues $\lambda$'s of the following linearized eigenvalue problem
Let us recall that the Morse index $m(u_\varepsilon)$ of the solution $u_\varepsilon$ to \eqref{1.1} is defined as
the number of negative eigenvalues $\lambda$ of the following linearized eigenvalue problem
\begin{equation*}%\label{aw2.18}
\begin{cases}
-\Delta\omega+V(x)\omega
-(p-\varepsilon)u_\varepsilon^{p-1-\varepsilon}
\omega=\lambda\omega \;\;\text{in}\;\;\mathbb{R}^N,\\[1mm]
\omega\in H_V^1(\mathbb{R}^N),
\end{cases}
\end{equation*}
where $H_V^1(\mathbb{R}^N)$ is the closure of $C_0^\infty(\mathbb{R}^N)$ with respect to the norm $$\|u\|_{H_V^1(\mathbb{R}^N)}=\left(\displaystyle\int_{\mathbb{R}^N}\big(|\nabla u|^2+V(x)|u|^2\big)\mathrm{d}x
\right)^{\frac{1}{2}}.$$
It is straightforward to observe  that the Morse index of the solution $u_\varepsilon$ to \eqref{1.1}
can alternatively be defined as the number of eigenvalues  strictly less than $1$ corresponding to  the weighted eigenvalue problem
\begin{equation}\label{w2.18}
\begin{cases}
-\Delta\omega+V(x)\omega-(p-\varepsilon)\lambda u_\varepsilon^{p-1-\varepsilon}\omega=0
\;\;\;\text{in}\;\;\;\mathbb{R}^N,\\[1mm]
\omega\in H_V^1(\mathbb{R}^N).
\end{cases}
\end{equation}

 \vskip 0.1cm
Under the same assumptions as in Theorem \ref{athm1.1}, we can obtain the Morse index of the solution $u_\varepsilon$ to  \eqref {1.1} as follows.
\begin{theorem}\label{athm1.2}
Let $N\geq6$  and  $\xi_{i}^{*}$, $i=1,\cdots,k$,  be the non-degenerate critical points of $V(x)$.
Assume that  $V(x)\in C^2\big(B_{2\delta}(\xi^*_{i})\big)$ and $V(\xi_i^*)>0$.
Then  there exists  $\varepsilon_0>0$ such
that for $\varepsilon\in(0,\varepsilon_0)$, the Morse index of the solution $u_\varepsilon$ to \eqref{1.1} of the form \eqref{1a} is
$$m(u_{\varepsilon})=\displaystyle\sum_{i=1}^{k}m\big(\nabla^{2}V(\xi_{i}^{*})\big)+k,$$
where
$m\big(\nabla^{2}V(\xi_{i}^{*})\big)$ is the Morse index  of  the matrix $\big(\nabla^{2}V(\xi_{i}^{*})\big)$, i.e., the number of the negative eigenvalues of the matrix $\big(\nabla^{2}V(\xi_{i}^{*})\big)$.
\end{theorem}

As a direct consequence of Theorem \ref{athm1.2}, the non-degeneracy of the multi-peak solutions of \eqref {1.1} constructed in Theorem A follows immediately.
 We say that $u_\varepsilon$
 is non-degenerate if $w$ satisfies the following linearized equation of \eqref{1.1} at the $k$-peak solution $u_\varepsilon$, that is,
 \begin{equation*}%\label{aw2.19}
\begin{cases}
-\Delta w+V(x)w-(p-\varepsilon)u_\varepsilon^{p-1-\varepsilon} w=0 \;\;\text{in}\;\;\mathbb{R}^N,\\[1mm]
w\in H_V^1(\mathbb{R}^N),
\end{cases}
\end{equation*}
 then $w\equiv0$.
  \vskip 0.1cm
\begin{theorem}\label{athm1.3}
Under the same assumptions as in Theorem \ref{athm1.2},  there exists  $\varepsilon_0>0$ such
that for $\varepsilon\in(0,\varepsilon_0)$, the solution $u_\varepsilon$ obtained in Theorem A is non-degenerate.
\end{theorem}

\begin{remark}
The condition $N\geq6$ is technical. For $N=5$, we cannot estimate accurately the eigenvalues $\lambda_{\varepsilon,j}$ of the  weighted eigenvalue problem  \eqref{w2.18}  $\big( see \ Lemma \ \ref{lem4.p4}\big)$.
\end{remark}

Local uniqueness  and the Morse index of solutions of nonlinear elliptic problems are  important and challenging research topics
and it is impossible to give a complete bibliography here.
We merely mention some interesting results concerning peak solutions of critical elliptic problems associated with  \eqref{1.1}.
Local uniqueness of single peak solutions for the Brezis-Nirenberg problem was initially investigated in \cite{Gl} via a degree-counting method.
 Subsequently, the authors in \cite{DLY,GPY} employed  the local Pohozaev identities to study the local uniqueness problem of
multi-peak solutions in $\mathbb{R}^N$. It should be pointed out that the use of the local Pohozaev identities not only simplifies the
estimates, but also yields refined ``feedback" estimates for multi-peak solutions, thereby uncovering some essential hidden relationships.
 Moreover, the Morse index of  peak solutions in bounded domains  was computed in terms
of the negative eigenvalues of the Hessian matrix of the Robin function  (see \cite{BLR,Gp}).
Different from bounded domains, we will derive refined estimates for peak solutions by an appropriate projection,
 and then utilize the blow-up technique based on Pohozaev identities to compute the Morse index of multi-peak solutions.
 Finally we note that  the qualitative properties of peak solutions for critical elliptic equations
have garnered significant attention in recent literature, see for instance \cite{CLP,CPY,CeG,DLY,G05,GMPY,WY1} and references therein.
We are confident that the techniques proposed in this paper can provide valuable insights and make important contributions to the investigation of the qualitative characteristics of peak solutions for more general critical elliptic problems.

 \vskip 0.1cm
This paper is organized as follows.  In Section 2, we derive some precise estimates for multi-peak solutions that will be needed in subsequent sections.
 In Section 3,  we use the blow-up analysis  based on local Pohozaev identities to establish Theorem \ref{athm1.1}.
 Section 4 is dedicated to computing the Morse index of multi-peak solutions and thereby  presenting
 the proofs of  Theorems \ref{athm1.2} and  \ref{athm1.3}.
 %Some key estimates are included in the appendix.

%\textcolor{red}{
%Note $\omega_{l,\varepsilon}$ is the $l$-th eigenfunction to the $l$-th eigenvalue of problem \eqref{w2.18} and we suppose
%$$\lambda_{1,\varepsilon} <\lambda_{2,\varepsilon}\leq\cdots\leq
%\lambda_{l,\varepsilon}\cdots.  $$
%Moreover, the eigenvalues admit the variational characterization, we have
%$$\lambda_{l,\varepsilon}=\inf_{\dim S=l}\sup_{u\in S,u\neq0}\frac{\|u\|^{2}}
%{(p-\varepsilon)\displaystyle\int_{\mathbb{R}^N}
%u_{\varepsilon}^{p-1+\varepsilon}u^{2}\mathrm{d}x
%}, \ \ l\in\mathbb{N}.$$}

\section{Preliminaries}

In this section, we give several local Pohozaev identities and  some  key estimates  needed in later sections.

\vskip 0.1cm
%\begin{lemma}\label{lem: the computation}
%Let  $\alpha>1$. Then there holds
%\begin{align*}
%\int_{\mathbb{R}^{N}}\frac{1}
%{|x-y|^{N-2}}\frac{1}
%{(1+|y|^{2})^{\alpha}}\mathrm{d}y\leq
%\begin{cases}
%\frac{C}{(1+|x|)^{2(\alpha-1)}}
%,\ \ &\text{if}\;\;N-2\alpha>0,\\[2mm]
%\frac{C(1+\ln|x|)}{(1+|x|)^{N-2}}, \ \ &\text{if}\;\;N-2\alpha=0,
%\\[2mm]
%\frac{C}{(1+|x|)^{N-2}},\ \ &\text{if}\;\;N-2\alpha<0.
%\end{cases}
%\end{align*}
%\end{lemma}
%\begin{proof}[\bf Proof.]
%The proof of the lemma can be found in Lemma 2.3 in \cite{BF2023}.
%\end{proof}

Let
\begin{equation}\label{a2.1}
-\Delta u +V(x)u=u^{p-\varepsilon}\quad\text{in}
\quad\mathbb{R}^N,
\end{equation}
and
\begin{equation}\label{a2.2}
-\Delta\omega+V(x)\omega-(p-\varepsilon)\lambda u_\varepsilon^{p-1-\varepsilon}\omega=0\quad\text{in}\quad\mathbb{R}^N.
\end{equation}
To obtain the  refined estimates of peak solutions in Theorem A, we need several local Pohozaev identities for \eqref{a2.1} and \eqref{a2.2}.
Let $\Omega$ be a smooth bounded domain in $\mathbb{R}^N$. Then we have the following  Pohozaev-type identities.
\begin{lemma}\label{lem5.1}
Let  $u\in C^2(\mathbb{R}^N)$ be a solution of \eqref{a2.1}. Then
\begin{equation}\label{5.2}
\begin{split}
&\left(\frac{N}{p+1-\varepsilon}-\frac{N-2}{2}\right)
\int_{\Omega} u^{p+1-\varepsilon}\mathrm{d}x
-\frac{1}{2}\int_{\Omega}\big\langle x-\xi,\nabla V(x)\big\rangle u^{2}\mathrm{d}x
-\int_{\Omega}V(x)u^{2}\mathrm{d}x\\
=&\frac{1}{p+1-\varepsilon}\int_{\partial \Omega} u^{p+1-\varepsilon}\langle x-\xi,\nu\rangle\mathrm{d}s
+\int_{\partial \Omega}\frac{\partial u}{\partial\nu}\big\langle x-\xi,\nabla u\big\rangle\mathrm{d}s
-\frac{1}{2}\int_{\partial\Omega}|\nabla u|^2 \langle x-\xi,\nu\rangle\mathrm{d}s
\\
&-\frac{1}{2}\int_{\partial \Omega}V(x)u^2 \langle x-\xi,\nu\rangle\mathrm{d}s
+\frac{N-2}{2}\int_{\partial \Omega}\frac{\partial u}{\partial\nu}u\mathrm{d}s,
\end{split}
\end{equation}
and, for $l=1,\cdots,N$,
\begin{align}\label{14.3}
\frac{1}{2}\int_{\Omega}\frac{\partial V(x)}{\partial x_{l}}u^{2}\mathrm{d}x
=&\frac{1}{2}\int_{\partial\Omega}|\nabla u|^{2}\nu_{l}\mathrm{d}s
+\frac{1}{2}\int_{\partial\Omega}V(x)u^{2}
\nu_{l}\mathrm{d}s
\notag
\\&
-\int_{\partial \Omega}\frac{\partial u}{\partial x_{l}}\frac{\partial u}{\partial\nu}\mathrm{d}s
-\frac{1}{p+1-\varepsilon}\int_{\partial \Omega} u^{p+1-\varepsilon} \nu_{l}\mathrm{d}s,
\end{align}
where  $\xi\in\mathbb{R}^N$, %$\nu_{l}$ is the $l$-th component of the unit outer  normal vector of $\partial\Omega$, $l=1,\cdots,N$,
$\nu=(\nu_1,\cdots,\nu_N)$ is the outer  unit normal vector on $\partial\Omega$.
\end{lemma}
\begin{proof}[\bf Proof.]
Multiplying \eqref{a2.1} by $\langle x-\xi,\nabla u\rangle$ and integrating on $\Omega$, we have
\begin{equation*}
-\int_{\Omega}\Delta u\langle x-\xi,\nabla u\rangle \mathrm{d}x+\int_{\Omega}V(x)u\langle x-\xi,\nabla u\rangle \mathrm{d}x=\int_{\Omega} u^{p-\varepsilon}\langle x-\xi,\nabla u\rangle \mathrm{d}x.
\end{equation*}
By integrating by parts, we obtain
\begin{equation}\label{p4}
\begin{split}
-\int_{\Omega}\Delta u\langle x-\xi,\nabla u\rangle \mathrm{d}x
%=&-\int_{\partial\Omega}\frac{\partial u}{\partial \nu}\langle x-\xi,\nabla u\rangle dx+ \int_{\Omega}\left(|\nabla u|^{2}+\frac{1}{2}\langle x-\xi,\nabla (|\nabla u|^{2})\rangle\right)\mathrm{d}x\\
=-\int_{\partial\Omega}\frac{\partial u}{\partial \nu}\langle x-\xi,\nabla u\rangle \mathrm{d}s
+\frac{1}{2}\int_{\partial\Omega}\langle x-\xi,\nu\rangle|\nabla u|^{2}\mathrm{d}s-\frac{N-2}{2}\int_{\Omega}|\nabla u|^{2}\mathrm{d}x,
\end{split}
\end{equation}
\begin{equation}\label{p5}
\begin{split}
\int_{\Omega}V(x)u\langle x-\xi,\nabla u\rangle \mathrm{d}x
=&\frac{1}{2}\int_{\partial \Omega}V(x)u^{2}\langle x-\xi,\nu\rangle \mathrm{d}s-\frac{1}{2}\int_{\Omega}\langle x-\xi,\nabla V(x)\rangle u^{2}\mathrm{d}x -\frac{N}{2}\int_{\Omega}V(x)u^{2}\mathrm{d}x
\end{split}
\end{equation}
and
\begin{equation}\label{p6}
\int_{\Omega} u^{p-\varepsilon}\langle x-\xi,\nabla u\rangle \mathrm{d}x
=\frac{1}{p+1-\varepsilon}\int_{\partial \Omega}u^{p+1-\varepsilon}\langle x-\xi,\nu\rangle \mathrm{d}s-\frac{N}{p+1-\varepsilon}\int_{\Omega}u^{p+1-\varepsilon}\mathrm{d}x.
\end{equation}
On the other hand, we find
\begin{equation}\label{p7}
\begin{split}
\int_{\Omega}|\nabla u|^{2}\mathrm{d}x
%&=-\int_{\Omega}u\Delta u\mathrm{d}x+\int_{\partial \Omega}u\frac{\partial u}{\partial \nu}\mathrm{d}x\\
=\int_{\Omega}u^{p+1-\varepsilon}\mathrm{d}x-\int_{\Omega}V(x)u^{2}\mathrm{d}x+\int_{\partial \Omega}u\frac{\partial u}{\partial \nu}\mathrm{d}s.
\end{split}
\end{equation}
From \eqref{p4}, \eqref{p5}, \eqref{p6} and \eqref{p7}, we  readily  obtain \eqref{5.2}. The identity \eqref{14.3} can be derived by multiplying both sides of \eqref{a2.1} by $\frac{\partial u}{\partial x_l}$.
%\begin{equation*}
%\begin{split}
%&\left(\frac{N}{p+1-\varepsilon}-\frac{N-2}{2}\right)\int_{\Omega}u^{p+1-\varepsilon}dx
%-\frac{1}{2}\int_{\Omega}\langle x-\xi,\nabla V(x)\rangle u^{2}dx-\int_{\Omega}V(x)u^{2}dx\\
%&=\frac{1}{p+1-\varepsilon}\int_{\partial\Omega}u^{p+1-\varepsilon}\langle x-\xi,\nu\rangle dx+\int_{\partial\Omega} \frac{\partial u}{\partial \nu}\langle x-\xi,\nabla u\rangle dx-\frac{1}{2}\int_{\Omega}|\nabla u|^{2}\langle x-\xi,\nu\rangle dx\\
%&-\frac{1}{2}\int_{\partial\Omega}V(x)u^{2}\langle x-\xi,\nu\rangle dx+\frac{N-2}{2}\int_{\partial\Omega}\frac{\partial u}{\partial \nu}dx.
%\end{split}
%\end{equation*}
\end{proof}

\begin{lemma}\label{lem4.3}
Assume that $u$ and $\omega$ are classical solutions of \eqref{a2.1} and \eqref{a2.2} respectively. Then we have
\begin{align}\label{4.3}
&(p-\varepsilon)(\lambda-1)
\int_{\Omega} u^{p-1-\varepsilon}
\omega
\big\langle x-\xi,\nabla u\big\rangle
\mathrm{d}x
+\frac{N-2}{2}(p-\varepsilon)
(\lambda-1)\int_{\Omega} u^{p-\varepsilon}\omega\mathrm{d}x
\notag\\&\quad
-\frac{N-2}{2}\varepsilon\int_{\Omega} u^{p-\varepsilon}\omega\mathrm{d}x
+2\int_{\Omega}V(x)u\omega\mathrm{d}x
+\int_{\Omega}\big\langle x-\xi,\nabla V(x)\big\rangle u\omega\mathrm{d}x\notag
\\
=&-\int_{\partial \Omega}\frac{\partial u}{\partial\nu}\langle x-\xi,\nabla\omega\rangle\mathrm{d}s
-\int_{\partial \Omega}\frac{\partial \omega}{\partial\nu}\langle x-\xi,\nabla u\rangle\mathrm{d}s\
+\int_{\partial \Omega}\langle\nabla u,\nabla\omega\rangle\langle x-\xi,\nu\rangle\mathrm{d}s
\\
& -\frac{N-2}{2}\int_{\partial \Omega}\bigg(\frac{\partial u}{\partial\nu}\omega+\frac{\partial \omega}{\partial\nu}u\bigg)\mathrm{d}s
-\int_{\partial \Omega} u^{p-\varepsilon}\omega \langle x-\xi,\nu\rangle\mathrm{d}s
+\int_{\partial \Omega}V(x)u\omega\langle x-\xi,\nu\rangle\mathrm{d}s,\notag
\end{align}
and
\begin{align}\label{4.003}
&(p-\varepsilon)(\lambda-1)
\int_{\Omega} u^{p-1-\varepsilon}\omega \frac{\partial u}{\partial x_{l}}\mathrm{d}x+\int_{\Omega}\frac{\partial V(x)}{\partial x_{l}}u\omega\mathrm{d}x\\
=&\int_{\partial\Omega}\langle\nabla u,\nabla\omega\rangle\nu_{l}\mathrm{d}s
-\int_{\partial \Omega}\frac{\partial u}{\partial x_{l}}\frac{\partial\omega}{\partial\nu}
\mathrm{d}s
-\int_{\partial \Omega}\frac{\partial\omega}{\partial x_{l}}\frac{\partial u}{\partial\nu}
\mathrm{d}s
-\int_{\partial \Omega} u^{p-\varepsilon}\omega \nu_{l}\mathrm{d}s
+\int_{\partial \Omega}V(x)u\omega\nu_{l}\mathrm{d}s,\notag
\end{align}
where  $\xi\in\mathbb{R}^N$, %$\nu_{l}$ is the $l$-th component of the unit outer  normal vector of $\partial\Omega$, $l=1,\cdots,N$,
$\nu=(\nu_1,\cdots,\nu_N)$ is the outer unit  normal vector on $\partial\Omega$.
\end{lemma}
\begin{proof}[\bf Proof.]
Note that
\begin{align}\label{aa4.4}
&\int_{\Omega}\Big(-\Delta u\big\langle x-\xi,\nabla\omega
\big\rangle+(-\Delta) \omega\big\langle x-\xi,\nabla u\big\rangle\Big)\mathrm{d}x\notag\\
=&\int_{\Omega} u^{p-\varepsilon}\big\langle x-\xi,\nabla\omega\big\rangle\mathrm{d}x
-\int_{\Omega}V(x)u\big\langle x-\xi,\nabla\omega\big\rangle\mathrm{d}x+(p-\varepsilon)
\int_{\Omega} u^{p-1-\varepsilon}\omega
\big\langle x-\xi,\nabla u\big\rangle\mathrm{d}x\\
%&+(p-\varepsilon)
%\int_{\Omega} u^{p-1-\varepsilon}\omega
%\big\langle x-\xi,\nabla u\big\rangle\mathrm{d}x-
%\int_{\Omega}
%V(x)\omega\big\langle x-\xi,\nabla u\big\rangle\mathrm{d}x\\
&-\int_{\Omega}
V(x)\omega\big\langle x-\xi,\nabla u\big\rangle\mathrm{d}x+(p-\varepsilon)(\lambda-1)
\int_{\Omega} u^{p-1-\varepsilon}\omega
\big\langle x-\xi,\nabla u\big\rangle\mathrm{d}x.\notag
\end{align}
Then it follows from  direct calculations that
\begin{align*}
(p+1-\varepsilon)&\int_{\Omega} u^{p-\varepsilon}\omega\mathrm{d}x
+(p-\varepsilon)(\lambda-1)
\int_{\Omega} u^{p-\varepsilon}\omega\mathrm{d}x
-2\int_{\Omega}V(x)u\omega\mathrm{d}x\\
=&\int_{\Omega}\big(-u\Delta\omega+\omega(-\Delta u)\big)\mathrm{d}x\\
=&-\int_{\partial\Omega}u\frac{\partial\omega}{\partial\nu}
\mathrm{d}s-\int_{\partial \Omega}\omega\frac{\partial u}{\partial\nu}\mathrm{d}s
+2\int_{\Omega}\langle\nabla u,\nabla\omega\rangle\mathrm{d}x.
\end{align*}
Thus we obtain
\begin{align*}
\int_{\Omega}\big\langle\nabla u,\nabla\omega
\big\rangle\mathrm{d}x
=&\frac{p+1-\varepsilon}{2}
\int_{\Omega} u^{p-\varepsilon}\omega\mathrm{d}x
+\frac{(p-\varepsilon)(\lambda-1)}{2}
\int_{\Omega} u^{p-\varepsilon}\omega\mathrm{d}x\\
&-\int_{\Omega}   V(x)u\omega\mathrm{d}x+
\frac{1}{2}\int_{\partial \Omega}\Big(u
\frac{\partial\omega}{\partial\nu}
+\omega\frac{\partial u}{\partial\nu}\Big)\mathrm{d}s.
\end{align*}
By integrating by part, we find
\begin{align*}
\text{LHS \;of \;} \eqref{aa4.4}
=&-\int_{\partial\Omega}\frac{\partial u}{\partial\nu}\langle x-\xi,\nabla\omega\rangle\mathrm{d}s
+\int_{\Omega}\nabla u\nabla\langle x-\xi,\nabla\omega\rangle\mathrm{d}x\\
&-\int_{\partial \Omega}\frac{\partial \omega}{\partial\nu}\langle x-\xi,\nabla u\rangle\mathrm{d}s
+\int_{\Omega}\nabla \omega\nabla\langle x-\xi,\nabla u\rangle\mathrm{d}x\\
=&-\int_{\partial\Omega}\frac{\partial u}{\partial\nu}\langle x-\xi,\nabla\omega\rangle\mathrm{d}s
-\int_{\partial \Omega}\frac{\partial \omega}{\partial\nu}\langle x-\xi,\nabla u\rangle\mathrm{d}s\\
&+\int_{\partial \Omega}\langle\nabla u,\nabla\omega\rangle\langle x-\xi,\nu\rangle\mathrm{d}s
+(2-N)\int_{\Omega}\langle\nabla u,\nabla\omega\rangle\mathrm{d}x.
\end{align*}
On the other hand, we have
\begin{align*}
\text{RHS \;of \;} \eqref{aa4.4}=&\int_{\Omega}\big\langle x-\xi,\nabla(u^{p-\varepsilon}
\omega)\big\rangle\mathrm{d}x
-\int_{\Omega}V(x)\big\langle x-\xi,\nabla(u\omega)
\big\rangle\mathrm{d}x\\
&+(p-\varepsilon)(\lambda-1)
\int_{\Omega} u^{p-1-\varepsilon}\omega
\big\langle x-\xi,\nabla u\big\rangle\mathrm{d}x\\
=&\int_{\partial \Omega} u^{p-\varepsilon}\omega\langle x-\xi,\nu\rangle\mathrm{d}s
-N\int_{\Omega} u^{p-\varepsilon}\omega\mathrm{d}x
+N\int_{\Omega}V(x)u\omega\mathrm{d}x\\
&-\int_{\partial \Omega}V(x)u \omega\langle x-\xi,\nu\rangle\mathrm{d}s
+\int_{\Omega}\big\langle x-\xi,\nabla V(x)\big\rangle u\omega\mathrm{d}x
\\
&+(p-\varepsilon)(\lambda-1)
\int_{\Omega} u^{p-1-\varepsilon}\omega
\big\langle x-\xi,\nabla u\big\rangle\mathrm{d}x.
\end{align*}
Hence the desired identity \eqref{4.3} follows. Following the analogous calculations as above, we can also obtain \eqref{4.003}.
\end{proof}

Let  $u_{\varepsilon}$ be the $k$-peak  solutions of  \eqref{1.1} with the form \eqref{1a} obtained in Theorem A.  To derive refined estimates for $u_\varepsilon$,
we introduce  the following weighted-$L^\infty$ norms as in \cite{CLl1}
\begin{align}\label{aa2.12}
\|u\|_*=\sup\limits_{x\in\mathbb{R}^N}\bigg(
\displaystyle\sum_{i=1}^{k}
\frac{\mu_{\varepsilon,i}^{\frac{N-2}{2}}}
{\big(1+\mu_{\varepsilon,i}^{2}|x-
\xi_{\varepsilon,i}|^{2}\big)
^{\frac{N-2}{4}+\frac{\sigma}{2}}}\bigg)^{-1}
|u(x)|
\end{align}
and
\begin{align*}
\|h\|_{**}=\sup\limits_{x\in\mathbb{R}^N}
\bigg(
\displaystyle\sum_{i=1}^{k}
\frac{\mu_{\varepsilon,i}^{\frac{N+2}{2}}}
{\big(1+\mu_{\varepsilon,i}^{2}|x-
\xi_{\varepsilon,i}|^{2}\big)
^{\frac{N+2}{4}+\frac{\sigma}{2}}}\bigg)^{-1}
|h(x)|,
\end{align*}
where $\sigma\in \left(0,\frac{1}{2}\right)$, $\mu_{\varepsilon,i}\sim\varepsilon^{-\frac{1}{2}}$, $\xi_{\varepsilon,i}\in\mathbb{R}^N$, $i=1,\cdots,k$.
%\vskip 0.1cm
Then we denote by $C_*(\mathbb{R}^N)$ and $C_{**}(\mathbb{R}^N)$  the spaces of continuous functions  on
$\mathbb{R}^N$ with finite norms  $\|\cdot\|_*$ and $\|\cdot\|_{**}$ respectively.

From \cite{CLl1}, we know that $u_\varepsilon=\boldsymbol{W_{\mu_{\varepsilon},\xi_{\varepsilon}}}+\phi_\varepsilon$ satisfies
 the following nonlinear projected problem:
\begin{equation}\label{2.18}
\begin{cases}
-\Delta\big(\boldsymbol{W_{\mu_{\varepsilon},\xi_{\varepsilon}}}+\phi_\varepsilon\big)
+V(x)\big(\boldsymbol{W_{\mu_{\varepsilon},\xi_{\varepsilon}}}+\phi_\varepsilon\big)=
\big(\boldsymbol{W_{\mu_{\varepsilon},\xi_{\varepsilon}}}+\phi_\varepsilon\big)_+^{p-\varepsilon}
\,\,\,\text{in}\,\,\,{\mathbb{R}^N},\\[2mm]
\phi_\varepsilon\in C_{*}(\mathbb{R}^N)\cap H_V^1(\mathbb{R}^N),\\[2mm]
\displaystyle\int_{\mathbb{R}^N} W_{\mu_{\varepsilon,i},\xi_{\varepsilon,i}}^{p-1}
Z_{\mu_{\varepsilon,i},
\xi_{\varepsilon,i}}^{l}\phi_\varepsilon\mathrm{d}x=0,\;i=1,2,\cdots,k, \ \ l=0,1,\cdots,N,
\end{cases}
\end{equation}
where
$Z_{\mu_{\varepsilon,i},
\xi_{\varepsilon,i}}^{0}=\frac{\partial W_{\mu_{\varepsilon,i},
\xi_{\varepsilon,i}}}
{\partial \mu_{\varepsilon,i}}$ and $Z_{\mu_{\varepsilon,i},
\xi_{\varepsilon,i}}^{l}=\frac{\partial W_{\mu_{\varepsilon,i},\xi_{\varepsilon,i}}}
{\partial \xi_{\varepsilon,i,l}}$,  $i=1,\cdots,k$,  $l=1,\cdots,N$.
It is easy to verify that problem \eqref{2.18} can be rewritten as
\begin{align*}
-\Delta\phi_\varepsilon+V(x)\phi_\varepsilon
-p\boldsymbol{W}^{p-1}_{\boldsymbol{\mu_{\varepsilon},\xi_{\varepsilon}}}\phi_\varepsilon
=\mathcal{N}(\phi_\varepsilon)+\mathcal{R}
\,\,\,\,\text{in}\,\,\,\,{\mathbb{R}^N},\\
\end{align*}
where $\mathcal{N}(\phi_\varepsilon)=
\mathcal{N}_1(\phi_\varepsilon)
+\mathcal{N}_2(\phi_\varepsilon)$, $\mathcal{R}=\mathcal{R}_1+\mathcal{R}_2+\mathcal{R}_3+\mathcal{R}_4$,
\begin{align*}
\mathcal{N}_1(\phi_\varepsilon)&=\left(\boldsymbol{W_{\mu_{\varepsilon},\xi_{\varepsilon}}}
+\phi_\varepsilon\right)_+^{p-\varepsilon}
-\boldsymbol{W}^{p-\varepsilon}_{\boldsymbol{\mu_{\varepsilon},\xi_{\varepsilon}}}
-(p-\varepsilon)\boldsymbol{W}^{p-1-\varepsilon}_{\boldsymbol{\mu_{\varepsilon},\xi_{\varepsilon}}}\phi_\varepsilon,
\\
\mathcal{N}_2(\phi_\varepsilon)&=(p-\varepsilon)
\boldsymbol{W}^{p-1-\varepsilon}_{\boldsymbol{\mu_{\varepsilon},\xi_{\varepsilon}}}\phi_\varepsilon-p \boldsymbol{W}^{p-1-\varepsilon}_{\boldsymbol{\mu_{\varepsilon},\xi_{\varepsilon}}},
\end{align*}
and
\begin{align*}
%\mathcal{N}_1(\phi_\varepsilon)&=\left(\boldsymbol{W_{\mu_{\varepsilon},\xi_{\varepsilon}}}
%+\phi_\varepsilon\right)^{p-\varepsilon}
%-\boldsymbol{W}^{p-\varepsilon}_{\boldsymbol{\mu_{\varepsilon},\xi_{\varepsilon}}}
%-(p-\varepsilon)\boldsymbol{W}^{p-1-\varepsilon}_{\boldsymbol{\mu_{\varepsilon},\xi_{\varepsilon}}}\phi_\varepsilon,
%\\
%\mathcal{N}_2(\phi_\varepsilon)&=(p-\varepsilon)
%\boldsymbol{W}^{p-1-\varepsilon}_{\boldsymbol{\mu_{\varepsilon},\xi_{\varepsilon}}}\phi_\varepsilon-p \boldsymbol{W}^{p-1-\varepsilon}_{\boldsymbol{\mu_{\varepsilon},\xi_{\varepsilon}}},
%\\
\mathcal{R}_1&=\boldsymbol{W}^{p-\varepsilon}_{\boldsymbol{\mu_{\varepsilon},\xi_{\varepsilon}}}
-\displaystyle\sum_{i=1}^{k}\eta_{i} U_{\mu_{\varepsilon,i},\xi_{\varepsilon,i}}^{p},
\,\,\,\,\,\,\,\,\,\,\,\,\mathcal{R}_2=
2\displaystyle\sum_{i=1}^{k}\nabla\eta_{i} \nabla U_{\mu_{\varepsilon,i},\xi_{\varepsilon,i}},\\
\mathcal{R}_3&=\displaystyle\sum_{i=1}^{k}
\Delta\eta_{i} U_{\mu_{\varepsilon,i},\xi_{\varepsilon,i}},\,\,\,
\quad\quad\quad\quad\quad\,
\mathcal{R}_4=
-V(x)\displaystyle\boldsymbol{W}_{\boldsymbol{\mu_{\varepsilon},\xi_{\varepsilon}}}.
\end{align*}

To study the qualitative properties of peak solutions $u_\varepsilon$ to  \eqref{1.1} of the form \eqref{1a}, we need more sharper estimates
for $\|\phi_{\varepsilon}\|_{*}$ and $\frac{\partial V}{\partial x_{l}}(\xi_{\varepsilon,i})$  than those in \cite{CLl1}.
\begin{proposition}\label{prop2.4}
Let $N\geq6$ and $\xi_i^*$ $(i=1,2,\dots,k)$ be the non-degenerate critical points of $V(x)$.
 Assume that $V(x)\in C^2\big(B_{2\delta}(\xi^*_{i})\big)$ and $V(\xi_i^*)>0$.
  Then we have
\begin{equation*}\label{1q}
\|\phi_{\varepsilon}\|
_{*}=O(\varepsilon^{1-\frac{\sigma}{2}})  \ \ \text{and} \ \ \frac{\partial V}{\partial x_{l}}(\xi_{\varepsilon,i})=
o\big(\varepsilon^{\frac{1}{2}}\big), \ \ l=1,\cdots,N.
\end{equation*}
\end{proposition}

\begin{proof}
Note that
\begin{equation*}
\left|(a+b)_+^{p}-a^p-pa^{p-1}b\right|
\leq C\begin{cases}
\min\{a^{p-2}|b|^{2},|b|^{p}\},\;\;&\text{if}\;\;p\leq2,\\[2mm]
a^{p-2}|b|^{2}+|b|^{p},\;\;&\text{if}\;\;p>2,
\end{cases}
\end{equation*}
where $a>0, b\in\mathbb{R}$.
If $N\geq6$, by the H\"{o}lder inequality, we have
\begin{align*}
|\mathcal{N}_1(\phi_{\varepsilon})|\leq & C|\phi_{\varepsilon}|^{p-\varepsilon}\leq C\|\phi_{\varepsilon}\|_*^{p-\varepsilon}\left(\sum_{i=1}^k
\frac{\mu_{\varepsilon,i}^{\frac{N-2}{2}}}{(1+\mu_{\varepsilon,i}
|x-\xi_{\varepsilon,i}|)^{\frac{N-2}{2}+\sigma}}\right)^{p-\varepsilon}\\
\leq&C\|\phi_{\varepsilon}\|_*^{p-\varepsilon}\sum_{i=1}^k
\frac{\mu_{\varepsilon,i}^{\frac{N+2}{2}}}
{(1+\mu_{\varepsilon,i}|x-
\xi_{\varepsilon,i}|)^{\frac{N+2}{2}+\sigma}}.
\end{align*}
If $N=5$, we find
\[\begin{split}
|\mathcal{N}_1(\phi_{\varepsilon})|\leq & C\big(\boldsymbol{W}^{\frac{1}{3}-\varepsilon}_{\boldsymbol{\mu_\varepsilon,\xi_\varepsilon}}\phi_\varepsilon^2
+|\phi_\varepsilon|^{\frac{7}{3}-\varepsilon}\big)\\
\leq& C\boldsymbol{W}^{\frac{1}{3}-\epsilon}_{\boldsymbol{\mu_\varepsilon,\xi_\varepsilon}}\|\phi_\varepsilon\|_*^{2}\left(\sum_{j=1}^k
\frac{\mu_{\varepsilon,j}^{\frac{3}{2}}}{(1+\mu_{\varepsilon,j}|x-\xi_{\varepsilon,j}|)^{\frac{3}{2}+\sigma}}\right)^{2}\\
&+C\|\phi_\varepsilon\|_*^{\frac{7}{3}-\epsilon}
\left(\sum_{j=1}^k\frac{\mu_{\varepsilon,j}^{\frac{3}{2}}}{(1+\mu_{\varepsilon,j}|x-\xi_{\varepsilon,j}|)^{\frac{3}{2}+\sigma}}\right)^{\frac{7}{3}-\varepsilon}\\
\leq&C(\|\phi_\varepsilon\|_*^2+\|\phi_\varepsilon\|_*^{p-\varepsilon})\sum_{j=1}^k\frac{\mu_{\varepsilon,j}^{\frac{7}{2}}}{(1+\mu_{\varepsilon,j}|x-\xi_{\varepsilon,j}|)^{\frac{7}{2}+\sigma}}.
\end{split}
\]
Then it follows that
$$\|\mathcal{N}_1(\phi_\varepsilon)\|_{**}\leq C\big(\|\phi_\varepsilon\|_*^2+\|\phi_\varepsilon\|_*^{p-\varepsilon}\big).$$
By the mean value theorem, we infer
\begin{align*}
|\mathcal{N}_2(\phi_{\varepsilon})|\leq & \sum_{i=1}^k\varepsilon(\eta_i
U_{\mu_{\varepsilon,i},
\xi_{\varepsilon,i}})^
{p-1-\theta\varepsilon}\Big(1+\big|
\log(\eta_iU_{\mu_{\varepsilon,i},
\xi_{\varepsilon,i}})\big|
\Big)|\phi_{\varepsilon}|\\
%\leq&C\varepsilon\|\phi_{\varepsilon}\|_*
%\sum_{i=1}^k\big|\log\mu_{i,\varepsilon}\big|
%\frac{\mu_i^{2-\frac{N-2}{2}\theta
%\varepsilon}}
%{\big(1+\mu_{i,\varepsilon}^{2}
%|x-\xi_{i,\varepsilon}|^{2}\big)^{2-\frac{N-2}{2}\theta\varepsilon}}
%\sum_{i=1}^k\frac{\mu_{i,\varepsilon}^{\frac{N-2}{2}}}
%{\big(1+\mu_{i,\varepsilon}|x-\xi_{i,\varepsilon}|\big)^{\frac{N-2}{2}+\sigma}}\\
\leq&C\varepsilon|\log\varepsilon|\|\phi_{\varepsilon}\|_*\sum_{i=1}^k
\frac{\mu_{\varepsilon,i}^{\frac{N+2}{2}}}{\big(1
+\mu_{\varepsilon,i}|x-\xi_{\varepsilon,i}|\big)
^{\frac{N+2}{2}+\sigma}},
\end{align*}
where $\theta\in(0,1)$. Then we get
\begin{align*}\|\mathcal{N}_2(\phi_{\varepsilon})\|_{**}\leq C\varepsilon|\log\varepsilon|\|
\phi_{\varepsilon}\|_*.\end{align*}
A direct calculation yields that
\begin{align*}
|\mathcal{R}_1|\leq&\sum_{i=1}^k
\chi_{B_\delta(\xi_{\varepsilon,i})}
\big|U_{\mu_{\varepsilon,i},\xi_{\varepsilon,
i}}^{p-\varepsilon}-U_{\mu_{\varepsilon,i},
\xi_{\varepsilon,i}}^p\big|
+\sum_{i=1}^{k}\chi_{
B_{2\delta}(\xi_{\varepsilon,i})\backslash B_\delta(\xi_{\varepsilon,i})}
\big(U_{\mu_{\varepsilon,i},
\xi_{\varepsilon,i}}^{p-\varepsilon}
+U_{\mu_{\varepsilon,i},\xi_{\varepsilon,i}}^p\big)\\
\leq&C\varepsilon\sum_{i=1}^k
U_{\mu_{\varepsilon,i},
\xi_{\varepsilon,i}}^{p-\theta\varepsilon}\log U_{\mu_{\varepsilon,i},\xi_{\varepsilon,i}}
+C\sum_{i=1}^{k}\chi_{B_{2\delta}
(\xi_{\varepsilon,i})\backslash B_\delta(\xi_{\varepsilon,i})}
\frac{\mu_{\varepsilon,i}^{\frac{N+2}{2}}}
{\big(1+\mu_{i,\varepsilon}|x-\xi_{\varepsilon,i}|\big)^{N+2}}\\
\leq&C\varepsilon|\log\varepsilon|\sum_{i=1}^k
\frac{\mu_{\varepsilon,i}^{\frac{N+2}{2}}}
{\big(1+\mu_{\varepsilon,i}
|x-\xi_{\varepsilon,i}|\big)^{\frac{N+2}{2}+\sigma}}
+C\sum_{i=1}^k\mu_{\varepsilon,i}^{-\frac{N+2}{2}+\sigma}
\frac{\mu_{\varepsilon,i}^{\frac{N+2}{2}}}
{\big(1+\mu_{\varepsilon,i}|x-
\xi_{\varepsilon,i}|\big)^{\frac{N+2}{2}
+\sigma}},
\end{align*}
which gives that
\begin{align*}
\|\mathcal{R}_1\|_{**}\leq C\big(\varepsilon|\log\varepsilon|+
\varepsilon^{\frac{N+2}{4}-\frac{\sigma}{2}}\big).
\end{align*}
In a similar way, we obtain
\begin{align*}
|\mathcal{R}_2|\leq& C\sum_{i=1}^k\chi_{B_{2\delta}
(\xi_{\varepsilon,i})\backslash B_\delta(\xi_{\varepsilon,i})}\mu_{\varepsilon,i}^{\frac{N+2}{2}}
\frac{|x-\xi_{\varepsilon,i}|}{\big
(1+\mu_{\varepsilon,i}^2
|x-\xi_{\varepsilon,i}|^2\big)
^{\frac{N}{2}}}\\
\leq&C\sum_{i=1}^k\mu_{\varepsilon,i}^{-\frac{N-2}{2}+\sigma}
\frac{\mu_{\varepsilon,i}^{\frac{N+2}{2}}}
{\big(1+\mu_{\varepsilon,i}|x-\xi_{\varepsilon,i}|\big)
^{\frac{N+2}{2}+\sigma}},
\end{align*}
\begin{align*}
|\mathcal{R}_3|\leq& C\sum_{i=1}^k \chi_{B_{2\delta}(\xi_{\varepsilon,i})\backslash B_\delta(\xi_{\varepsilon,i})}\frac{\mu_{\varepsilon,i}^{\frac{N-2}{2}}}
{\big(1+\mu_{\varepsilon,i}|x-\xi_{\varepsilon,i}|\big)^{N-2}}\\
\leq&C\sum_{i=1}^k \mu_{\varepsilon,i}^{-\frac{N}{2}+1+\sigma}
\frac{\mu_{\varepsilon,i}^{\frac{N+2}{2}}}
{\big(1+\mu_{\varepsilon,i}|x-\xi_{\varepsilon,i}|\big)^{\frac{N+2}{2}+\sigma}},
%\leq&C\varepsilon^{\frac{N-2}{4}-\frac{\sigma}{2}}
%\sum_{i=1}^k\frac{\mu_{\varepsilon,i}^{\frac{N+2}{2}}}
%{\big(1+\mu_{\varepsilon,i}|x-
%\xi_{\varepsilon,i}|\big)^{\frac{N+2}{2}
%+\sigma}}
\end{align*}
and
\begin{align*}
|\mathcal{R}_4|&\leq C
\displaystyle\sum_{i=1}^{k}\chi_{B_{2\delta}
(\xi_{\varepsilon,i})}
\frac{\mu_{\varepsilon,i}^{\frac{N+2}{2}}}
{\big(1+\mu_{\varepsilon,i}|x-\xi_{\varepsilon,i}|\big)
^{\frac{N+2}{2}+\sigma}}
\frac{\mu_{\varepsilon,i}^{-2}}
{\big(1+\mu_{\varepsilon,i}|x-\xi_{\varepsilon,i}|\big)
^{\frac{N}{2}-3-\sigma}}\\[0.02mm]&
\leq C
\displaystyle\sum_{i=1}^{k}\chi_{B_{2\delta}
(\xi_{\varepsilon,i})}
\mu_{\varepsilon,i}^{-2+\sigma}
\frac{\mu_{\varepsilon,i}^{\frac{N+2}{2}}}
{\big(1+\mu_{\varepsilon,i}|x-\xi_{\varepsilon,i}|\big)
^{\frac{N+2}{2}+\sigma}}.
\end{align*}
 By Proposition 2.3 in \cite{CLl1}, we conclude
\begin{align*}
\|\phi_{\varepsilon}\|_{*}\leq C\big\|\mathcal{N}(\phi_\varepsilon)+\mathcal{R}
\big\|_{**}
\leq C\varepsilon^{1-\frac{\sigma}{2}} .
\end{align*}
Next, we estimate $\frac{\partial V}{\partial x_{l}}(\xi_{\varepsilon,i})$.
Choosing $\Omega=B_\rho(\xi_{\varepsilon,i})$ in Lemma \ref{lem5.1},  we have
\begin{align}\label{676}
\int_{B_\rho(\xi_{\varepsilon,i})}\frac{\partial V(x)}
{\partial x_{l}} u_\varepsilon^{2}&=
\int_{\partial B_\rho(\xi_{\varepsilon,i})}V(x)u_\varepsilon^{2}
\nu_{l}
\mathrm{d}s-\frac{2}{p+1-\varepsilon}
\int_{\partial B_\rho(\xi_{\varepsilon,i})}
u_\varepsilon^{p+1-\varepsilon}\nu_{l}
\mathrm{d}s\notag\\
&\quad+\int_{\partial B_\rho(\xi_{\varepsilon,i})}\big|\nabla u_\varepsilon\big|^{2}\nu_{l}\mathrm{d}s-2
\int_{\partial B_\rho(\xi_{\varepsilon,i})}\frac{\partial u_\varepsilon}{\partial x_{l}}\frac{\partial u_\varepsilon}{\partial\nu}\mathrm{d}s\notag
\\
&=O\bigg(\int_{\partial B_\rho(\xi_{\varepsilon,i})}\Big(|\nabla\phi_{\varepsilon}|^{2}
+|\phi_{\varepsilon}|^{2}+|\phi_{\varepsilon}|
^{p+1-\varepsilon}\Big)\mathrm{d}s.
\end{align}
Define $S_{1}=B_{5\delta}(\xi_{\varepsilon,i})\backslash B_{2\delta}(\xi_{\varepsilon,i})$ and $S_{2}=B_{4\delta}(\xi_{\varepsilon,i})
\backslash B_{3\delta}(\xi_{\varepsilon,i})$. Let $\varphi$ be a smooth cut-off function satisfying
$\varphi(x)=1$, if $x\in S_{2}$; $\varphi(x)=0$ if $x\in\mathbb{R}^N\backslash S_{1}$ and $\big|\nabla \varphi(x)\big|\leq\frac{2}{\delta}$.
Since
\begin{equation*}
-\Delta\phi_\varepsilon
=\Delta\boldsymbol{W_{\mu_{\varepsilon},\xi_{\varepsilon}}}
-V(x)\boldsymbol{W_{\mu_{\varepsilon},\xi_{\varepsilon}}}+
\big(\boldsymbol{W_{\mu_{\varepsilon}
,\xi_{\varepsilon}}}+\phi_\varepsilon
\big)_+^{p-\varepsilon}-V(x)\phi_\varepsilon
\,\,\,\text{in}\,\,\,{\mathbb{R}^N}.
\end{equation*}
Multiplying the above equation by $\varphi^{2}\phi_\varepsilon$ and integrating by parts, we have
%\[
%\begin{split}
%\int_{S_{1}} \varphi^{2} |\nabla\phi_\varepsilon|^2\mathrm{d}x&=
%\int_{S_{1}}\varphi^{2}
%\phi_\varepsilon^{p+1-\varepsilon}\mathrm{d}x-\int_{S_{1}}
%V(x)\phi_\varepsilon^{2}\varphi^{2}\mathrm{d}x-2
%\int_{S_{1}}\varphi\phi_\varepsilon
%\nabla\phi_\varepsilon\nabla\varphi\mathrm{d}x
%\\
%&
%\leq\int_{S_{1}}|\phi_\varepsilon|^{p
%+1-\varepsilon}\mathrm{d}x+\int_{S_{1}}|
%\phi_\varepsilon|^{2}\mathrm{d}x+\frac{1}{2}
%\int_{S_{1}} \varphi^{2} |\nabla\phi_\varepsilon|^2\mathrm{d}x
%+2\int_{S_{1}}|
%\phi_\varepsilon|^{2}
%|\nabla\varphi|^2\mathrm{d}x
%,
%\end{split}
%\]
%which shows that
\begin{align*}
\int_{S_{2}}  |\nabla\phi_\varepsilon|^2\mathrm{d}x
%\int_{S_{1}} \varphi^{2} |\nabla\phi_\varepsilon|^2\mathrm{d}x
\leq& C\int_{S_{1}}|\phi_\varepsilon|^{p
+1-\varepsilon}\mathrm{d}x+C\int_{S_{1}}|
\phi_\varepsilon|^{2}\mathrm{d}x
\\
\leq &C\mu_{\varepsilon,i}^{-\frac{2N\sigma}{N-2}}
\|\phi_\varepsilon
\|_*^{p+1-\varepsilon}
+C\mu_{\varepsilon,i}^{-2\sigma}
\|\phi_\varepsilon\|_*^2
\leq
C\varepsilon^{2}.
\end{align*}
Then there exists  $\rho\in(3\delta,4\delta)$ such that
\begin{align*}
\int_{\partial B_\rho(\xi_{\varepsilon,i})}\Big(|\nabla\phi_{\varepsilon}|^{2}
+|\phi_{\varepsilon}|^{2}+|\phi_{\varepsilon}|
^{p+1-\varepsilon}\Big)\mathrm{d}s\leq C\varepsilon^{2}.
\end{align*}
Taking advantage of \eqref{676}, we derive
\begin{align}\label{677}
\int_{B_\rho(\xi_{\varepsilon,i})}\frac{\partial V(x)}
{\partial x_{l}} u_\varepsilon^{2}\mathrm{d}x=O\big(\varepsilon^{2}\big).
\end{align}
On the other hand,  we have
\begin{align}\label{678}
\int_{B_\rho(\xi_{\varepsilon,i})}\frac{\partial V(x)}
{\partial x_{l}} u_\varepsilon^{2}\mathrm{d}x
%&=\int_{B_\rho(\xi_{i,\varepsilon})}\frac{\partial V(x)}
%{\partial x_{l}} \eta_{i}^{2}U_{\mu_{i,\varepsilon},
%\xi_{i,\varepsilon}}^{2}\mathrm{d}x
%+O\bigg(\int_{B_\rho(\xi_{i,\varepsilon})}
%U_{\mu_{i,\varepsilon},\xi_{i,\varepsilon}}
%|\phi_{\varepsilon}|+|\phi_{\varepsilon}|^{2}\mathrm{d}x
%\bigg)\notag
%\\[0.02mm]
&=\int_{B_\delta(\xi_{\varepsilon,i})}\frac{\partial V(x)}
{\partial x_{l}}U_{\mu_{\varepsilon,i},\xi_{\varepsilon,i}
}^{2}\mathrm{d}x
+O\big(\mu_{\varepsilon,i}^{2-N}\big)
+O\big(\mu_{\varepsilon,i}^{-2}
\|\phi_{\varepsilon}\|_{*}\big)
+O\big(\mu_{\varepsilon,i}^{-2\sigma}
\|\phi_{\varepsilon}\|_{*}^{2}\big)\notag
\\[0.02mm]
&=\int_{B_{\varepsilon^{\tau_{0}}}
(\xi_{\varepsilon,i})}\frac{\partial V(x)}
{\partial x_{l}}U_{\mu_{\varepsilon,i},\xi_{\varepsilon,i}
}^{2}\mathrm{d}x
%%+
%%\int_{B_{\delta}(\xi_{\varepsilon,i})\backslash
%%B_{\varepsilon^{\tau_{0}}}
%%(\xi_{\varepsilon,i})}\frac{\partial V(x)}
%%{\partial x_{l}}U_{\mu_{\varepsilon,i},\xi_{\varepsilon,i}
%%}^{2}\mathrm{d}x
+
O(\varepsilon^{2-\frac{\sigma}{2}})\notag
\\[0.02mm]
%&=\int_{B_{\varepsilon^{\tau_{0}}}
%(\xi_{i,\varepsilon})}
%\Big(
%\frac{\partial V(\xi_{i,\varepsilon})}
%{\partial x_{l}}+o\big(|x-\xi_{i,\varepsilon}|\big)\Big)
%U_{\mu_{i,\varepsilon},\xi_{i,\varepsilon}
%}^{2}\mathrm{d}x+
%O(\varepsilon^{2-\frac{\sigma}{2}})\notag
%\\[0.02mm]
&=A \mu_{\varepsilon,i}^{-2}
\frac{\partial V(\xi_{\varepsilon,i})}
{\partial x_{l}}+
o\big(\varepsilon^{-\frac{3}{2}}\big),
\end{align}
where $\tau_{0}>0$ is a small constant, $A=\int_{\mathbb{R}^N}U_{1,0}^2\mathrm{d}x$.
Combining \eqref{677} with \eqref{678}, we deduce that $\frac{\partial V(\xi_{\varepsilon,i})}{\partial x_{l}}
=o\big(\varepsilon^{\frac{1}{2}}\big)$ and the proof of  Proposition \ref{prop2.4} is complete.
\end{proof}

\vskip 0.1cm
Let $\omega_\varepsilon$ be a solution of problem \eqref{w2.18}.  We next  establish  the decay estimates of the $k$-peak solutions $u_\varepsilon$ and
$\omega_\varepsilon$, which is crucial for applying the local Pohozaev-type identities.
Without loss of generality, we assume that $\|\omega_\varepsilon\|_*=1$. For this purpose, we need the following well-known
integral estimate.

\begin{lemma}\label{lem: the computation}
Let  $\alpha>1$. Then there holds
\begin{align*}
\int_{\mathbb{R}^{N}}\frac{1}
{|x-y|^{N-2}}\frac{1}
{(1+|y|^{2})^{\alpha}}\mathrm{d}y\leq
\begin{cases}
\frac{C}{(1+|x|)^{2(\alpha-1)}}
,\ \ &\text{if}\;\;N-2\alpha>0,\\[2mm]
\frac{C(1+\log|x|)}{(1+|x|)^{N-2}}, \ \ &\text{if}\;\;N-2\alpha=0,
\\[2mm]
\frac{C}{(1+|x|)^{N-2}},\ \ &\text{if}\;\;N-2\alpha<0.
\end{cases}
\end{align*}
\end{lemma}
\begin{proof}[\bf Proof.]
The proof of the lemma can be found in Lemma 2.3 in \cite{BF2023}.
\end{proof}

\begin{lemma}\label{lem4.1}
There exists a positive constant $C$ independent of $\varepsilon$ such that
\begin{align*}
|u_\varepsilon(x)|\leq C\displaystyle\sum_{i=1}^{k}
\frac{\mu_{\varepsilon,i}^{\frac{N-2}{2}}}
{\left(1+\mu_{\varepsilon,i}^{2}|x-
\xi_{\varepsilon,i}|^{2}\right)^{\frac{N-2}{2}}},\;\;
|\nabla u_\varepsilon(x)|\leq C\displaystyle\sum_{i=1}^{k}\frac{\mu_{\varepsilon,i}^{\frac{N}{2}}}
{\left(1+\mu_{\varepsilon,i}^{2}|x-
\xi_{\varepsilon,i}|^{2}\right)^{\frac{N-1}{2}}},\\
|\omega_\varepsilon(x)|\leq C\displaystyle\sum_{i=1}^{k}\frac{
\mu_{\varepsilon,i}^{\frac{N-2}{2}}}
{\left(1+\mu_{\varepsilon,i}^{2}|x-\xi_{\varepsilon,i}|^{2}\right)^{\frac{N-2}{2}}},\;\;
|\nabla \omega_\varepsilon(x)|\leq C\displaystyle\sum_{i=1}^{k}\frac{\mu_{\varepsilon,i}^{\frac{N}{2}}}
{\left(1+\mu_{\varepsilon,i}^{2}|x-\xi_{\varepsilon,i}|^{2}\right)^{\frac{N-1}{2}}}.\\
\end{align*}
\begin{proof}[\bf Proof.]% We only give the estimates for $u_\varepsilon$.
Note that
\begin{align*}
|u_\varepsilon(x)|&\leq C\displaystyle\sum_{i=1}^{k}
\frac{\mu_{\varepsilon,i}^{\frac{N-2}{2}}}
{\big(1+\mu_{\varepsilon,i}^{2}|x-\xi_{\varepsilon,i}|^{2}\big)
^{\frac{N-2}{4}+\frac{\sigma}{2}}}.
\end{align*}
By the Green's representation formula, we see that
\[
u_\varepsilon(x)=
\int_{\mathbb{R}^N}G(x,y) u_\varepsilon^{p-\varepsilon}(y)\mathrm{d}y,
\]
where $G(x,y)$ is the Green's function of the operator $-\Delta+V(x)$ in $\mathbb{R}^N$.
Using  Lemma  \ref{lem: the computation}, we deduce
\begin{align*}
|u_\varepsilon(x)|&\leq C\int_{\mathbb{R}^N}\frac{1}{|x-y|^{N-2}} |u_\varepsilon|^{p-\varepsilon}\mathrm{d}y\\
&\leq C\sum_{i=1}^k\int_{\mathbb{R}^N}\frac{1}{|x-y|^{N-2}} \frac{\mu_{\varepsilon,i}^{\frac{N+2}{2}}}
{(1+\mu_{\varepsilon,i}^{2}|y-
\xi_{\varepsilon,i}|^{2})^{\frac{N+2+2p\sigma}{4}-\frac{N-2+2\sigma}{4}\varepsilon}}\mathrm{d}y\\
&\leq C\sum_{i=1}^k\int_{\mathbb{R}^N}\frac{1}{
|y-\mu_{\varepsilon,i}(x-\xi_{\varepsilon,i})|^{N-2}} \frac{\mu_{\varepsilon,i}^{\frac{N-2}{2}}}{(1+|y|^{2})^{\frac{N+2+2p\sigma}{4}-\frac{N-2+2\sigma}{4}\varepsilon}}\mathrm{d}y\\
&\leq C\sum_{i=1}^k \frac{\mu_{\varepsilon,i}^{\frac{N-2}{2}}}
{(1+\mu_{\varepsilon,i}^{2}
|x-\xi_{\varepsilon,i}|^{2})^{
\frac{N-2+2p\sigma}{4}-\frac{N-2+2\sigma}{4}
\varepsilon}}.
%\\|\nabla u_\varepsilon(x)|&=
%\Big|\int_{\mathbb{R}^N}\nabla_{x}G(x,y) u_\varepsilon^{p-\varepsilon}(y)\mathrm{d}y
%\Big|\leq C\int_{\mathbb{R}^N}\frac{1}{|x-y|^{N-1}} \big|u_\varepsilon\big|^{p-\varepsilon}
%\mathrm{d}y
\end{align*}
Since
$$
\frac{N-2+2p\sigma}{4}-\frac{N-2+2\sigma}{4}
\varepsilon>\frac{N-2}{4}+\frac{\sigma}{2}\;\;\;\text{for}\;\;\varepsilon\;\;\text{small}.
$$
Thus we can continue the above procedure to conclude  the  estimate for $u_\varepsilon$.

On the other hand, we have
\begin{align*}
|\nabla u_\varepsilon(x)|
%&=\Big|\int_{\mathbb{R}^N}
%\nabla G(x,y) u_\varepsilon^{p-\varepsilon}(y)\mathrm{d}y
%\Big|
\leq C\int_{\mathbb{R}^N}\frac{1}{|x-y|^{N-1}} |u_\varepsilon|
^{p-\varepsilon}\mathrm{d}y.\end{align*}
Using the similar arguments as above, we  derive  the desired estimate for $\nabla u_\varepsilon$.
The estimates for $\omega_\varepsilon$ can be obtained in a similar way and the proof of Lemma \ref{lem4.1} is complete.
\end{proof}
\end{lemma}

\section{Local uniqueness of the peak solutions: Proof of Theorem \ref{athm1.1}}
In this section, we investigate the local uniqueness of $k$-peak solutions and present the proof of  Theorem \ref{athm1.1}.
For this purpose, we assume that $u_\varepsilon^{(1)}$ and $u_\varepsilon^{(2)}$ are two distinct
$k$-peak solutions  to problem \eqref{1.1} with the form \eqref{1a}, i.e.,
\[
u_\varepsilon^{(m)}=\boldsymbol{W}_{{\boldsymbol{\mu}_\varepsilon^{(m)}},
{\boldsymbol{\xi}^{(m)}_\varepsilon}}+\phi_\varepsilon^{(m)},\;\;m=1,2,
\]
where $\boldsymbol{\mu}^{(m)}_\varepsilon=
(\mu_{\varepsilon,1}^{(m)},\cdots,\mu_{\varepsilon,k}^{(m)})$,
$\boldsymbol{\xi}^{(m)}_\varepsilon=
(\xi_{\varepsilon,1}^{(m)},\cdots,\xi_{\varepsilon,k}^{(m)})$,
$\mu_{\varepsilon,i}^{(m)}\sim\varepsilon^{-\frac{1}{2}}$,
$\xi_{\varepsilon,i}^{(m)}\rightarrow\xi_{i}^*$,   $i=1,\cdots,k$   as $\varepsilon\rightarrow0$ and  $\|\phi_\varepsilon^{(m)}\|_{*}=O(\varepsilon^{1-\frac{\sigma}{2}})$.

Define \[
\psi_\varepsilon=\frac{u_\varepsilon^{(1)}-u_\varepsilon^{(2)}}{\big\|u_\varepsilon^{(1)}-u_\varepsilon^{(2)}\big\|_*}.
\]
Then $\psi_\varepsilon$ satisfies $\|\psi_\varepsilon\|_*=1$ and
\begin{equation*}
-\Delta\psi_\varepsilon+V(x)\psi_\varepsilon
=(p-\varepsilon) a_\varepsilon(x)\psi_\varepsilon,
\end{equation*}
where
$$
a_\varepsilon(x)=\int_0^1\left(tu_\varepsilon^{(1)}+(1-t)u_\varepsilon^{(2)}\right)
^{p-1-\varepsilon}\mathrm{d}t.
$$

Firstly, we  utilize the  Pohozaev  identity  \eqref{5.2} for $u_\varepsilon^{(m)}$  on $B_{\delta}(\xi^{(1)}_{\varepsilon,i})$ to obtain
the refined estimates of the concentration parameters $\mu_{\varepsilon,i}^{(m)}$, $m=1,2$, $i=1,2,\cdots,k$.
\begin{lemma}\label{lem5.3}
Let $N\geq6$ and $\xi_i^*$ $(i=1,2,\dots,k)$ be the non-degenerate critical points of $V(x)$.
 Assume that $V(x)\in C^2\big(B_{2\delta}(\xi^*_{i})\big)$ and $V(\xi_i^*)>0$. Then there holds
\begin{align*}
\mu_{\varepsilon,i}^{(m)}=\Lambda\varepsilon^{-\frac{1}{2}}
+O\left(\varepsilon^{\frac{1-\sigma}{2}}
\right), \ \ m=1,2,\; i=1,\cdots,k,
\end{align*}
where  $\Lambda$ is a positive constant independent of $\varepsilon$.
\end{lemma}
\begin{proof}[\bf Proof.]For simplicity, we omit the indices  $(m)$. By Proposition \ref{prop2.4}, we have
 \begin{align*}
 \int_{B_\delta(\xi_{\varepsilon,i})} u_\varepsilon^{p+1-\varepsilon}\mathrm{d}x
 =&\int_{B_\delta(\xi_{\varepsilon,i})} U_{\mu_{\varepsilon,i},\xi_{\varepsilon,i}}^{p+1-\varepsilon}\mathrm{d}x
 +O\bigg(\int_{B_\delta(\xi_{\varepsilon,i})} U_{\mu_{\varepsilon,i},\xi_{\varepsilon,i}}^{p-\varepsilon}|\phi_\varepsilon|\mathrm{d}x
 +\int_{B_\delta(\xi_{\varepsilon,i})} |\phi_\varepsilon|^{p+1-\varepsilon}\mathrm{d}x
 \bigg)\\
 =&\int_{B_\delta(\xi_{\varepsilon,i})} U_{\mu_{\varepsilon,i},\xi_{\varepsilon,i}}^{p+1}\mathrm{d}x
 +O(\varepsilon|\log\varepsilon|)+O\big(\|\phi_\varepsilon\|_*\big)
 +O\big(\|\phi_\varepsilon\|_*
 ^{p+1-\varepsilon}\big)\\
 =&B+O\left(\varepsilon\big|\log\varepsilon\big|\right)
 +O\left(\varepsilon^{1-\frac{\sigma}{2}}\right),
 \end{align*}
 where $B=\int_{\mathbb{R}^N}U_{1,0}^{p+1}.$
Observe that
 $$
 \frac{1}{\mu_{\varepsilon,i}}\leq\frac{C}
 {(1+\mu_{\varepsilon,i}^{2}|x-
 \xi_{\varepsilon,i}|^{2})^{
 \frac{1}{2}}}\;\;
 \text{in}\;\;B_\delta(\xi_{\varepsilon,i})
 ,\;\;i=1,\cdots,k.
 $$
 Then we deduce
 \begin{align*}
 \int_{B_\delta(\xi_{\varepsilon,i})}V(x)u_\varepsilon^2\mathrm{d}x
 =&
 \int_{B_\delta(\xi_{\varepsilon,i})}V(x)
 U_{\mu_{\varepsilon,i},\xi_{\varepsilon,i}}^2\mathrm{d}x
 +O\left(\mu_{\varepsilon,i}^{-2}\|\phi_\varepsilon\|_*\right)
 +O\left(\mu_{\varepsilon,i}^{-2\sigma}
 \|\phi_\varepsilon\|_*^2\right)
 \\
 =&V(\xi_{\varepsilon,i})\int_{B_\delta(\xi_{\varepsilon,i})}
 U_{\mu_{\varepsilon,i},\xi_{\varepsilon,i}}^2\mathrm{d}x+O(\varepsilon^2|\log\varepsilon|)
 +O\left(\mu_{\varepsilon,i}^{-2}\|\phi_\varepsilon\|_*\right)
 +O\left(\mu_{\varepsilon,i}^{-2\sigma}
 \|\phi_\varepsilon\|_*^2\right)\\
 =&A\mu_{\varepsilon,i}^{-2}V(\xi_{\varepsilon,i})
 +O\left(\varepsilon^{2-\frac{\sigma}{2}}\right),
 \end{align*}
 and
 \begin{align*}
&\int_{B_\delta(\xi_{\varepsilon,i})}\big\langle x-\xi_{\varepsilon,i},\nabla V(x)\big\rangle u_\varepsilon^2\mathrm{d}x\\
%=
%\int_{B_\delta(\xi_{\varepsilon,i})}\big\langle x-\xi_{\varepsilon,i},\nabla V(x)\big\rangle\big(
% U_{\mu_{\varepsilon,i},\xi_{\varepsilon,i}}+\phi_{\varepsilon}\big)
% ^{2}\mathrm{d}x
%\\[0.02mm]
 =&
\int_{B_\delta(\xi_{\varepsilon,i})}\big\langle x-\xi_{\varepsilon,i},\nabla V(x)\big\rangle U_{\mu_{\varepsilon,i},\xi_{\varepsilon,i}}^2\mathrm{d}x
 +O\left(\mu_{\varepsilon,i}^{-2}\|\phi_\varepsilon\|_*\right)
 +O\left(\mu_{\varepsilon,i}^{-2\sigma}
 \|\phi_\varepsilon\|_*^2\right)\\[0.02mm]
 =&O\bigg(\int_{B_\delta(\xi_{\varepsilon,i})} \big|x-\xi_{\varepsilon,i}\big |^2 U_{\mu_{\varepsilon,i},\xi_{\varepsilon,i}}^2\mathrm{d}x\bigg)
 +O\left(\varepsilon^{2-\frac{\sigma}{2}}\right)
 \\[0.02mm]
 =&O\left(\varepsilon^{2-\frac{\sigma}{2}}\right).
 \end{align*}
It follows from Lemma \ref{lem4.1} that
\begin{align*}
\text{RHS\;of}\;\eqref{5.2}=O\left(\mu_{\varepsilon,i}^{2-N}\right).
\end{align*}
Hence, by virtue of \eqref{5.2}, we obtain
$$
\mu_{\varepsilon,i}=\Lambda\varepsilon^{-\frac{1}{2}}
+O\Big(\varepsilon^{\frac{1-\sigma}{2}}\Big).
$$
\end{proof}
As a direct consequence of Lemma \ref{lem5.3}, we can compare $u_\varepsilon^{(1)}$ and $u_\varepsilon^{(2)}$.
\begin{lemma}\label{cor5.4}
Let $N\geq6$ and $\xi_i^*$ $(i=1,2,\dots,k)$ be the non-degenerate critical points of $V(x)$.
 Assume that $V(x)\in C^2\big(B_{2\delta}(\xi^*_{i})\big)$ and $V(\xi_i^*)>0$.
 Then we have
\begin{align*}
  \left|u_\varepsilon^{(1)}-
  u_\varepsilon^{(2)}\right|&=
  o(1)\displaystyle\sum_{i=1}^{k}
  \eta^{(1)}_{i}(x) U_{\mu_{\varepsilon,i}^{(1)},
  \xi_{\varepsilon,i}^{(1)}
  }+O\bigg(\displaystyle\sum_{i=1}^{k}
\chi_{B_{\frac{5\delta}{2}}
\left(\xi_{\varepsilon,
i}^{(1)}\right)\backslash B_{\frac{\delta}{2}}\left(\xi_{\varepsilon,
i}^{(1)}\right)}
\left(\mu_{\varepsilon,i}^{(1)}
\right)^{\frac{2-N}{2}}
  \bigg)
  +O\big(
  |\phi^{(1)}_{\varepsilon}|
  +|\phi^{(2)}_{\varepsilon}|\big),
\end{align*}
where $\eta^{(m)}_{i}(x)=\eta(x-\xi_{\varepsilon,i}^{(m)})$, $m=1,2$.
\end{lemma}
\begin{proof}[\bf Proof.]
Since $\xi_{i}^{*}$  are the non-degenerate critical points of $V(x)$, then  we have
\begin{align*}
\big|
\xi^{(m)}_{\varepsilon,i}- \xi_{i}^{*}\big|\leq\gamma\left|\nabla V(\xi^{(m)}_{\varepsilon,i})-\nabla V(\xi_{i}^{*})\right|,\;\;m=1,2,\;i=1,\cdots,k,
\end{align*}
where $\gamma>0$ is a  constant independent of $\varepsilon$.
Then, by Proposition \ref{prop2.4}, we obtain
$$\left|
\xi^{(1)}_{\varepsilon,i}- \xi^{(2)}_{\varepsilon,i}\right|
=o\big(\varepsilon^{\frac{1}{2}}\big),\ \ i=1,2,\cdots,k.$$
Hence we derive
\begin{align*}
  \left|u_\varepsilon^{(1)}-u_\varepsilon^{(2)}
  \right|
&\leq\displaystyle\sum_{i=1}^{k}\eta^{(1)}_{i}(x) \left|U_{\mu_{\varepsilon,i}^{(1)},
  \xi_{\varepsilon,i}^{(1)}
  }-U_{\mu_{\varepsilon,i}^{(2)},
  \xi_{\varepsilon,i}^{(2)}
  }
  \right|+
  \displaystyle\sum_{i=1}^{k}\left|
  \eta^{(1)}_{i}(x)-
  \eta^{(2)}_{i}(x) \right| U_{\mu_{\varepsilon,i}^{(2)},
  \xi_{\varepsilon,i}^{(2)}
  }+|\phi^{(1)}_{\varepsilon}|
  +|\phi^{(2)}_{\varepsilon} |
  \\
&\leq C
  \displaystyle\sum_{i=1}^{k}\eta^{(1)}_{i}(x) U_{\mu_{\varepsilon,i}^{(1)},
  \xi_{\varepsilon,i}^{(1)}
  }\Big(\mu_{\varepsilon,i}^{(1)}\big|
  \xi_{\varepsilon,i}^{(1)}-
  \xi_{\varepsilon,i}^{(2)}
  \big|+\big(\mu_{\varepsilon,i}^{(1)}
\big)^{-1}\big|\mu_{\varepsilon,i}^{(1)}-
\mu_{\varepsilon,i}^{(2)} \big|
\Big)
  \\
&\quad+O\bigg(\displaystyle\sum_{i=1}^{k}
\chi_{B_{\frac{5\delta}{2}}
\left(\xi_{\varepsilon,
i}^{(1)}\right)\backslash B_{\frac{\delta}{2}}\left(\xi_{\varepsilon,
i}^{(1)}\right)}
\left(\mu_{\varepsilon,i}^{(1)}
\right)^{\frac{2-N}{2}}
  \bigg)+|\phi^{(1)}_{\varepsilon}|
  +|\phi^{(2)}_{\varepsilon} |
  \\
&=
  o(1)\displaystyle\sum_{i=1}^{k}
  \eta^{(1)}_{i}(x) U_{\mu_{\varepsilon,i}^{(1)},
  \xi_{\varepsilon,i}^{(1)}
  }+O\bigg(\displaystyle\sum_{i=1}^{k}
\chi_{B_{\frac{5\delta}{2}}
\left(\xi_{\varepsilon,i}^{(1)}\right)\backslash B_{\frac{\delta}{2}}\left(\xi_{\varepsilon,
i}^{(1)}\right)}
\left(\mu_{\varepsilon,i}^{(1)}
\right)^{\frac{2-N}{2}}
  \bigg)
  +O\big(
  |\phi^{(1)}_{\varepsilon}|
  +|\phi^{(2)}_{\varepsilon}|\big)
 .
\end{align*}
\end{proof}
To obtain a more refined estimate for $\psi_\varepsilon$,
 we next decompose $\psi_\varepsilon$ into
\begin{align*}
\psi_\varepsilon(x)=\displaystyle\sum_{i=1}^{k}
\rho_{i,0}^{(1)}
\mu_{\varepsilon,i}^{(1)}Z_{\mu^{(1)}_{\varepsilon,i},
\xi_{\varepsilon,i}^{(1)}}^{0}
+\displaystyle\sum_{i=1}^{k}\displaystyle\sum_{l=1}
^{N}\rho_{i,l}^{(1)}
\left(\mu_{\varepsilon,i}^{(1)}\right)^{-1}
Z_{\mu^{(1)}_{\varepsilon,i},
\xi_{\varepsilon,i}^{(1)}}^{l}
+\psi_\varepsilon^*(x),
\end{align*}
where $\psi_\varepsilon^*(x)$ satisfies
\[
\int_{\mathbb{R}^N} W_{\mu_{\varepsilon,i}^{(1)},
\xi_{\varepsilon,i}^{(1)}
}^{p-1}
Z_{\mu^{(1)}_{\varepsilon,i},
\xi_{\varepsilon,i}^{(1)}}^{l}
\psi_\varepsilon^*(x)\mathrm{d}x=0,\;i=1,2,\cdots,k,\ l=0,1,\cdots,N.
\]
It is easy to check that   $\rho_{i,l}^{(1)}$ are bounded.  We first estimate the remainder term  $\psi_\varepsilon^*$.
Note that  $\psi_\varepsilon^*$ satisfies
\begin{small}
\begin{align*}
&-\Delta \psi^*_\varepsilon+V(x)\psi^*_\varepsilon
-(p-\varepsilon) a_\varepsilon(x)\psi^*_\varepsilon \\
=&(p-\varepsilon)\sum_{i=1}^{k}
\rho_{i,0}^{(1)}\eta_{i}^{(1)}
\mu_{\varepsilon,i}^{(1)}Z_{i,0}^{(1)}
\Big(a_\varepsilon(x)- \big(u^{(1)}_\varepsilon\big)^{p-1-\varepsilon}
\Big) \\
&+\sum_{i=1}^{k}
\rho_{i,0}^{(1)}\eta_{i}^{(1)}
\mu_{\varepsilon,i}^{(1)}Z_{i,0}^{(1)}
\Big((p-\varepsilon) \big(u^{(1)}_\varepsilon\big)^{p-1-\varepsilon}-p \big(u^{(1)}_\varepsilon\big)^{p-1}\Big) \\
&+p\sum_{i=1}^{k}
\rho_{i,0}^{(1)}\eta_{i}^{(1)}
\mu_{\varepsilon,i}^{(1)}Z_{i,0}^{(1)}
\Big(\big(u^{(1)}_\varepsilon\big)^{p-1}-  U_{\mu_{\varepsilon,i}^{(1)}, \xi_{\varepsilon,i}^{(1)}}^{p-1}\Big)
+2\sum_{i=1}^{k}
\rho_{i,0}^{(1)}\mu_{\varepsilon,i}^{(1)}
\nabla\eta_{i}^{(1)}\nabla
Z_{i,0}^{(1)} \\
&+\sum_{i=1}^{k}
\rho_{i,0}^{(1)}\mu_{\varepsilon,i}^{(1)}Z_{i,0}^{(1)}
\Delta\eta_{i}^{(1)}
-V(x)\sum_{i=1}^{k}
\rho_{i,0}^{(1)}\mu_{\varepsilon,i}^{(1)}
\eta_{i}^{(1)}
Z_{i,0}^{(1)} \\
&+(p-\varepsilon)
 \sum_{i=1}^{k}\sum_{l=1}
^{N}\rho_{i,l}^{(1)}\eta_{i}^{(1)}
\Big(\mu_{\varepsilon,i}^{(1)}\Big)^{-1}Z_{i,l}^{(1)}
\Big(a_\varepsilon(x)- \big(u^{(1)}_\varepsilon\big)^{p-1-\varepsilon}
\Big) \\
&+
\sum_{i=1}^{k}\sum_{l=1}
^{N}\rho_{i,l}^{(1)}\eta_{i}^{(1)}
\Big(\mu_{\varepsilon,i}^{(1)}\Big)^{-1}
Z_{i,l}^{(1)}
\Big((p-\varepsilon) \big(u^{(1)}_\varepsilon\big)^{p-1-\varepsilon}-p \big(u^{(1)}_\varepsilon\big)^{p-1}\Big) \\
&
+p
\sum_{i=1}^{k}\sum_{l=1}
^{N}\rho_{i,l}^{(1)}\eta_{i}^{(1)}
\Big(\mu_{\varepsilon,i}^{(1)}\Big)^{-1}
Z_{i,l}^{(1)}
\Big(\big(u^{(1)}_\varepsilon\big)^{p-1}- U_{\mu_{\varepsilon,i}^{(1)
},\xi_{\varepsilon,i}^{(1)}}^{p-1}\Big)
+2\sum_{i=1}^{k}\sum_{l=1}
^{N}\rho_{i,l}^{(1)}
\Big(\mu_{\varepsilon,i}^{(1)}\Big)^{-1}
\nabla\eta_{i}^{(1)}
\nabla Z_{i,l}^{(1)} \\
&
+\sum_{i=1}^{k}\sum_{l=1}
^{N}\rho_{i,l}^{(1)}
\Big(\mu_{\varepsilon,i}^{(1)}\Big)^{-1}
Z_{i,l}^{(1)}
\Delta\eta_{i}^{(1)}
-V(x)\sum_{i=1}^{k}\sum_{l=1}^{N}
\rho_{i,l}^{(1)}
\Big(\mu_{\varepsilon,i}^{(1)}\Big)^{-1}
\eta_{i}^{(1)}Z_{i,l}^{(1)} \\
&-(p-\varepsilon)
\sum_{i=1}^{k}\sum_{l=1}
^{N}\rho_{i,l}^{(1)}
\Big(\mu_{\varepsilon,i}^{(1)}\Big)^{-1}
\nabla\eta_{i}^{(1)}
U_{\mu_{\varepsilon,i}^{(1)
},\xi_{\varepsilon,i}^{(1)}}
\Big(a_\varepsilon(x)- \big(u^{(1)}_\varepsilon\big)^{p-1-\varepsilon}
\Big) \\
&-(p-\varepsilon)
\sum_{i=1}^{k}\sum_{l=1}
^{N}\rho_{i,l}^{(1)}
\Big(\mu_{\varepsilon,i}^{(1)}\Big)^{-1}
\nabla\eta_{i}^{(1)}
U_{\mu_{\varepsilon,i}^{(1)
},\xi_{\varepsilon,i}^{(1)}}
\big(u^{(1)}_\varepsilon
\big)^{p-1-\varepsilon} \\
&+
\sum_{i=1}^{k}\sum_{l=1}
^{N}\rho_{i,l}^{(1)}
\Big(\mu_{\varepsilon,i}^{(1)}\Big)^{-1}
\nabla\eta_{i}^{(1)} U_{\mu_{\varepsilon,i}^{(1)
},\xi_{\varepsilon,i}^{(1)}}^{p}
-2\sum_{i=1}^{k}\sum_{l=1}
^{N}\rho_{i,l}^{(1)}
\Big(\mu_{\varepsilon,i}^{(1)}\Big)^{-1}
\nabla^{2}\eta_{i}^{(1)}
\nabla U_{\mu_{\varepsilon,i}^{(1)
},\xi_{\varepsilon,i}^{(1)}} \\
&
-\sum_{i=1}^{k}\sum_{l=1}
^{N}\rho_{i,l}^{(1)}
\Big(\mu_{\varepsilon,i}^{(1)}\Big)^{-1}
U_{\mu_{\varepsilon,i}^{(1)
},\xi_{\varepsilon,i}^{(1)}}
\nabla^{3}\eta_{i}^{(1)}
+V(x)\sum_{i=1}^{k}\sum_{l=1}^{N}
\rho_{i,l}^{(1)}
\Big(\mu_{\varepsilon,i}^{(1)}\Big)^{-1}
\nabla\eta_{i}^{(1)} U_{\mu_{\varepsilon,i}^{(1)
},\xi_{\varepsilon,i}^{(1)}} \\
:=& L_1+L_2+\cdots+L_6+P_1+P_2+\cdots+P_6
+J_1+J_2+\cdots+J_6,
\end{align*}
\end{small}
where $Z_{i,0}^{(1)}=\frac{\partial U_{\mu_{\varepsilon,i}^{(1)
},\xi_{\varepsilon,i}^{(1)}}
}{\partial \mu_{\varepsilon,i}^{(1)
}}$ and $Z_{i,l}^{(1)}=\frac{\partial U_{\mu_{\varepsilon,i}^{(1)
},\xi_{\varepsilon,i}^{(1)}} }{\partial \xi_{\varepsilon,i,l}^{(1)}}$, $i=1,2,\cdots,k$, $l=1,\cdots,N$.
%\textcolor{blue}{
\begin{lemma}\label{lem5.05}
Under the same assumptions as in Lemma \ref{lem5.3},
it holds that
\begin{align*}
|\psi_\varepsilon^*|=
  o(1)
\displaystyle\sum_{i=1}^{k}
\frac{\big(\mu_{\varepsilon,i}^{(1)}
\big)^{\frac{N-2}{2}}}
{\left(1+\mu_{\varepsilon,i}^{(1)}
\big|x-\xi_{\varepsilon,i}^{(1)}\big |\right)
^{\frac{N+2\sigma}{2}}}.
\end{align*}
\end{lemma}
\begin{proof}[\bf Proof.]
Note that
\begin{align*}
\left|(a+b)^{p}-a^p\right|
\leq C\begin{cases}
\min\{a^{p-1}|b|,|b|^{p}\},\;\;&\text{if}\;\;p\leq1,\\[2mm]
a^{p-1}|b|+|b|^{p},\;\;&\text{if}\;\;p>1,
\end{cases}
\end{align*}
where $a>0, b\in\mathbb{R}$.
Then we have
\begin{align*}
|L_1|\leq& %C\displaystyle\sum_{i=1}^{k}
%\eta_{i}^{(1)}U_{\mu_{\varepsilon,i}^{(1)
%},\xi_{\varepsilon,i}^{(1)}}
%\left|a_\varepsilon(x)- \big(u^{(1)}_\varepsilon\big)
%^{p-1-\varepsilon}\right|
%\\
C\displaystyle\sum_{i=1}^{k}
\eta_{i}^{(1)}
U_{\mu_{\varepsilon,i}^{(1)
},\xi_{\varepsilon,i}^{(1)}}
\big|u^{(1)}_{\varepsilon}-
u^{(2)}_{\varepsilon}\big|^{p-1-\varepsilon}
\\&
= o(1)
\bigg(\displaystyle\sum_{i=1}^{k}
\frac{\big(\mu_{\varepsilon,i}^{(1)}\big)^{\frac{N-2}{2}}}
{\left(1+\big(\mu_{\varepsilon,i}^{(1)}\big)^{2}\big
|x-\xi_{\varepsilon,i}^{(1)}\big|^{2}\right)^{\frac{N-2}{2}}}
\bigg)^{p-\varepsilon}
+O\bigg( \displaystyle\sum_{i=1}^{k}
\eta_{i}^{(1)}
U_{\mu_{\varepsilon,i}^{(1)
},\xi_{\varepsilon,i}^{(1)}}
\left(\big|\phi^{(1)}_{\varepsilon}\big|
+\big|\phi^{(2)}_{\varepsilon}\big
|
\right)^{p-1-\varepsilon}\bigg)\\&\quad+O\bigg(\displaystyle\sum_{i=1}^{k}
\chi_{B_{2\delta}
\left(\xi_{\varepsilon,
i}^{(1)}\right)\backslash B_{\frac{\delta}{2}}\left(\xi_{\varepsilon,
i}^{(1)}\right)}
\left(\mu_{\varepsilon,i}^{(1)}
\right)^{-2}U_{\mu_{\varepsilon,i}^{(1)
},\xi_{\varepsilon,i}^{(1)}}
  \bigg)
\\
=&o(1)
\displaystyle\sum_{i=1}^{k}
\frac{\big(\mu_{\varepsilon,i}^{(1)}\big)^{\frac{N+2}{2}}}
{\left(1+\big(\mu_{\varepsilon,i}^{(1)}\big)^{2}
\big|x-\xi_{\varepsilon,i}^{(1)}\big |^{2}\right)
^{\frac{N+4}{4}+\frac{\sigma}{2}}}.
\end{align*}
By the mean value theorem, we find
\begin{align*}
|L_2|\leq& C\varepsilon\displaystyle\sum_{i=1}^{k}\eta_{i}^{(1)}
U_{\mu_{\varepsilon,i}^{(1)},\xi_{\varepsilon,i}^{(1)}}
(u_\varepsilon^{(1)})^{p-1-\theta\varepsilon}
\Big(1+\big|\log u_\varepsilon^{(1)}\big|\Big)\\
\leq&C\varepsilon|\log\varepsilon|
\displaystyle\sum_{i=1}^{k}\eta_{i}^{(1)}
U_{\mu_{\varepsilon,i}^{(1)},\xi_{\varepsilon,i}^{(1)}}^{p-\theta\varepsilon}
\bigg(1+\log\Big(1+\big(\mu_{i,
\varepsilon}^{(1)}\big)^{2}
\big|x-\xi_{\varepsilon,i}^{(1)}\big |^{2}\Big)\bigg)\\
\leq& C\varepsilon|\log\varepsilon|
\displaystyle\sum_{i=1}^{k}
\frac{\big(\mu_{\varepsilon,i}^{(1)}\big)^{\frac{N+2}{2}}}
{\left(1+\big(\mu_{\varepsilon,i}^{(1)}\big)^{2}
\big|x-\xi_{\varepsilon,i}^{(1)}\big |^{2}\right)
^{\frac{N+4}{4}+\frac{\sigma}{2}}},
\end{align*}
where $\theta\in(0,1)$.
In a similar way, we get
\begin{align*}
|L_3|\leq& C\displaystyle\sum_{i=1}^{k}\eta_{i}^{(1)} U_{\mu_{\varepsilon,i}^{(1)},\xi_{\varepsilon,i}^{(1)}}
\big|\phi^{(1)}_{\varepsilon}\big|^{p-1}
\leq C\big\|\phi^{(1)}_{\varepsilon}\big\|^{p-1}_*
\displaystyle\sum_{i=1}^{k}
\frac{\big(\mu_{\varepsilon,i}^{(1)}\big)^{\frac{N+2}{2}}}
{\left(1+\big(\mu_{\varepsilon,i}^{(1)}\big)^{2}
\big|x-\xi_{\varepsilon,i}^{(1)}\big |^{2}\right)
^{\frac{N+4}{4}+\frac{\sigma}{2}}}.
\end{align*}
For $\delta<\big|x-\xi_{\varepsilon,i}^{(1)}
\big|<2\delta$, $i=1,\cdots,k $, we infer
\begin{align*}
|L_4|\leq&
C
\displaystyle\sum_{i=1}^{k}
\frac{\big(\mu_{\varepsilon,i}^{(1)}\big)^{\frac{N+2}{2}}}
{\left(1+\big(\mu_{\varepsilon,i}^{(1)}\big)^{2}
\big|x-\xi_{\varepsilon,i}^{(1)}\big |^{2}\right)
^{\frac{N}{2}}}
\leq C\displaystyle\sum_{i=1}^{k}
\frac{1}{\big(\mu_{\varepsilon,i}^{(1)}
\big)^{\frac{N-4}{2}-\sigma}}
\frac{\big(\mu_{\varepsilon,i}^{(1)}\big)^{\frac{N+2}{2}}}
{\left(1+\big(\mu_{\varepsilon,i}^{(1)}
\big)^{2}
\big|x-\xi_{\varepsilon,i}^{(1)}\big |^{2}\right)
^{\frac{N+4}{4}+\frac{\sigma}{2}}},
\end{align*}
%For $\delta<\big|x-\xi_{\varepsilon,i}^{(1)}\big|<2\delta$, $i=1,\cdots,k $, we have
and
\begin{align*}
|L_5|\leq&
C
\displaystyle\sum_{i=1}^{k}
\frac{\big(\mu_{\varepsilon,i}^{(1)}\big)^{\frac{N-2}{2}}}
{\left(1+\big(\mu_{\varepsilon,i}^{(1)}\big)^{2}
\big|x-\xi_{\varepsilon,i}^{(1)}\big |^{2}\right)
^{\frac{N-2}{2}}}
\leq C\displaystyle\sum_{i=1}^{k}\frac{1}
{\big(\mu_{\varepsilon,i}^{(1)}\big)^{\frac{N-4}{2}-\sigma}}
\frac{\big(\mu_{\varepsilon,i}^{(1)}\big)^{\frac{N+2}{2}}}
{\left(1+\big(\mu_{\varepsilon,i}^{(1)}\big)^{2}
\big|x-\xi_{\varepsilon,i}^{(1)}\big|^{2}\right)
^{\frac{N+4}{4}+\frac{\sigma}{2}}
}.
\end{align*}
Moreover, we have
\begin{align*}
|L_6|\leq&
C
\displaystyle\sum_{i=1}^{k}
\frac{\big(\mu_{\varepsilon,i}^{(1)}\big)^{\frac{N-2}{2}}}
{\left(1+\big(\mu_{\varepsilon,i}^{(1)}\big)^{2}
\big|x-\xi_{\varepsilon,i}^{(1)}\big |^{2}\right)
^{\frac{N-2}{2}}}\leq C\displaystyle\sum_{i=1}^{k}\frac{1}
{\big(\mu_{\varepsilon,i}^{(1)}\big)^{1-\sigma}}
\frac{\big(\mu_{\varepsilon,i}^{(1)}\big)^{\frac{N+2}{2}}}
{\left(1+\big(\mu_{\varepsilon,i}^{(1)}\big)^{2}
\big|x-\xi_{\varepsilon,i}^{(1)}\big |^{2}\right)
^{\frac{N+4}{4}+\frac{\sigma}{2}}}.
\end{align*}
Hence we  derive  that
 \begin{align*}
\sum^6_{i=1}|L_i|=o(1)
\displaystyle\sum_{i=1}^{k}
\frac{\big(\mu_{\varepsilon,i}^{(1)}
\big)^{\frac{N+2}{2}}}
{\left(1+\big(\mu_{\varepsilon,i}^{(1)}\big)^{2}
\big|x-\xi_{\varepsilon,i}^{(1)}\big |^{2}\right)
^{\frac{N+4}{4}+\frac{\sigma}{2}}}.
\end{align*}
By similar arguments, we obtain
 \begin{align*}
\sum^6_{i=1}
|P_i|=o(1)
\displaystyle\sum_{i=1}^{k}
\frac{\big(\mu_{\varepsilon,i}^{(1)}\big)^{\frac{N+2}{2}}}
{\left(1+\big(\mu_{\varepsilon,i}^{(1)}\big)^{2}
\big|x-\xi_{\varepsilon,i}^{(1)}\big |^{2}\right)
^{\frac{N+4}{4}+\frac{\sigma}{2}}},
\end{align*}
and
\begin{align*}
\sum^6_{i=1}
|J_i|=o(1)
\displaystyle\sum_{i=1}^{k}
\frac{\big(\mu_{\varepsilon,i}^{(1)}\big)^{\frac{N+2}{2}}}
{\left(1+\big(\mu_{\varepsilon,i}^{(1)}\big)^{2}
\big|x-\xi_{\varepsilon,i}^{(1)}\big |^{2}\right)
^{\frac{N+4}{4}+\frac{\sigma}{2}}}
.
\end{align*}
As a result, we arrive at
$$
\sum^6_{i=1}\|L_i\|_{**}+\sum^6_{i=1}
\|P_i\|_{**}+\sum^6_{i=1}
\|J_i\|_{**}=o(1),
$$
which, together with Proposition 2.3 in \cite{CLl1},  gives that
\begin{align*}
\|\psi_\varepsilon^*\|_*\leq C\bigg(\sum^6_{i=1}\|L_i\|_{**}+\sum^6_{i=1}
\|P_i\|_{**}+\sum^6_{i=1}
\|J_i\|_{**}\bigg)=o(1).
\end{align*}
By Lemma \ref{lem: the computation}, we have
\begin{align*}
|\psi_\varepsilon^*|\leq&
%\bigg|\int_{\mathbb{R}^N}G(x,y) \left((p-\varepsilon) a_\varepsilon(x)\psi_\varepsilon+
%\sum^6_{i=1}L_i+\sum^6_{i=1}
%P_i+\sum^6_{i=1}
%J_i
% \right)\mathrm{d}y
%\bigg|\\
o(1)
\displaystyle\sum_{i=1}^{k}\int_{\mathbb{R}^N}\frac{1}
{|x-y|^{N-2}}
\frac{\big(\mu_{\varepsilon,i}^{(1)}\big)^{\frac{N+2}{2}}}
{\left(1+\big(\mu_{\varepsilon,i}^{(1)}\big)^{2}
\big|y-\xi_{\varepsilon,i}^{(1)}\big |^{2}\right)
^{\frac{N+4}{4}+\frac{\sigma}{2}}}\mathrm{d}y
\\
\leq&o(1)
\displaystyle\sum_{i=1}^{k}
\frac{\big(\mu_{\varepsilon,i}^{(1)}\big)^{\frac{N-2}{2}}}
{\left(1+\big(\mu_{\varepsilon,i}^{(1)}\big)
\big|x-\xi_{\varepsilon,i}^{(1)}\big |\right)
^{\frac{N+2\sigma}{2}}}.\end{align*}
\end{proof}

It follows from Lemma \ref{lem5.1} that $\psi_\varepsilon$  satisfies the following Pohozaev type identities
\begin{align}\label{a3.1}
&\bigg(\frac{N}{p+1-\varepsilon}-\frac{N-2}{2}
\bigg)
\int_{B_\delta(\xi_{\varepsilon,i}^{(1)})} b_\varepsilon(x)\psi_\varepsilon\mathrm{d}x
-\frac{1}{2}\int_{B_\delta(\xi_{\varepsilon,i}^{(1)})}\big\langle x-\xi_{\varepsilon,i}^{(1)},\nabla V(x)\big\rangle\Big(u_\varepsilon^{(1)}
+u_\varepsilon^{(2)}\Big)\psi_\varepsilon\mathrm{d}x
\notag\\
&\quad-\int_{B_\delta(\xi_{\varepsilon,i}^{(1)})}V(x)
\Big(u_\varepsilon^{(1)}+u_\varepsilon^{(2)}
\Big)\psi_\varepsilon\mathrm{d}x\notag
\\
&=
\int_{\partial B_\delta(\xi_{\varepsilon,i}^{(1)})}\frac{\partial\psi_\varepsilon}
{\partial\nu}\big\langle x-\xi_{\varepsilon,i}^{(1)},\nabla u_\varepsilon^{(1)}\big\rangle\mathrm{d}s
+\int_{\partial B_\delta(\xi_{\varepsilon,i}^{(1)})}\frac{\partial u_\varepsilon^{(2)}}{\partial\nu}\big\langle x-\xi_{\varepsilon,i}^{(1)},\nabla\psi_\varepsilon
\big\rangle\mathrm{d}s
\notag\\
&\quad+
\frac{1}{p+1-\varepsilon}\int_{\partial B_\delta(\xi_{\varepsilon,i}^{(1)})} b_\varepsilon(x)\psi_\varepsilon\langle x-\xi_{\varepsilon,i}^{(1)},\nu\rangle
\mathrm{d}s
+\frac{N-2}{2}\int_{\partial B_\delta(\xi_{\varepsilon,i}^{(1)})}\bigg(\frac{\partial\psi_\varepsilon}{\partial\nu}u_\varepsilon^{(1)}+\frac{\partial u_\varepsilon^{(2)}}{\partial\nu}
\psi_\varepsilon\bigg)\mathrm{d}s\notag\\
&\quad-\frac{1}{2}\int_{\partial B_\delta(\xi_{\varepsilon,i}^{(1)})}\nabla\Big(u_\varepsilon^{(1)}
+u_\varepsilon^{(2)}\Big)\nabla\psi_\varepsilon \langle x-\xi_{\varepsilon,i}^{(1)},\nu\rangle
\mathrm{d}s
-\frac{1}{2}\int_{\partial B_\delta(\xi_{\varepsilon,i}^{(1)})}V(x)\Big(u_\varepsilon^{(1)}
+u_\varepsilon^{(2)}\Big)\psi_\varepsilon \langle x-\xi_{\varepsilon,i}^{(1)},\nu\rangle
\mathrm{d}s,
\end{align}
and
\begin{equation}\label{a3.2}
\begin{split}
&\frac{1}{2}\int_{B_\delta(\xi_{\varepsilon,i}^{(1)})}\frac{\partial V(x)}{\partial x_{l}}\Big(u_\varepsilon^{(1)}
+u_\varepsilon^{(2)}\Big)
\psi_\varepsilon\mathrm{d}x
=\frac{1}{2}\int_{\partial B_\delta(\xi_{\varepsilon,i}^{(1)})}\nabla\Big(u_\varepsilon^{(1)}
+u_\varepsilon^{(2)}\Big)\nabla\psi_\varepsilon \nu_{l}\mathrm{d}s
-\int_{\partial B_\delta(\xi_{\varepsilon,i}^{(1)})}\frac{\partial u_\varepsilon^{(1)}}{\partial x_{l}}\frac{\partial \psi_\varepsilon}{\partial\nu}\mathrm{d}s
\\&
-\int_{\partial B_\delta(\xi_{\varepsilon,i}^{(1)})}\frac{\partial \psi_\varepsilon
}{\partial x_{l}}\frac{\partial u_\varepsilon^{(2)} }{\partial\nu}\mathrm{d}s
+\frac{1}{2}\int_{\partial B_\delta(\xi_{\varepsilon,i}^{(1)})}V(x)\Big(u_\varepsilon^{(1)}
+u_\varepsilon^{(2)}\Big)
\psi_\varepsilon
\nu_{l}\mathrm{d}s
-\frac{1}{p+1-\varepsilon}\int_{\partial B_\delta(\xi_{\varepsilon,i}^{(1)})} b_\varepsilon(x)\psi_\varepsilon \nu_{l}\mathrm{d}s,
\end{split}
\end{equation}
where $i=1,\cdots,k$, $l=1,\cdots,N$,
$$
b_\varepsilon(x)=(p+1-\varepsilon)\int_0^1
\Big(tu_\varepsilon^{(1)}+(1-t)u_\varepsilon^{(2)}
\Big)^{p-\varepsilon}\mathrm{d}t.
$$

Now  we  use  the Pohozaev identities \eqref {a3.1} and \eqref{a3.2}  to estimate the coefficients $\rho_{i,0}^{(1)}$ and $\rho_{i,j}^{(1)}$, $i=1,\cdots,k$, $j=1,\cdots,N$.
\begin{lemma}\label{lem5.6}
Let $N\geq6$ and $\xi_i^*$ $(i=1,2,\dots,k)$ be the non-degenerate critical points of $V(x)$.
 Assume that $V(x)\in C^2\big(B_{2\delta}(\xi^*_{i})\big)$ and $V(\xi_i^*)>0$.
Then we have
\begin{align}\label{tta}
\big|\rho_{i,0}^{(1)}\big|=o(1)\ \,\text{for}\;\;
i=1,\cdots,k.
\end{align}
\begin{proof}[\bf Proof.]
By a direct computation, we have
\begin{align*}
\int_{B_\delta(\xi_{\varepsilon,i}^{(1)})}  b_\varepsilon(x)\psi_\varepsilon\mathrm{d}x
%=(p+1-\varepsilon)\int_{B_\delta(\xi_{\varepsilon,i}^{(1)})}
%\int_0^1\left(tu_\varepsilon^{(1)}+(1-t)
%u_\varepsilon^{(2)}\right)^{p-\varepsilon}\mathrm{d}t
%\psi_\varepsilon\mathrm{d}x
=&(p+1-\varepsilon)\int_{B_\delta(\xi_{\varepsilon,i}^{(1)})} \big(u_\varepsilon^{(1)}\big)^{p-\varepsilon}
\psi_\varepsilon\mathrm{d}x
+O\bigg(\int_{B_\delta(\xi_{\varepsilon,i}^{(1)})} \big(u_\varepsilon^{(1)}\big)^{p-1-\varepsilon}
\big|u_\varepsilon^{(1)}-u_\varepsilon^{(2)}
\big|
\big|\psi_\varepsilon\big|\mathrm{d}x\bigg)\\
&+O\bigg(\int_{B_\delta(\xi_{i,
\varepsilon}^{(1)})}\big|u_\varepsilon^{(1)}
-u_\varepsilon^{(2)}\big|^{p-\varepsilon}
\big|\psi_\varepsilon\big|\mathrm{d}x\bigg).
\end{align*}
It follows from Lemma \ref{lem4.1}  and Lemma \ref{lem5.05} that
\begin{align*}
&\int_{B_\delta(\xi_{\varepsilon,i}^{(1)})}\big(u_\varepsilon^{(1)}\big)^{p-\varepsilon}
\psi_\varepsilon\mathrm{d}x\\%=\int_{B_\delta(\xi_{\varepsilon,i}^{(1)})}
%\big(u_\varepsilon^{(1)}\big)^{p-\varepsilon}
%\bigg(\rho_{i,0}^{(1)}\mu_{\varepsilon,i}^{(1)}Z_{i,0}^{(1)}+
%\displaystyle\sum_{l=1}^{N}
%\rho_{i,l}^{(1)}
%\big(\mu_{\varepsilon,i}^{(1)}\big)^{-1}Z_{i,l}^{(1)}+\psi^{\ast}_\varepsilon
%\bigg)
%\mathrm{d}x
=&\rho_{i,0}^{(1)}\mu_{\varepsilon,i}^{(1)}
\int_{B_\delta(\xi_{\varepsilon,i}^{(1)})}  \big(u_\varepsilon^{(1)}\big)^{p-\varepsilon}
Z_{i,0}^{(1)}\mathrm{d}x
+
\displaystyle\sum_{l=1}^{N}
\rho_{i,l}^{(1)}
\big(\mu_{\varepsilon,i}^{(1)}\big)^{-1}
\int_{B_\delta(\xi_{\varepsilon,i}^{(1)})}
\big(u_\varepsilon^{(1)}\big)
^{p-\varepsilon}
Z_{i,l}^{(1)}\mathrm{d}x+
o(1)
\\
=&\rho_{i,0}^{(1)}\mu_{\varepsilon,i}^{(1)}
\int_{B_\delta(\xi_{\varepsilon,i}^{(1)})} \bigg[\big(U_{\mu_{\varepsilon,i}^{(1)},
\xi_{\varepsilon,i}^{(1)}}
+\phi^{(1)}_{\varepsilon}\big)^{p-\varepsilon}
-U_{\mu_{\varepsilon,i}^{(1)},\xi_{\varepsilon,i}^{(1)}}
^{p-\varepsilon}\bigg]Z_{i,0}^{(1)}\mathrm{d}x
%\frac{\partial U_{\mu_{\varepsilon,i}^{(1)},\xi_{\varepsilon,i}^{(1)}
%}}
%{\partial\mu_{\varepsilon,i}^{(1)}}\mathrm{d}x
\\
&+\rho_{i,0}^{(1)}\mu_{\varepsilon,i}^{(1)}
\int_{B_\delta(\xi_{\varepsilon,i}^{(1)})} \bigg[U_{\mu_{\varepsilon,i}^{(1)},
\xi_{\varepsilon,i}^{(1)}}
^{p-\varepsilon}
-U_{\mu_{\varepsilon,i}^{(1)},
\xi_{\varepsilon,i}^{(1)}}
^{p}\bigg]
Z_{i,0}^{(1)}\mathrm{d}x\\
&+\displaystyle\sum_{l=1}^{N}
\rho_{i,l}^{(1)}
\big(\mu_{\varepsilon,i}^{(1)}\big)^{-1}
\int_{B_\delta(\xi_{\varepsilon,i}^{(1)})} \bigg[\big(U_{\mu_{\varepsilon,i}^{(1)},
\xi_{\varepsilon,i}^{(1)}}
+\phi^{(1)}_{\varepsilon}\big)^{p-\varepsilon}
-U_{\mu_{\varepsilon,i}^{(1)},
\xi_{\varepsilon,i}^{(1)}}
^{p-\varepsilon}\bigg]Z_{i,l}^{(1)}\mathrm{d}x
\\
&-\rho_{i,0}^{(1)}\mu_{\varepsilon,i}^{(1)}
\int_{\mathbb{R}^N\backslash B_\delta(\xi_{\varepsilon,i}^{(1)})} U_{\mu_{\varepsilon,i}^{(1)},\xi_{\varepsilon,i}^{(1)}
}
^{p}Z_{i,0}^{(1)}\mathrm{d}x
+
o(1)
\\
:=&\Gamma_{1}+\Gamma_{2}+\Gamma_{3}
%+O\Big(\big(\mu_{\varepsilon,i}^{(1)}
%\big)^{-N}\Big)
+o(1).
\end{align*}
Then we obtain
\begin{align*}
|\Gamma_1|\leq& C\int_{B_\delta(\xi_{\varepsilon,i}^{(1)})} U_{\mu_{\varepsilon,i}^{(1)},\xi_{\varepsilon,i}^{(1)}
}
^{p-\varepsilon}\big|
\phi^{(1)}_{\varepsilon}\big|\mathrm{d}x+C
\int_{B_\delta(\xi_{\varepsilon,i}^{(1)})} U_{\mu_{\varepsilon,i}^{(1)},\xi_{\varepsilon,i}^{(1)}
}
\big|\phi^{(1)}_{\varepsilon}
\big|^{p-\varepsilon}\mathrm{d}x\\
\leq& C\Big(\big\|\phi^{(1)}_{\varepsilon}\big\|_*
+\big\|\phi^{(1)}_{\varepsilon}
\big\|_*^{p-\varepsilon}\Big).
\end{align*}
By the mean value theorem, we find
$$
|\Gamma_2|\leq C\varepsilon \int_{B_\delta(\xi_{\varepsilon,i}^{(1)})} U_{\mu_{\varepsilon,i}^{(1)},\xi_{\varepsilon,i}^{(1)}
}
^{p+1-\theta\varepsilon}\big|\log U_{\mu_{\varepsilon,i}^{(1)},\xi_{\varepsilon,i}^{(1)}
}\big|\mathrm{d}x\leq C\varepsilon|\log\varepsilon|,
$$
where $\theta\in(0,1)$.
Similarly,  we can get
$$|\Gamma_3|\leq  C\big\|\phi^{(1)}_{\varepsilon}\big\|_*
+C\big\|\phi^{(1)}_{\varepsilon}\big
\|_*^{p-\varepsilon}.$$
In view of Lemmas \ref{lem4.1} and  \ref{cor5.4}, we  infer
\begin{align*}
 &\int_{B_\delta(\xi_{\varepsilon,i}^{(1)})} (u_\varepsilon^{(1)})^{p-1-\varepsilon}
 \big|u_\varepsilon^{(1)}-u_\varepsilon^{(2)}
 \big||\psi_\varepsilon|\mathrm{d}x\\
=&o(1)\int_{
B_\delta(\xi_{\varepsilon,i}^{(1)})} U_{\mu_{\varepsilon,i}^{(1)},\xi_{\varepsilon,i}^{(1)}
}
^{p-\varepsilon}
|\psi_\varepsilon|\mathrm{d}x+O
\bigg(\frac{1}{\big(\mu_{\varepsilon,i}^{(1)}
\big)^{\frac{N+2}{2}+\sigma}}
\bigg)
+O\bigg(\int_{B_\delta(\xi_{\varepsilon,i}^{(1)})} U_{\mu_{\varepsilon,i}^{(1)},\xi_{\varepsilon,i}^{(1)}
}^{p-1-\varepsilon}
|\psi_\varepsilon|\left(\big|
\phi^{(1)}_{\varepsilon}\big|
+\big|\phi^{(2)}_{\varepsilon}\big|
\right)\mathrm{d}x\bigg)
\\
=&o(1)+O\left(
\big\|\phi^{(1)}_{\varepsilon}\big\|_*
\right)+O\left(\big
\|\phi^{(2)}_{\varepsilon}\big\|_*\right),
\end{align*}
and
\begin{align*}
 &\int_{B_\delta(\xi_{\varepsilon,i}^{(1)})} \big|u_\varepsilon^{(1)}-u_\varepsilon^{(2)}\big|
 ^{p-\varepsilon}|\psi_\varepsilon|
 \mathrm{d}x
=o(1)+O\left(
\big\|\phi^{(1)}_{\varepsilon}\big\|_*
\right)+O\left(\big
\|\phi^{(2)}_{\varepsilon}\big\|_*\right).
\end{align*}
Hence, we deduce
$$
\int_{B_\delta(\xi_{\varepsilon,i}^{(1)})} b_\varepsilon(x)\psi_\varepsilon\mathrm{d}x
=o(1).
$$
On the other hand, we have
\begin{align*}
&-\int_{B_\delta(\xi_{\varepsilon,i}^{(1)})}  V(x)\big(u_\varepsilon^{(1)}
+u_\varepsilon^{(2)}\big)\psi_\varepsilon\mathrm{d}x
\\
=&-2\rho_{i,0}^{(1)}\mu_{\varepsilon,i}^{(1)}
\int_{B_\delta(\xi_{\varepsilon,i}^{(1)})}V(x)u_\varepsilon^{(1)}Z_{i,0}^{(1)}\mathrm{d}x
+\rho_{i,0}^{(1)}\mu_{\varepsilon,i}^{(1)}
\int_{B_\delta(\xi_{\varepsilon,i}^{(1)})}V(x)
\big(u_\varepsilon^{(1)}-
u_\varepsilon^{(2)}\big)
Z_{i,0}^{(1)}\mathrm{d}x
\\
&-2\big(\mu_{\varepsilon,i}^{(1)}\big)^{-1}\displaystyle\sum_{l=1}
^{N}\rho_{i,l}^{(1)}\int_{B_\delta(\xi_{\varepsilon,i}^{(1)})}V(x)
u_\varepsilon^{(1)}Z_{i,l}^{(1)}\mathrm{d}x
+\big(\mu_{\varepsilon,i}^{(1)}\big)^{-1}\displaystyle\sum_{l=1}
^{N}\rho_{i,l}^{(1)}\int_{B_\delta(\xi_{\varepsilon,i}^{(1)})}
V(x)\big(u_\varepsilon^{(1)}-u_\varepsilon^{(2)}
\big)Z_{i,l}^{(1)}\mathrm{d}x
\\
&-
\int_{B_\delta(\xi_{\varepsilon,i}^{(1)})}V(x)
\big(u_\varepsilon^{(1)}
+u_\varepsilon^{(2)}\big)\psi_\varepsilon^*\mathrm{d}x\\
:=&T_{1}+T_{2}+T_{3}+T_{4}+T_{5}.
\end{align*}
Next, we estimate them term by term.
A direct calculation gives
\begin{align*}
T_{1}=&-2\rho_{i,0}^{(1)}\mu_{\varepsilon,i}^{(1)}
\int_{B_\delta(\xi^{(1)}_{\varepsilon,i})}V(x)
U_{\mu_{\varepsilon,i}^{(1)},\xi_{\varepsilon,i}^{(1)}
}
Z_{i,0}^{(1)}\mathrm{d}x-
2\rho_{i,0}^{(1)}
\mu_{\varepsilon,i}^{(1)}
\int_{B_\delta(\xi_{\varepsilon,i}^{(1)})}V(x)Z_{i,0}^{(1)}
\phi^{(1)}_{\varepsilon}\mathrm{d}x\\
%=&-2\rho_{i,0}^{(1)}\mu_{\varepsilon,i}^{(1)}
%\int_{B_\delta(\xi_{\varepsilon,i}^{(1)})}
%\Big(V(\xi_{\varepsilon,i}^{(1)})+
%\big\langle\nabla V(\xi_{\varepsilon,i}^{(1)}),x-\xi_{\varepsilon,i}^{(1)}
%\big\rangle
%+
%O\left
%(\big|x-\xi_{\varepsilon,i}^{(1)}\big|^{2}\right)\Big)
%U_{\mu_{\varepsilon,i}^{(1)},\xi_{\varepsilon,i}^{(1)}
%}
%Z_{i,0}^{(1)}\mathrm{d}x\\&
%+O\bigg(\int_{B_\delta(\xi_{\varepsilon,i}^{(1)})}
%U_{\mu_{\varepsilon,i}^{(1)},\xi_{\varepsilon,i}^{(1)}
%}
%\big|\phi^{(1)}_{\varepsilon}\big|
%\mathrm{d}x\bigg)\\
=&D\rho_{i,0}^{(1)}
(\mu_{\varepsilon,i}^{(1)})^{-2}
V(\xi_{\varepsilon,i}^{(1)})
+O\left(\varepsilon^2|\log\varepsilon|
\right)+
O\left(\varepsilon\big\|\phi^{(1)}_{\varepsilon}\big\|_*
\right)\\
=&D\rho_{i,0}^{(1)}
(\mu_{\varepsilon,i}^{(1)})^{-2}
V(\xi_{\varepsilon,i}^{(1)})
+O\left(\varepsilon^{2-\frac{\sigma}{2}}
\right),
%%(2-N)\rho_{i,0}^{(1)}\alpha_{N}^2
%%(\mu_{\varepsilon,i}^{(1)})^{-2}
%%V(\xi_{\varepsilon,i}^{(1)})
%%\int_{\mathbb{R}^N}
%%\frac{1-|x|^{2}}{(1+|x|^{2})
%%^{N-1}}\mathrm{d}x
%%+O\left(\varepsilon^2|\log\varepsilon|
%%\right)+
%%O\left(\varepsilon\big\|\phi^{(1)}_{\varepsilon}\big\|_*
%%\right).
\end{align*}
where $D=-2\int_{\mathbb{R}^N}U_{1,0}\frac{\partial U_{\mu,0}}{\partial\mu}|_{\mu=1}\mathrm{d}x.$
By Proposition \ref{prop2.4} and Lemma \ref{cor5.4}, we derive
\begin{align*}
T_{2}= &
%C\int_{B_\delta(\xi_{\varepsilon,i}^{(1)})}\big|
%u_\varepsilon^{(1)}-u_\varepsilon^{(2)}
%\big|
%U_{\mu_{\varepsilon,i}^{(1)},\xi_{\varepsilon,i}^{(1)}
%}\mathrm{d}x
% \\
o(1)\int_{B_\delta
(\xi_{\varepsilon,i}^{(1)})}
U_{\mu_{\varepsilon,i}^{(1)},\xi_{\varepsilon,i}^{(1)}
}^2\mathrm{d}x
+O\bigg(\int_{B_\delta(\xi_{\varepsilon,i})}
U_{\mu_{\varepsilon,i}^{(1)},\xi_{\varepsilon,i}^{(1)}
}
\left(\big|\phi^{(1)}_{\varepsilon}\big|
+\big|\phi^{(2)}_{\varepsilon}\big|
\right)\mathrm{d}x\bigg)+
O\left( \mu_{\varepsilon,i}^{2-N}
\right)\\
=&o(\varepsilon)
+O\left(\varepsilon
\big\|\phi^{(1)}_{\varepsilon}\big\|_*
\right)+O\left(\varepsilon\big
\|\phi^{(2)}_{\varepsilon}\big\|_*\right)
=o(\varepsilon),
\end{align*}
and
\begin{align*}
T_{3}=&-2\big(\mu_{\varepsilon,i}^{(1)}\big)^{-1}\displaystyle\sum_{l=1}
^{N}\rho_{i,l}^{(1)}
\int_{B_\delta(\xi_{\varepsilon,i}^{(1)})}
V(x)\left(
U_{\mu_{\varepsilon,i}^{(1)},\xi_{\varepsilon,i}^{(1)}
}+\phi^{(1)}_{\varepsilon}\right)
Z_{i,l}^{(1)}\mathrm{d}x\\
%=&-2\big(\mu_{\varepsilon,i}^{(1)}\big)^{-1}
%\displaystyle\sum_{l=1}
%^{N}\rho_{i,l}^{(1)}
%\int_{B_\delta(\xi_{\varepsilon,i}^{(1)})}\Big(\big\langle\nabla V(\xi_{\varepsilon,i}^{(1)}),x-\xi_{\varepsilon,i}^{(1)}\big\rangle+O\big(|x-\xi_{i}|^{2}\big)\Big)
%U_{\mu_{\varepsilon,i}^{(1)},\xi_{\varepsilon,i}^{(1)}
%}
%Z_{i,l}^{(1)}\mathrm{d}x\\
%&
%+
%O\left(\big(\mu_{\varepsilon,i}^{(1)}\big)^{-2}
%\big\|\phi^{(1)}_{\varepsilon}\big\|_*
%\right)\\
=&o(\varepsilon^2)
+
O\left(\varepsilon^{2-\frac{\sigma}{2}}
\right)
%=&\frac{(4-2N)}{N}\frac{\alpha_{N}^2}
%{\big(\mu_{\varepsilon,i}^{(1)}\big)^{3}}
%\displaystyle\sum_{l=1}
%^{N}\rho_{i,l}^{(1)}
%\frac{\partial V(\xi_{\varepsilon,i}^{(1)})}{\partial x_{l}}
%\int_{\mathbb{R}^N}
%\frac{|x|^{2}}{(1+|x|^{2})
%^{N-1}}\mathrm{d}x
%+
%O\left(\varepsilon^{2-\frac{\sigma}{2}}
%\right)\\
=O\left(\varepsilon^{2-\frac{\sigma}{2}}
\right).
\end{align*}
It follows from analogous  arguments that
$$
T_{4}=o(\varepsilon)\;\;\text{and}\;\;T_{5}=o(\varepsilon).
$$
%Lemma \ref{lem4.1} and Lemma \ref{lem5.05} allow us to obtain
%\begin{align*}
%T_{5}
%=o_{\varepsilon}(1)
%\int_{B_\delta(\xi_{\varepsilon,i}^{(1)})}
%\frac{\big(\mu_{\varepsilon,i}^{(1)}\big)^{N-2}}
%{\left(1+\big(\mu_{\varepsilon,i}^{(1)}\big)
%^{2}|x-\xi_{\varepsilon,i}^{(1)}|^{2}
%\right)^{\frac{3N-4}{4}
%+\frac{\sigma}{2}}}\mathrm{d}x
%=o_{\varepsilon}(\varepsilon).
%\end{align*}
Moreover, we  have
\begin{align*}
&\int_{B_\delta(\xi_{\varepsilon,i}^{(1)})} \big\langle x-\xi_{\varepsilon,i}^{(1)},\nabla V(x)\big\rangle\big(u_\varepsilon^{(1)}
+u_\varepsilon^{(2)}\big)\psi_\varepsilon\mathrm{d}x
\\
=&\int_{B_\delta(\xi_{\varepsilon,i}^{(1)})}\big\langle x-\xi_{\varepsilon,i}^{(1)},
\nabla V(x)\big\rangle\left(2u_\varepsilon^{(1)}+ \big(u_\varepsilon^{(2)}-u_\varepsilon^{(1)}
\big)\right)\psi_\varepsilon\mathrm{d}x
\\
=&\rho_{i,0}^{(1)}\mu_{\varepsilon,i}^{(1)}\int_{B_\delta(\xi_{\varepsilon,i}^{(1)})}
\big\langle x-\xi_{\varepsilon,i}^{(1)},\nabla V(x)\big\rangle
\left(2u_\varepsilon^{(1)}+ \big(u_\varepsilon^{(2)}-u_\varepsilon^{(1)}
\big)\right)
Z_{i,0}^{(1)}\mathrm{d}x\\
&+\big(\mu_{\varepsilon,i}^{(1)}\big)^{-1}\displaystyle\sum_{l=1}
^{N}\rho_{i,l}^{(1)}\int_{B_\delta(\xi_{\varepsilon,i}^{(1)})}
\big\langle x-\xi_{\varepsilon,i}^{(1)},\nabla V(x)\big\rangle
\left(2u_\varepsilon^{(1)}+ \big(u_\varepsilon^{(2)}-u_\varepsilon^{(1)}
\big)\right)
Z_{i,l}^{(1)}\mathrm{d}x\\
&
+\int_{B_\delta(\xi_{\varepsilon,i}^{(1)})}
\big\langle x-\xi_{\varepsilon,i}^{(1)},\nabla V(x)\big\rangle
\left(2u_\varepsilon^{(1)}+ \big(u_\varepsilon^{(2)}-u_\varepsilon^{(1)}
\big)\right)\psi^*_\varepsilon\mathrm{d}x\\
=&o(\varepsilon).
\end{align*}
Using Lemma \ref{lem4.1}, we conclude that
\begin{align*}
\text{RHS \;of \;}\eqref{a3.1}=O\big(\varepsilon^{\frac{N-2}{2}}\big).
\end{align*}
Thus, in view of \eqref{a3.1}, we obtain   \eqref{tta} and  the proof of this lemma is complete.
\end{proof}\end{lemma}
\begin{lemma}\label{lem5.61}
Let $N\geq6$ and $\xi_i^*$ $(i=1,2,\dots,k)$ be the non-degenerate critical points of $V(x)$.
 Assume that $V(x)\in C^2\big(B_{2\delta}(\xi^*_{i})\big)$ and $V(\xi_i^*)>0$.
 Then
$$\rho_{i,j}^{(1)}=
  o(1),\ \
i=1,\cdots,k,\ \ j=1,\cdots,N.
$$
\end{lemma}
\begin{proof}[\bf Proof.]
According to Proposition \ref{prop2.4} and Lemma \ref{lem5.05}, we have
\begin{align*}
&\frac{1}{2}
\displaystyle\int_{B_{\delta}
(\xi_{\varepsilon,i}^{(1)})}\frac{\partial V(x)}{\partial x_{l}}\Big(u_\varepsilon^{(1)}+
u_\varepsilon^{(2)}\Big)
\psi_\varepsilon\mathrm{d}x\\
=&
\displaystyle\int_{B_{\delta}
(\xi_{\varepsilon,i}^{(1)})}\frac{\partial V(x)}{\partial x_{l}}\Big(u_\varepsilon^{(1)}
+\frac{u_\varepsilon^{(2)}-
u_\varepsilon^{(1)}}{2}\Big)\Big(
\rho_{i,0}^{(1)}
\mu_{\varepsilon,i}^{(1)}Z_{i,0}^{(1)}
+\displaystyle\sum_{j=1}
^{N}{\rho_{i,j}^{(1)}\big(\mu_{\varepsilon,i}^{(1)}\big)^{-1}
Z_{i,j}^{(1)}}
\Big)
\mathrm{d}x\\
&
+\frac{1}{2}\displaystyle\int_{B_{\delta}
(\xi_{\varepsilon,i}^{(1)})}\frac{\partial V(x)}{\partial x_{l}}\Big(u_\varepsilon^{(1)}
+u_\varepsilon^{(2)}\Big)
\psi^{\ast}_\varepsilon\mathrm{d}x
\\
=&
\int_{B_{\delta}(\xi_{\varepsilon,i}^{(1)})}\frac{\partial V(x)}{\partial x_{l}}u_\varepsilon^{(1)}
\Big(
\rho_{i,0}^{(1)}
\mu_{\varepsilon,i}^{(1)}Z_{i,0}^{(1)}
+\displaystyle\sum_{j=1}
^{N}{\rho_{i,j}^{(1)}\big(\mu_{\varepsilon,i}^{(1)}\big)^{-1}
Z_{i,j}^{(1)}}
\Big)
\mathrm{d}x
\\
&+
\frac{1}{2}
\int_{B_{\delta}(\xi_{\varepsilon,i}^{(1)})}
\Big(\frac{\partial V(x)}{\partial x_{l}}
-\frac{\partial V\big(\xi_{\varepsilon,i}^{(1)}\big)}{\partial x_{l}}+\frac{\partial V\big(\xi_{\varepsilon,i}^{(1)}\big)}{\partial x_{l}}
\Big)\Big(u_\varepsilon^{(1)}+
u_\varepsilon^{(2)}\Big)
\psi_\varepsilon^*(x)\mathrm{d}x\\
&+O\Bigg(
\int_{B_\delta(\xi_{\varepsilon,i}^{(1)})}
\bigg|\frac{\partial V(x)}{\partial x_{l}}
-\frac{\partial V\big(\xi_{\varepsilon,i}^{(1)}\big)}{\partial x_{l}}+\frac{\partial V\big(\xi_{\varepsilon,i}^{(1)}\big)}{\partial x_{l}}
\bigg
|
\big|u_\varepsilon^{(2)}-u_\varepsilon^{(1)}
\big|
U_{\mu_{\varepsilon,i}^{(1)},\xi_{\varepsilon,i}^{(1)}
}
\mathrm{d}x\Bigg)
\\
=&\int_{B_{\delta}(\xi_{\varepsilon,i}^{(1)})} \frac{\partial V(x)}{\partial x_{l}}
\Big(U_{\mu_{\varepsilon,i}^{(1)},
\xi_{\varepsilon,i}^{(1)}}
+\phi^{(1)}_{\varepsilon}\Big)
\Big(
\rho_{i,0}^{(1)}\mu_{\varepsilon,i}^{(1)}Z_{i,0}^{(1)}
+\displaystyle\sum_{j=1}
^{N}{\rho_{i,j}^{(1)}\big(\mu_{\varepsilon,i}^{(1)}\big)^{-1}
Z_{i,j}^{(1)}}
\Big)\mathrm{d}x+o\big(\varepsilon^{\frac{3}{2}
}
\big)\\
=&\rho_{i,0}^{(1)}\mu_{\varepsilon,i}^{(1)}
\int_{B_{\delta}(\xi_{\varepsilon,i}^{(1)})}
\frac{\partial V(x)}{\partial x_{l}}U_{\mu_{\varepsilon,i}^{(1)},
\xi_{\varepsilon,i}^{(1)}}
Z_{i,0}^{(1)}\mathrm{d}x
+\displaystyle\sum_{j=1}
^{N}\rho_{i,j}^{(1)}\big(\mu_{\varepsilon,i}^{(1)}\big)^{-1}
\int_{B_{\delta}(\xi_{\varepsilon,i}^{(1)})} \frac{\partial V(x)}{\partial x_{l}}
U_{\mu_{\varepsilon,i}^{(1)},\xi_{\varepsilon,i}^{(1)}
}
Z_{i,j}^{(1)}\mathrm{d}x
+
o\big(\varepsilon^{\frac{3}{2}}
\big)\\
%&
%+\displaystyle\sum_{j=1}
%^{N}\rho_{i,j}^{(1)}\big(\mu_{\varepsilon,i}^{(1)}\big)^{-1}
%\int_{B_{\delta}(\xi_{\varepsilon,i}^{(1)})} \frac{\partial V(x)}{\partial x_{l}}
%U_{\mu_{\varepsilon,i}^{(1)},\xi_{\varepsilon,i}^{(1)}
%}
%Z_{i,j}^{(1)}\mathrm{d}x
%\\
:=&F_{1}+F_{2}+o\big(\varepsilon^{\frac{3}{2}
}
\big).
\end{align*}
By  Proposition \ref{prop2.4} and Lemma \ref{lem5.6}, we  infer
\begin{align*}
F_{1}
%\rho_{i,0}^{(1)}\mu_{\varepsilon,i}^{(1)}
%\displaystyle\int_{
%B_{\delta}(\xi_{\varepsilon,i}^{(1)})}
%\bigg(
%\frac{\partial V\big(\xi_{\varepsilon,i}^{(1)}\big)}{\partial x_{l}}+
%\frac{\partial V(x)}{\partial x_{l}}-
%\frac{\partial V\big(\xi_{\varepsilon,i}^{(1)}\big)}{\partial x_{l}}
%\bigg)U_{\mu_{\varepsilon,i}^{(1)},
%\xi_{\varepsilon,i}^{(1)}}
%Z_{i,0}^{(1)}\mathrm{d}x
%\\
=&\rho_{i,0}^{(1)}\mu_{\varepsilon,i}^{(1)}\frac{\partial V\big(\xi_{\varepsilon,i}^{(1)}\big)}{\partial x_{l}}\int_{B_{\delta}(\xi_{\varepsilon,i}^{(1)})}U_{\mu_{\varepsilon,i}^{(1)},\xi_{\varepsilon,i}^{(1)}
}
Z_{i,0}^{(1)}\mathrm{d}x+O\left(\big|\rho_{i,0}^{(1)}\big|\big(\mu_{\varepsilon,i}^{(1)}\big)^{-3}\right)\\
=&\frac{\partial V\big(\xi_{\varepsilon,i}^{(1)}\big) }{\partial x_{l}} O\left(\big(\mu_{\varepsilon,i}^{(1)}\big)^{-2}\right)
+O\left(\big(\mu_{\varepsilon,i}^{(1)}\big)^{-3}\right)\\
=&o\big(\varepsilon^{\frac{3}{2}
}
\big)
\end{align*}and
\begin{align*}
F_{2}
=&\displaystyle\sum_{j=1}^{N}\rho_{i,j}^{(1)}\big(\mu_{\varepsilon,i}^{(1)}\big)^{-1}
\int_{B_{\varepsilon^{\tau_{0}}}
(\xi_{\varepsilon,i}^{(1)})}
\frac{\partial^{2} V\big(\xi_{\varepsilon,i}^{(1)}\big)}{\partial x_{l}\partial
x_{j}}
\big(x_{j}-\xi_{i,j,\varepsilon}^{(1)}\big) U_{\mu_{\varepsilon,i}^{(1)},\xi_{\varepsilon,i}^{(1)}
}
Z_{i,j}^{(1)}\mathrm{d}x +o\big(
\varepsilon^{\frac{3}{2}}\big)
\\
=&\displaystyle\sum_{j=1}
^{N} D_j
\rho_{i,j}^{(1)}\big(\mu_{\varepsilon,i}^{(1)}\big)^{-3}
\frac{\partial^{2} V(\xi_{i}^{\ast})}{\partial x_{l}\partial
x_{j}}
+o\big(
\varepsilon^{\frac{3}{2}}\big)
,\\
%=&\alpha_{N}^{2}(N-2)\displaystyle\sum_{j=1}
%^{N}
%\rho_{i,j}^{(1)}\big(\mu_{\varepsilon,i}^{(1)}\big)^{-3}
%\frac{\partial^{2} V(\xi_{i}^{\ast})}{\partial x_{l}\partial
%x_{j}}
%\int_{\mathbb{R}^N}
%\frac{x_{j}^{2}}
%{(1+|x|^{2})^{N-1}}\mathrm{d}x
%+o\bigg(\displaystyle\sum_{l=1}
%^{N}\big|\rho_{i,j}^{(1)}\big|
%\varepsilon^{\frac{3}{2}}\bigg)
%,
\end{align*}
where $D_j=\int_{\mathbb{R}^N}U_{1,0}x_j\frac{\partial U_{1,0}}{\partial x_j}\mathrm{d}x.$
Using Lemma \ref{lem4.1}, we find  that
\[
\text{RHS \;of \;}\eqref{a3.2}=O\Big(\varepsilon^{
\frac{N-2}{2}}
\Big).
\]
As a consequence,  by \eqref{a3.2},  we obtain
$$\rho_{i,j}^{(1)}=
  o(1)
,\ \
i=1,\cdots,k,\ \ j=1,\cdots,N.
$$
This concludes  the proof of  Lemma \ref{lem5.61}.
\end{proof}
\begin{proof}[\bf Proof of Theorem \ref{athm1.1}]
It follows from Lemmas \ref{lem: the computation} and  \ref{lem4.1} that
$$
|\psi_\varepsilon|\leq C\int_{\mathbb{R}^N}\frac{1}{|x-y|^{N-2}} a_\varepsilon(y)\psi_\varepsilon(y)\mathrm{d}y
\leq C\|\psi_\varepsilon\|_*
\displaystyle\sum_{i=1}^{k}\frac{\mu_{\varepsilon,i}^{\frac{N-2}{2}}}
{\big(1+\mu_{\varepsilon,i}^{2}
|x-\xi_{\varepsilon,i}|^{2}\big)^{\frac{\beta}{2}}},
$$
where $\beta\in(\frac{N-2}{2}+\sigma,N-2)$. Then  we have

\begin{align*}
\bigg(
\displaystyle\sum_{i=1}^{k}
\frac{\mu_{\varepsilon,i}^{\frac{N-2}{2}}}
{\big(1+\mu_{\varepsilon,i}^{2}|x-\xi_{\varepsilon,i}|^{2}\big)
^{\frac{N-2}{4}+\frac{\sigma}{2}}}\bigg)^{-1}
|\psi_\varepsilon(x)|
\leq  C\|\psi_\varepsilon\|_*
\displaystyle\sum_{i=1}^{k}\frac{1}
{\big(1+\mu_{\varepsilon,i}^{2}
|x-\xi_{\varepsilon,i}|^{2}\big)
^{\frac{\beta}
{2}-\frac{N-2}{4}-\frac{\sigma}{2}}}.
\end{align*}
By Lemmas \ref{lem5.05},  \ref{lem5.6} and  \ref{lem5.61}, we deduce that for any $R>0$,
$$\mu_{\varepsilon,i}^{-\frac{N-2}{2}}
\psi_\varepsilon\rightarrow0\,\,\,\mbox{in}\,\,\, B_{R\mu_{\varepsilon,i}^{-1}}
(\xi_{\varepsilon,i}),\ \  i=1,\cdots,k,$$
which, together with $\|\psi_\varepsilon\|_*=1$,  yields that  $\bigg(\displaystyle\sum_{i=1}^{k}
\frac{\mu_{\varepsilon,i}^{\frac{N-2}{2}}}
{\big(1+\mu_{\varepsilon,i}^{2}|x-\xi_{\varepsilon,i}|^{2}\big)
^{\frac{N-2}{4}+\frac{\sigma}{2}}}\bigg)^{-1}|\psi_\varepsilon(x)|$
cannot attain its maximum in $B_{R\mu_{\varepsilon,i}^{-1}}(\xi_{\varepsilon,i})$, $i=1,2,\cdots,k$. Choosing $R$ large enough such that
\[
\frac{C}{(1+R^{2})^{\frac{\beta}{2}-\frac{N-2}{4}
+\frac{\sigma}{2}}}\leq\frac{1}{2},
\]
then we conclude
\[
\|\psi_\varepsilon\|_*\leq\frac{1}{2}\|\psi_\varepsilon\|_*.
\]
Thus we obtain a contradiction and complete the proof of Theorem \ref{athm1.1}.
\end{proof}
\section{Morse index of the peak solutions: Proof of Theorems \ref{athm1.2} and \ref{athm1.3}}
In this section, we aim to compute the Morse index of the peak solutions $u_\varepsilon$ by  Pohozaev type identities.
To this end, we need to consider  the weighted eigenvalue problem
\begin{equation}\label{aa4.1}
\begin{cases}
-\Delta\omega+V(x)\omega-(p-\varepsilon)\lambda u_\varepsilon^{p-1-\varepsilon}\omega=0
\;\;\text{in}\;\;\mathbb{R}^N,\\[1mm]
\omega\in H_V^1(\mathbb{R}^N).
\end{cases}
\end{equation}
It is well-known that problem \eqref{aa4.1} admits a sequence of non-decreasing eigenvalues $\{\lambda_{\varepsilon,l}\}$ along with their  associated eigenfunctions $\{\omega_{\varepsilon,l}\}$,
$l=1,2\cdots$. Without loss of generality, we may assume that $\omega_{\varepsilon,l}$ satisfies
\[
\int_{\mathbb{R}^N} u_\varepsilon^{p-1-\varepsilon}\omega_{\varepsilon,l} \omega_{\varepsilon,k}\mathrm{d}x=0\;\;\text{for}\;l\neq k.
\]
Recall that
\begin{align*}
u_{\varepsilon}=\sum_{i=1}^k W_{\mu_{\varepsilon,i},\xi_{\varepsilon,i}}
+\phi_\varepsilon,
\end{align*}
where $W_{\mu_{\varepsilon,i},\xi_{\varepsilon,i}}=\eta_i U_{\mu_{\varepsilon,i},\xi_{\varepsilon,i}}$,
$\|\phi_{\varepsilon}\|_{*}=O(\varepsilon^{1-\frac{\sigma}{2}})$.
In order to compute the Morse index of  $u_\varepsilon$, it suffices to determine the number of eigenvalues of \eqref{aa4.1} that are smaller than $1$.

We start with some crucial estimates of the eigenvalues $\lambda_{\varepsilon,i}$ and the corresponding  eigenfunctions $\omega_{\varepsilon,i}$ of \eqref{aa4.1}.
Without loss of generality, we assume that $\|\omega_{\varepsilon,i}\|_*=1$.
Set $u_{\varepsilon,i}=\eta_{i}
u_{\varepsilon}$, $i=1,\cdots,k$. It is then evident that they are linearly independent.
\begin{proposition}\label{Prop4.2}
Let $u_\varepsilon$ be the $k$-peak solution to problem \eqref{1.1} with the form \eqref{1a}.
Then we have $$\lambda_{\varepsilon , i }<1,\;\; i=1,\cdots,k.$$
\end{proposition}
\begin{proof}[\bf Proof.]
Letting $S_{k}=span\big\{u_{\varepsilon,1},\cdots,
u_{\varepsilon,k}\big\}$, we see that  $\dim{S_{k}}=k$.
Then it follows from  the Courant-Weyl min-max variational principle that
\begin{align}\label{p098ju4.1}
\lambda_{\varepsilon , k}=&\inf_{\dim S=k}\sup_{u\in S\backslash \{0\}}J(u)
\leq\sup_{u\in S_{k}\backslash \{0\}} J(u),
\end{align}
where $$J(u)=\frac{\|u\|^{2}}
{(p-\varepsilon)\displaystyle\int_{\mathbb{R}^N}
u_{\varepsilon}^{p-1-\varepsilon}u^{2}\mathrm{d}x}.$$
Since $u\in S_{k}\backslash\{0\}$, then there exists some  $\big(\alpha_{\varepsilon ,1},\cdots,
\alpha_{\varepsilon ,k}\big)\in\mathbb{R}^k$ such that $u=\displaystyle\sum_{i=1}^{k}
\alpha_{\varepsilon ,i}u_{\varepsilon,i}$.
Thus
\begin{align}\label{ju4.1}
J(u)=\frac{\displaystyle\sum_{i=1}^{k}
\alpha_{\varepsilon ,i}^{2}
\int_{\mathbb{R}^N}\big(|\nabla u_{\varepsilon,i}|^{2}+V(x)
u_{\varepsilon,i}^{2}\big)
\mathrm{d}x}
{(p-\varepsilon)\displaystyle\sum_{i=1}^{k}
\alpha_{\varepsilon ,i}^{2}
\int_{\mathbb{R}^N}
u_{\varepsilon}^{p-1-\varepsilon}
u_{\varepsilon,i}^{2}\mathrm{d}x}\leq
\sup_{1\leq i\leq k}\frac{
\displaystyle\int_{\mathbb{R}^N}
\big(|\nabla u_{\varepsilon,i}|^{2}+V(x)
u_{\varepsilon,i}^{2}\big)
\mathrm{d}x}
{(p-\varepsilon)
\displaystyle\int_{\mathbb{R}^N}
u_{\varepsilon}^{p-1-\varepsilon}
u_{\varepsilon,i}^{2}\mathrm{d}x}:=\sup_{1\leq i\leq k}
\frac{A_{\varepsilon ,i}}{B_{\varepsilon ,i}}.
\end{align}

Next, we estimate $A_{\varepsilon ,i}$ and
$B_{\varepsilon ,i}$ separately.

Multiplying  \eqref{1.1} by $\eta_{i}^{2}
u_{\varepsilon}$ and integrating on $\mathbb{R}^N$, we obtain
\begin{align}\label{gl4.1}
\int_{\mathbb{R}^N}\eta_{i}^{2}
u_{\varepsilon}^{p+1-\varepsilon}
\mathrm{d}x
&=\int_{\mathbb{R}^N}\eta_{i}^{2}
 ( \nabla u_{\varepsilon})^{2}\mathrm{d}x
+\int_{\mathbb{R}^N}V(x)\eta_{i}^{2}u_{\varepsilon}^{2}
\mathrm{d}x
+2\int_{\mathbb{R}^N}\eta_{i}u_{\varepsilon}
\nabla\eta_{i} \nabla u_{\varepsilon}
\mathrm{d}x\notag
\\
&=\int_{\mathbb{R}^N}|\nabla(\eta_{i}
u_{\varepsilon})|^{2}\mathrm{d}x
-\int_{\mathbb{R}^N}|\nabla\eta_{i}|^{2}
u_{\varepsilon}^{2}\mathrm{d}x
+\int_{\mathbb{R}^N}V(x)\eta_{i}^{2}u_{\varepsilon}^{2}
\mathrm{d}x,
\end{align}
which gives that
\begin{align*}%\label{l4.1}
A_{\varepsilon ,i}&
=
\int_{\mathbb{R}^N}\eta_{i}^{2}
u_{\varepsilon}^{p+1-\varepsilon}\mathrm{d}x+\int_{\mathbb{R}^N}|\nabla\eta_{i}|^{2}
u_{\varepsilon}^{2}\mathrm{d}x%\notag
%\\
%&=\int_{ \mathbb{R}^N}
%\eta_{i}^{2}
%u_{\varepsilon}^{p+1-\varepsilon}\mathrm{d}x+\int_{B_{2\delta(\xi_{\varepsilon,i})\backslash B_\delta(\xi_{\varepsilon,i})}}|\nabla\eta_{i}|^{2}
%u_{\varepsilon}^{2}\mathrm{d}x
%\notag\\
%&
=\int_{ \mathbb{R}^N}
\eta_{i}^{2}
u_{\varepsilon}^{p+1-\varepsilon}\mathrm{d}x+O\Big(
\mu_{\varepsilon,i}^{2-N}\Big).
\end{align*}
By Proposition \ref{prop2.4} and Lemma \ref{lem4.1}, we get
\begin{align*}%\label{l4.1121}
\int_{\mathbb{R}^N}\eta_{i}^{2}
u_{\varepsilon}^{p+1-\varepsilon}
\mathrm{d}x
&=\int_{ B_\delta(\xi_{\varepsilon,i})}
\big(U_{\mu_{\varepsilon,i},\xi_{\varepsilon,i}}
+\phi_\varepsilon\big)^{p+1}
\mathrm{d}x+O\big(\varepsilon|\log\varepsilon|\big)+O\Big(
\mu_{\varepsilon,i}^{-N}\Big)%\notag
\\
&=\int_{ B_\delta(\xi_{\varepsilon,i})}
U_{\mu_{\varepsilon,i},\xi_{\varepsilon,i}}
^{p+1}\mathrm{d}x+O\big(
\|\phi_\varepsilon\|_*
\big)
+O\Big(
\|\phi_\varepsilon\|_*^{p+1}
\Big)+O\big(\varepsilon|\log
\varepsilon|\big)%\notag
\\
&=\int_{ \mathbb{R}^N}
U_{1,0}
^{p+1}\mathrm{d}x+O\Big(\varepsilon^{
1-\frac{\sigma}{2}}\Big).
\end{align*}
Then we deduce
\begin{align}\label{loppol4.1}
A_{\varepsilon ,i}
=\int_{ \mathbb{R}^N}
U_{1,0}
^{p+1}\mathrm{d}x+O\Big(\varepsilon^{
1-\frac{\sigma}{2}}\Big),
\end{align}
and
\begin{align}\label{ol4.1}
B_{\varepsilon ,i}&
%=(p-\varepsilon)
%\int_{ \mathbb{R}^N}
%U_{1,0}
%^{p+1}\mathrm{d}x+O\Big(\varepsilon^{
%1-\frac{\sigma}{2}}\Big)
=p\int_{ \mathbb{R}^N}
U_{1,0}
^{p+1}\mathrm{d}x+O\Big(\varepsilon^{
1-\frac{\sigma}{2}}\Big).
\end{align}
As a result,  from \eqref{p098ju4.1}, \eqref{ju4.1}, \eqref{loppol4.1} and \eqref{ol4.1}, we obtain
$$
J(u)\leq
\frac{1}{p}+o(1)<1.
$$
Hence we conclude that  $\lambda_{\varepsilon ,i}<1$ for $i=1,\cdots,k$ and the proof of Proposition \ref{Prop4.2} is complete.
\end{proof}

Next, we  analyze the eigenvalues $\lambda_{\varepsilon ,l}$  for $l=k+1,\cdots,k(N+1)$.
For this purpose, let $$S_{k(N+1)}=span\Big\{
\eta_{i}u_{\varepsilon},\, \eta_{i}\frac{\partial u_{\varepsilon}}
{\partial x_{j}},\, i=1,\cdots,k, \,j=1,\cdots,N
\Big\}.$$
Then $\dim S_{k(N+1)}=k(N+1)$.

\begin{proposition}\label{Propyu94.2}
Let $u_\varepsilon$ be the $k$-peak solution to problem \eqref{1.1} with the form \eqref{1a}. Then
for $l=k+1,\cdots,k(N+1)$, we have
\begin{align}\label{po0e2025}
\lambda_{\varepsilon,l}\leq 1+O(\varepsilon^{2 })\ \ \ \text{and} \ \
\
\lim_{\varepsilon\rightarrow0}
\lambda_{\varepsilon,l}
=1.
\end{align}
\end{proposition}
\begin{proof}[\bf Proof]
The proof is divided into  two steps.

\vskip 0.1cm

\noindent $\boldsymbol{Step \ 1.}$  We show  that
\begin{align}\label{e20257}
\lambda_{\varepsilon ,k(N+1)}\leq1+O(\varepsilon^{2 }).
\end{align}
By the Courant-Weyl min-max
variational principle, we have
\begin{align}\label{6.lp1}
\lambda_{\varepsilon ,k(N+1)}&=\inf_{\dim S=k(N+1)}\sup_{u\in S\backslash \{0\}}J(u)
\leq\sup_{u\in S_{k(N+1)}\backslash \{0\}}J(u).
\end{align}
Letting  $u\in S_{k(N+1)}\backslash \{0\}$, then there exist $c_{\varepsilon,i,j}$, $i=1,\cdots,k$, $j=1,\cdots,N$ such that
 $$u=\displaystyle\sum_{i=1}^{k}c_{\varepsilon,i,0}
u_{\varepsilon ,i}+
\displaystyle\sum_{i=1}^{k}\eta_{i}g_{\varepsilon,i},$$
where  $g_{\varepsilon,i}=\displaystyle\sum_{j=1}
^{N}c_{\varepsilon,i,j}\frac{\partial u_{\varepsilon}}
{\partial x_{j}}$.
Without loss of generality, we  can assume that
$$\displaystyle\sum_{i=1}^{k}\displaystyle\sum_{j=0}
^{N}c_{\varepsilon,i,j}^{2}=1.$$
By virtue of \eqref{gl4.1}, we have
\[
\int_{\mathbb{R}^N}
\left(\left|\nabla(\eta_{i}
u_{\varepsilon})\right|^{2}+V(x)(\eta_{i}
u_{\varepsilon})^{2}\right)\mathrm{d}x
=
\displaystyle\int_{\mathbb{R}^N}\eta_{i}^{2}
u_{\varepsilon}^{p+1-\varepsilon}
\mathrm{d}x
+
\displaystyle\int_{\mathbb{R}^N}|\nabla\eta_{i}|^{2}
u_{\varepsilon}^{2}\mathrm{d}x.
\]
Note that $g_{\varepsilon,i}$ satisfies
\begin{align}\label{mpp2}
-\Delta g_{\varepsilon ,i}+V(x)g_{\varepsilon ,i}+
\displaystyle\sum_{j=1}
^{N}c_{\varepsilon ,i,j}\frac{\partial V(x)}
{\partial x_{j}}u_{\varepsilon}
=(p-\varepsilon)
u_{\varepsilon}^{p-1-\varepsilon}g_{\varepsilon ,i}.
\end{align}
Multiplying \eqref{mpp2} by  $\eta_{i}^{2}g_{\varepsilon ,i}$ and integrating on $\mathbb{R}^N$, we  get
\begin{align*}
&(p-\varepsilon)\int_{\mathbb{R}^N}
\eta_{i}^{2}
u_{\varepsilon}^{p-1-\varepsilon}
g_{\varepsilon ,i}^{2}
\mathrm{d}x\\
%=
%\int_{\mathbb{R}^N}\nabla g_{\varepsilon ,i}\nabla (\eta_{i}^{2} g_{\varepsilon ,i})
%+V(x)\eta_{i}^{2}g_{\varepsilon ,i}^{2}
%+
%\displaystyle\sum_{j=1}
%^{N}c_{\varepsilon ,i,j}\frac{\partial V(x)}
%{\partial x_{j}}u_{\varepsilon}
%\eta_{i}^{2}g_{\varepsilon ,i}^{2}
%\mathrm{d}x
%\\
=&\int_{\mathbb{R}^N}\left(\eta_{i}^{2}
 |\nabla g_{\varepsilon ,i}|^{2}
+V(x)\eta_{i}^{2}g_{\varepsilon ,i}^{2}
+2\eta_{i}\nabla\eta_{i} g_{\varepsilon ,i}\nabla g_{\varepsilon ,i}\right)\mathrm{d}x+
\displaystyle\sum_{j=1}
^{N}c_{\varepsilon ,i,j}\int_{\mathbb{R}^N}
\frac{\partial V(x)}
{\partial x_{j}}u_{\varepsilon}
\eta_{i}^{2}g_{\varepsilon ,i}^{2}
\mathrm{d}x
\\
=&\int_{\mathbb{R}^N}\left(|\nabla(\eta_{i}
g_{\varepsilon ,i})|^{2}-|\nabla\eta_{i}|^{2}
g_{\varepsilon ,i}^{2}
+V(x)\eta_{i}^{2}g_{\varepsilon ,i}^{2}\right)
\mathrm{d}x+
\displaystyle\sum_{j=1}
^{N}c_{\varepsilon ,i,j}\int_{\mathbb{R}^N}
\frac{\partial V(x)}
{\partial x_{j}}u_{\varepsilon}
\eta_{i}^{2}g_{\varepsilon ,i}
\mathrm{d}x.
\end{align*}
Thus we derive
\begin{align*}
&\int_{\mathbb{R}^N}
\left(\left|\nabla(\eta_{i}
g_{\varepsilon,i})\right|^{2}
+V(x)(\eta_{i}
g_{\varepsilon,i})^{2}\right)\mathrm{d}x
\\
=&
(p-\varepsilon)
\int_{\mathbb{R}^N}
\eta_{i}^{2}u_{\varepsilon}^{p-1-\varepsilon}
g_{\varepsilon ,i}^{2}\mathrm{d}x
+
\int_{\mathbb{R}^N}|\nabla\eta_{i}|^{2}
g_{\varepsilon ,i}^{2}
\mathrm{d}x-
\displaystyle\sum_{j=1}
^{N}c_{\varepsilon ,i,j}
\int_{\mathbb{R}^N}\frac{\partial V(x)}
{\partial x_{j}}u_{\varepsilon}g_{\varepsilon ,i}\eta_{i}^{2}
\mathrm{d}x
.
\end{align*}
Similarly,
\[
\begin{split}
&2\int_{\mathbb{R}^N}
\left(\nabla(\eta_{i}u_{\varepsilon})\nabla(\eta_{i}
g_{\varepsilon,i})+V(x)\eta_{i}^{2}
u_{\varepsilon}g_{\varepsilon,i}\right)
\mathrm{d}x
\\
=&(p+1-\varepsilon)\int_{\mathbb{R}^N}
\eta_{i}^{2}
u_{\varepsilon}^{p-\varepsilon}g_{\varepsilon ,i}\mathrm{d}x
+2\int_{\mathbb{R}^N}
|\nabla\eta_{i}|^{2}u_{\varepsilon}
g_{\varepsilon ,i}\mathrm{d}x
-
\displaystyle\sum_{j=1}
^{N}c_{\varepsilon ,i,j}
\int_{\mathbb{R}^N}\frac{\partial V(x)}
{\partial x_{j}}u_{\varepsilon}^{2}\eta_{i}^{2}
\mathrm{d}x.
\end{split}
\]
Hence we have
\begin{align*}%\label{123pl}
\|u\|^{2}&=\displaystyle\sum_{i=1}
^{k}\displaystyle\int_{\mathbb{R}^N}
\left(\left|\nabla\big(\eta_{i}(c_{\varepsilon ,i,0}
u_{\varepsilon}+g_{\varepsilon ,i}
)
\big)\right|^{2}
+V(x)\eta_{i}^{2}
\big(c_{\varepsilon,i,0}
u_{\varepsilon}+g_{\varepsilon ,i}
\big)^{2}\right)
\mathrm{d}x\\&=\displaystyle\sum_{i=1}
^{k}c_{\varepsilon,i,0}^{2}\int_{\mathbb{R}^N}
\left(\left|\nabla(\eta_{i}
u_{\varepsilon})\right|^{2}+V(x)(\eta_{i}
u_{\varepsilon})^{2}\right)\mathrm{d}x+
\displaystyle\sum_{i=1}
^{k}\int_{\mathbb{R}^N}
\left(\left|\nabla(\eta_{i}
g_{\varepsilon,i})\right|^{2}
+V(x)(\eta_{i}
g_{\varepsilon,i})^{2}\right)\mathrm{d}x
\\&
\quad+2\displaystyle\sum_{i=1}
^{k}c_{\varepsilon,i,0}\int_{\mathbb{R}^N}
\left(\nabla(\eta_{i}u_{\varepsilon})\nabla(\eta_{i}
g_{\varepsilon,i})+V(x)\eta_{i}^{2}
u_{\varepsilon}g_{\varepsilon,i}\right)
\mathrm{d}x\\:&=\Upsilon_{1}+\Upsilon_{2}+\Upsilon_{3}.
\end{align*}
 By %$\mu_{\varepsilon,i}\sim\varepsilon^{-\frac{1}{2}}$,
Lemma \ref{lem4.1}, we  get
\begin{align*}
\Upsilon_{1}&=\displaystyle\sum_{i=1}
^{k}c_{\varepsilon,i,0}^{2}
\int_{B_{\delta}(\xi_{\varepsilon,i})}
u_{\varepsilon}^{p+1-\varepsilon}
\mathrm{d}x+
O(\varepsilon^{\frac{N-2}{2}})
\\&=(p-\varepsilon)\displaystyle\sum_{i=1}
^{k}c_{\varepsilon ,i,0}^{2}
\int_{B_{\delta}(\xi_{\varepsilon,i})}
u_{\varepsilon}^{p+1-\varepsilon}
\mathrm{d}x
+
\big(1-p\big)\displaystyle\sum_{i=1}
^{k}c_{\varepsilon ,i,0}^{2}
\int_{B_{\delta}(\xi_{\varepsilon,i})}
u_{\varepsilon}^{p+1-\varepsilon}
\mathrm{d}x
+O(\varepsilon)
,
\end{align*}
\begin{align*}
\Upsilon_{2}=(p-\varepsilon)\displaystyle\sum_{i=1}
^{k}
\int_{B_{\delta}(\xi_{\varepsilon,i})}
u_{\varepsilon}^{p-1-\varepsilon}
g_{\varepsilon ,i}^{2}
\mathrm{d}x+O(\varepsilon)
\end{align*}
and
\[
\Upsilon_{3}=O(\varepsilon).
\]
%\begin{align*}
%\Upsilon_{3}&=O\Big(\displaystyle\sum_{i=1}
%^{k}\displaystyle\sum_{j=1}
%^{N}|c_{\varepsilon ,i,j}|\cdot|c_{\varepsilon ,i,0}|
%\mu_{\varepsilon,i}^{-N }+
%\displaystyle\sum_{i=1}
%^{k}\displaystyle\sum_{j=1}
%^{N}|c_{\varepsilon ,i,j}|\cdot|c_{\varepsilon ,i,0}|
%\mu_{\varepsilon,i}^{-2 }+
%\displaystyle\sum_{i=1}
%^{k}\displaystyle\sum_{j=1}
%^{N}|c_{\varepsilon ,i,j}|\cdot|c_{\varepsilon ,i,0}|
%\mu_{\varepsilon,i}^{2-N }\Big)
%\\&=O\Big(\displaystyle\sum_{i=1}
%^{k}\displaystyle\sum_{j=1}
%^{N}|c_{\varepsilon ,i,j}|\cdot|c_{\varepsilon ,i,0}|
%\mu_{\varepsilon,i}^{-2 }\Big).
%\end{align*}
As a result,  we conclude that
\begin{align}\label{a6.163}
\|u\|^{2}=(p-\varepsilon)\displaystyle\sum_{i=1}
^{k}
\int_{B_{\delta}(\xi_{\varepsilon,i})}
u_{\varepsilon}^{p-1-\varepsilon}
\Big(
c_{\varepsilon ,i,0}^{2}
u_{\varepsilon}^{2}+g_{\varepsilon,i}^{2}
\Big)
\mathrm{d}x
+\big(1-p\big)\displaystyle\sum_{i=1}
^{k}c_{\varepsilon ,i,0}^{2}
\int_{B_{\delta}(\xi_{\varepsilon,i})}
u_{\varepsilon}^{p+1-\varepsilon}
\mathrm{d}x+O(\varepsilon)
.
\end{align}
On the other hand, by a direct calculation,  we have
\begin{align*}
(p-\varepsilon)\int_{\mathbb{R}^N}
u_{\varepsilon}^{p-1-\varepsilon}u^{2}\mathrm{d}x
&=(p-\varepsilon)\displaystyle\sum_{i=1}
^{k}
\int_{B_{\delta}(\xi_{\varepsilon,i})}
u_{\varepsilon}^{p-1-\varepsilon}
\Big(
c_{\varepsilon ,i,0}^{2}
u_{\varepsilon}^{2}+g_{\varepsilon ,i}^{2}
\Big)
\mathrm{d}x+Q_{1},
\end{align*}
where
\[
Q_{1}=(p-\varepsilon)\displaystyle\sum_{i=1}
^{k}\bigg(2c_{\varepsilon,i,0}
\int_{B_{2\delta}(\xi_{\varepsilon,i})}
u_{\varepsilon}^{p-\varepsilon}\eta_{i}^{2}
g_{\varepsilon,i}
\mathrm{d}x+
\int_{B_{2\delta}(\xi_{\varepsilon,i})\backslash
B_{\delta}(\xi_{\varepsilon,i})}
\eta_{i}^{2}u_{\varepsilon}^{p-1-\varepsilon}
\Big(
c_{\varepsilon ,i,0}^{2}
u_{\varepsilon}^{2}+g_{\varepsilon ,i}^{2}
\Big)
\mathrm{d}x\bigg).
\]
By virtue of Lemma \ref{lem4.1},  we  derive
\begin{align*}
Q_{1}&=2(p-\varepsilon)\displaystyle\sum_{i=1}
^{k}\displaystyle\sum_{j=1}
^{N}c_{\varepsilon ,i,0}c_{\varepsilon ,i,j}
\int_{B_{\delta}(\xi_{\varepsilon,i})} u_{\varepsilon}^{p-\varepsilon}
\frac{\partial u_{\varepsilon}}
{\partial x_{j}}
\mathrm{d}x+O(\varepsilon^{\frac{N}{2}})=O(\varepsilon^{\frac{N}{2}}).
\end{align*}
Thus we infer
\begin{align}\label{apop6.163}
(p-\varepsilon)\int_{\mathbb{R}^N}
u_{\varepsilon}^{p-1-\varepsilon}u^{2}
\mathrm{d}x
=(p-\varepsilon)\displaystyle\sum_{i=1}
^{k}
\int_{B_{\delta}(\xi_{\varepsilon,i})}
u_{\varepsilon}^{p-1-\varepsilon}
\Big(
c_{\varepsilon ,i,0}^{2}
u_{\varepsilon}^{2}+g_{\varepsilon ,i}^{2}
\Big)
\mathrm{d}x+O(\varepsilon^{\frac{N}{2}}).
\end{align}
It follows from  \eqref{a6.163} and \eqref{apop6.163} that
\begin{align}\label{0ip.lp1}
J(u)=1+\frac{\big(1-p\big)\displaystyle\sum_{i=1}
^{k}c_{\varepsilon ,i,0}^{2}
\int_{B_{\delta}(\xi_{\varepsilon,i})}
u_{\varepsilon}^{p+1-\varepsilon}
\mathrm{d}x+O(\varepsilon)
}
{(p-\varepsilon)\displaystyle\sum_{i=1}
^{k}
\int_{B_{\delta}(\xi_{\varepsilon,i})}
u_{\varepsilon}^{p-1-\varepsilon}
\Big(
c_{\varepsilon ,i,0}^{2}
u_{\varepsilon}^{2}+g_{\varepsilon,i}^{2}
\Big)
\mathrm{d}x+O(\varepsilon^{\frac{N}{2}})}.
\end{align}
Choosing $c_{\varepsilon, i,0}=0$, $c_{\varepsilon,i,j}=\frac{1}{\sqrt{kN}}$, $j=1,\cdots,N$, $i=1,\cdots,k$ in \eqref{0ip.lp1}, we obtain
\begin{align}\label{lp8u1}
\sup_{u\in S_{k(N+1)}\backslash \{0\}}J(u)\geq 1+ O(\varepsilon^{2 }).
\end{align}
Let $d_{\varepsilon,i,j}$, $i=1,\cdots,k$, $j=0,1,\cdots,N$ be the  maximum point of $J(u)$ on $S_{k(N+1)}$. Then
by direct computations, we have
\begin{align}\label{ip.lp1}
\sup_{u\in S_{k(N+1)}\backslash \{0\}}J(u)&\leq1+
\sup_{1\leq i\leq k}
\frac{\big(1-p\big)
d_{\varepsilon,i,0}^{2}
\displaystyle\int_{B_{\delta}(\xi_{\varepsilon,i})}
u_{\varepsilon}^{p+1-\varepsilon}
\mathrm{d}x+O(\varepsilon)}
{(p-\varepsilon)\displaystyle
\int_{B_{\delta}(\xi_{\varepsilon,i})}
u_{\varepsilon}^{p-1-\varepsilon}
\Big(
d_{\varepsilon,i,0}^{2}
u_{\varepsilon}^{2}+\widehat{g}_{\varepsilon,i}^{2}
\Big)
\mathrm{d}x+O\big(\displaystyle\varepsilon^{\frac{N}{2}}
\big)}\notag
\\
:&=1+J_{1}(u),
\end{align}
where $\widehat{g}_{\varepsilon ,i}=\displaystyle\sum_{j=1}
^{N}d_{\varepsilon ,i,j}\frac{\partial u_{\varepsilon}}
{\partial x_{j}}$.
Thus by Lemma \ref{lem4.1},
there exists $i_{0}\in\{1,\cdots,k\}$  such that \begin{align*}
J_{1}(u)&=
\frac{\big(1-p\big)d_{\varepsilon ,i_{0},0}^{2}
\displaystyle\int_{B_{\delta}
(\xi_{\varepsilon,i_{0}})}
u_{\varepsilon}^{p+1-\varepsilon}
\mathrm{d}x+O(\varepsilon)}
{(p-\varepsilon)\displaystyle
\int_{B_{\delta}(\xi_{\varepsilon,i_{0}})}
u_{\varepsilon}^{p-1-\varepsilon}
\Big(
d_{\varepsilon,i_{0},0}^{2}
u_{\varepsilon}^{2}+\widehat{g}_{\varepsilon
,i_{0}}^{2}
\Big)
\mathrm{d}x+O\Big(
\varepsilon^{\frac{N}{2} }\Big)}\\
&=
\frac{\big(1-p\big)(1+o(1))d_{\varepsilon ,i_{0},0}^{2}
\displaystyle\int_{B_{\delta}
(\xi_{\varepsilon,i_{0}})}
U_{\mu_{\varepsilon,i},\xi_{\varepsilon,i}}^{p+1
}
\mathrm{d}x}
{(p-\varepsilon)(1+o(1))
\displaystyle
\int_{B_{\delta}(\xi_{\varepsilon,i_{0}})}
u_{\varepsilon}^{p-1-\varepsilon}
\Big(
d_{\varepsilon,i_{0},0}^{2}
u_{\varepsilon}^{2}
+\displaystyle\sum_{j=1}
^{N}d_{\varepsilon ,i,j}^{2}\Big(\frac{\partial u_{\varepsilon}}
{\partial x_{j}}\Big)^{2}
\Big)
\mathrm{d}x}\\
&\leq\frac{\big(1-p\big)(1+o(1))
d_{\varepsilon ,i_{0},0}^{2}
\mu_{\varepsilon,i_{0}}^{2}
\displaystyle\int_{\mathbb{R}^N}
U_{1,0}^{p+1}
\mathrm{d}x}
{C(1+o(1))\mu_{\varepsilon,i_{0}}^{4}\Big(d_{\varepsilon,i_{0},0}^{2}
\mu_{\varepsilon,i_{0}}^{-2}
+ \displaystyle\sum_{j=1}
^{N}d_{\varepsilon ,i,j}^{2}\Big)},
\end{align*}
where $C>0$ is a constant independent of $\varepsilon$.

We claim that
\begin{align}\label{lp8u2}
d_{\varepsilon ,i_{0},0}^{2}
\mu_{\varepsilon,i_{0}}^{2}\leq
K\displaystyle\sum_{j=1}
^{N}d_{\varepsilon ,i_{0},j}^{2},
\end{align}
where $K$ is a positive constant independent of $\varepsilon$. Assume, on the contrary, that  for arbitrary large positive constant $\tilde{K}$, there are
$d_{\varepsilon ,i_{0},j}$, $j=0,1,\cdots,N$ satisfying
\begin{align*}
d_{\varepsilon ,i_{0},0}^{2}
\mu_{\varepsilon,i_{0}}^{2}>
\tilde{K}\displaystyle\sum_{j=1}
^{N}d_{\varepsilon ,i_{0},j}^{2}
.
\end{align*}
Then we derive
\begin{align*}
J_{1}(u)\leq
\frac{\big(1-p\big)(1+o(1))
\displaystyle\int_{\mathbb{R}^N}
U_{1,0}^{p+1}
\mathrm{d}x}
{C(1+o(1))\mu_{\varepsilon,i_{0}}^{4}
\Big(
\mu_{\varepsilon,i_{0}}^{-4}
+ d_{\varepsilon ,i_{0},0}^{-2}
\mu_{\varepsilon,i_{0}}^{-2}\displaystyle\sum_{j=1}
^{N}d_{\varepsilon ,i,j}^{2}\Big)}\leq-\tilde{K}_1\varepsilon^{2},
\end{align*}
 where $\tilde{K}_{1}$ is  a arbitrary large positive constant.
This is a contradiction with \eqref{lp8u1} and our claim follows.
As a consequence, in view  of  \eqref{6.lp1}, \eqref{lp8u1}, \eqref{ip.lp1} and \eqref{lp8u2},  we obtain  \eqref{e20257}.

\vskip 0.1cm

\noindent $\boldsymbol{Step \ 2.}$  We claim that
 \begin{align}\label{e20251}
\lambda_{\varepsilon ,k+1}\geq1+o(1).
\end{align}
%Without loss of generality, we let $\|\omega_{k+1,\varepsilon}\|_*=1$.
On the contrary, we assume that $
\lambda_{k+1}<1$, where $\lambda_{k+1}=\displaystyle
\lim_{\varepsilon\rightarrow0}
\lambda_{\varepsilon ,k+1}$.
Note that $\omega_{\varepsilon ,k+1}$  satisfies $\|\omega_{\varepsilon ,k+1}\|_*=1$ and
$$
-\Delta \omega_{\varepsilon ,k+1}+V(x)
\omega_{\varepsilon ,k+1}=
(p-\varepsilon)\lambda_{\varepsilon ,k+1} u_\varepsilon^{p-1-\varepsilon}
\omega_{\varepsilon ,k+1}\;\;\text{in}\;\;\mathbb{R}^N.
$$
Using Lemmas \ref{lem: the computation} and \ref{lem4.1}, we have
\begin{equation}\label{aa4.18}
|\omega_{\varepsilon ,k+1}|\leq C\int_{\mathbb{R}^N}\frac{1}{|x-y|^{N-2}} u_\varepsilon^{p-1-\varepsilon}(y)
\omega_{\varepsilon ,k+1}(y)\mathrm{d}y
\leq C\|\omega_{\varepsilon ,k+1}\|_*
\displaystyle\sum_{i=1}^{k}\frac{\mu_{\varepsilon,i}^{\frac{N-2}{2}}}
{(1+\mu_{\varepsilon,i}^{2}|x-\xi_{\varepsilon,i}|^{2})^{\frac{\beta}{2}}},
\end{equation}
where $\beta\in(\frac{N-2}{2}+\sigma,N-2)$. Then we find
$$
\bigg(
\displaystyle\sum_{i=1}^{k}
\frac{\mu_{\varepsilon,i}^{\frac{N-2}{2}}}
{\big(1+\mu_{\varepsilon,i}^{2}|x-\xi_{\varepsilon,i}|^{2}\big)
^{\frac{N-2}{4}+\frac{\sigma}{2}}}\bigg)^{-1}
|\omega_{\varepsilon ,k+1}(x)|
\leq  C
\displaystyle\sum_{i=1}^{k}\frac{\|\omega_{
\varepsilon ,k+1}\|_*}
{(1+\mu_{\varepsilon,i}^{2}|x-\xi_{\varepsilon,i}|^{2})
^{\frac{\beta}
{2}-\frac{N-2}{4}-\frac{\sigma}{2}}},
$$
which yields that
there exists some $i_{0}\in\{1,\cdots,k\}$ such that $\bigg(
\displaystyle\sum_{i=1}^{k}
\frac{\mu_{\varepsilon,i}^{\frac{N-2}{2}}}
{\big(1+\mu_{\varepsilon,i}^{2}|x-\xi_{\varepsilon,i}|^{2}\big)
^{\frac{N-2}{4}+\frac{\sigma}{2}}}\bigg)^{-1}
|\omega_{\varepsilon ,k+1}(x)|$
 attains its maximum in $ B_{R\mu_{i_{0},\varepsilon}^{-1}}(\xi_{i_{0}})$ for some large constant $R>0$.

Let $\tilde{\omega}_{\varepsilon,k+1,i_0}=
\mu_{\varepsilon,i_{0}}
^{-\frac{N-2}{2}}\omega_{\varepsilon,k+1}
(\xi_{\varepsilon,i_{0}}+x/\mu_{\varepsilon,i_{0}})$ . Then we have
\begin{align*}
-\Delta \tilde{\omega}_{\varepsilon,k+1,i_0}+
\mu_{\varepsilon,i_{0}}
^{-2}V(\xi_{\varepsilon,i_{0}}+x/\mu_{\varepsilon,i_{0}})
\tilde{\omega}_{\varepsilon,k+1,i_0}
=(p-\varepsilon)\lambda_{\varepsilon,k+1}
\mu_{\varepsilon ,i_{0}}^{-2} u_\varepsilon^{p-1-\varepsilon}(\xi_{\varepsilon,i_{0}}+x/\mu_{\varepsilon,i_{0}})
\tilde{\omega}_{\varepsilon,k+1,i_0}\;\;\text{in}\;\;
\mathbb{R}^N.
\end{align*}
By the regularity theory of elliptic equation,  we know that there is $\omega_{k+1,i_0}^*\neq0$ such that
$$\displaystyle
\tilde{\omega}_{\varepsilon,k+1,i_0}
\rightarrow\omega_{k+1,i_0}^*\;\;\text{in}\;\; C^{2}_{loc}(\mathbb{R}^N)$$
and
\begin{align*}
-\Delta \omega_{k+1,i_0}^*=p\lambda_{k+1} U_{1,0}^{p-1}
\omega_{k+1,i_0}^*\;\;\text{in}\;\;\mathbb{R}^N.
\end{align*}
Since $\lambda_{k+1}<1$, by Lemma A.1 in \cite{CGB}, we get $\lambda_{k+1}=
\frac{N-2}{N+2}$ and $\omega_{k+1,i_0}^*=\Upsilon U_{1,0}$, where $\Upsilon\neq0$. On the other hand, we have
\begin{equation*}
\int_{\mathbb{R}^N} u_{\varepsilon}^{p-\varepsilon}(y)
\omega_{\varepsilon,k+1}(y)\mathrm{d}y=0.
\end{equation*}
Letting $\varepsilon\rightarrow0$, we obtain
\begin{equation*}
\int_{\mathbb{R}^N} U_{1,0}^{p+1}(y)\mathrm{d}y=0,
\end{equation*}
which leads to a contradiction with  $\Upsilon\neq0$.  Therefore,  $\lambda_{k+1}\geq1$ and the desired estimates follows from  \eqref{e20257} and \eqref{e20251}.
\end{proof}

 Now we proceed to estimate the eigenvalues $\lambda_{\varepsilon,l}$ for $l=k(N+1)+1,\cdots,k(N+2)$.
\begin{proposition}\label{Prop0094.2}
For $l=k(N+1)+1,\cdots,k(N+2)$, there holds
\begin{align*}%\label{7e2025}
\lambda_{\varepsilon, l}\rightarrow1\ \ \text{as} \ \ \varepsilon\rightarrow0.
\end{align*}
\end{proposition}
\begin{proof}[\bf Proof.]
Let $\lambda_{l}=\lim\limits_{\varepsilon\rightarrow0}\lambda_{\varepsilon, l}$ for $ l=k(N+1)+1,\cdots,k(N+2)$.
 By Proposition \ref{Propyu94.2}, we have
 \begin{align*}
\lambda_{l}\geq1.
\end{align*}
Then we only need to prove
\begin{align}\label{e2025}
\lambda_{k(N+2)}\leq1.
\end{align}
Define
$$S_{k(N+2)}=span\Big\{
\eta_{i}u_{\varepsilon},\ \, \eta_{i}\frac{\partial u_{\varepsilon}}
{\partial x_{j}},\ \ \eta_{i}\varphi_{\varepsilon ,i},  \ \ j=1,\cdots,N,\ \  i=1,\cdots,k
\Big\},$$
where $\varphi_{\varepsilon ,i}=
(x-\xi_{\varepsilon,i})\nabla u_{\varepsilon}+\frac{2}{p-1-\varepsilon}
u_{\varepsilon}$.
Then $\dim S_{k(N+2)}=k(N+2)$. By the Courant-Weyl min-max principle, we have
\begin{align}\label{6.161}
\lambda_{\varepsilon ,k(N+2)}&=\inf_{\dim S=k(N+2)}\sup_{u\in S \backslash \{0\}}J(u)
\leq\sup_{u\in S_{k(N+2)}\backslash \{0\}}J(u).
\end{align}
Let $u\in S_{k(N+2)}\backslash \{0\}$. Then there exist  $\beta_{\varepsilon,i,j}$, $i=1,\cdots,k$, $j=0,1,\cdots,N+1$ such that  \[
u=\sum_{i=1}^k\left(\beta_{\varepsilon,i,0,}\eta_{i}
u_{\varepsilon}+\eta_{i}h_{\varepsilon ,i}+
\beta_{\varepsilon,i,N+1}\eta_{i}
\varphi_{\varepsilon,i}\right),
\]
where $h_{\varepsilon,i}=\displaystyle\sum_{j=1}
^{N}\beta_{\varepsilon,i,j}\frac{\partial u_{\varepsilon}}
{\partial x_{j}}$.
 Without loss of generality, we  assume that
$$\displaystyle\sum_{i=1}^{k}\displaystyle\sum_{j=0}
^{N+1}\beta_{\varepsilon,i,j}^{2}=1.$$\\
On the one hand, by Proposition \ref{Propyu94.2}, we  obtain
\begin{align*}%\label{123p2}
\|u\|^{2}&=
\displaystyle\sum_{i=1}
^{k}\displaystyle\int_{\mathbb{R}^N}\bigg
(\left|\nabla\big(\eta_{i}(\beta_{\varepsilon ,i,0,
}
u_{\varepsilon}+h_{\varepsilon ,i}+
\beta_{\varepsilon ,i,N+1}
\varphi_{\varepsilon ,i}
)
\big)\right|^{2}
+V(x)\eta_{i}^{2}
\big(\beta_{\varepsilon ,i,0}
u_{\varepsilon}+h_{\varepsilon ,i}+
\beta_{\varepsilon ,i,N+1}
\varphi_{\varepsilon ,i}
\big)^{2}\bigg)
\mathrm{d}x
\\&=\Upsilon_{0}+\displaystyle\sum_{i=1}
^{k}\beta_{\varepsilon ,i,N+1}^{2}\int_{\mathbb{R}^N}\Big(
\left|\nabla(\eta_{i}
\varphi_{\varepsilon ,i})\right|^{2}+V(x)(\eta_{i}
\varphi_{\varepsilon ,i})^{2}\Big)\mathrm{d}x
\\&
\quad+2\displaystyle\sum_{i=1}
^{k}\beta_{\varepsilon ,i,0}\beta_{\varepsilon ,i,N+1}
\int_{\mathbb{R}^N}\Big(
\nabla(\eta_{i}u_{\varepsilon})\nabla(\eta_{i}
\varphi_{\varepsilon,i})+V(x)\eta_{i}^{2}u_{\varepsilon}
\varphi_{\varepsilon,i}\Big)
\mathrm{d}x\\&
\quad+2\displaystyle\sum_{i=1}
^{k}\beta_{\varepsilon,i,N+1}\int_{\mathbb{R}^N}\Big(
\nabla(\eta_{i}\varphi_{\varepsilon ,i})
\nabla(\eta_{i}
h_{\varepsilon ,i})+V(x)\eta_{i}^{2}
\varphi_{\varepsilon ,i}h_{\varepsilon ,i}\Big)
\mathrm{d}x\notag\\:&=\Upsilon_{0}
+\Upsilon_{4}+\Upsilon_{5}+\Upsilon_{6},
\end{align*}
where (see \eqref{a6.163})
 \begin{align*}\Upsilon_{0}:&=(p-\varepsilon)
 \displaystyle\sum_{i=1}
^{k}
\int_{B_{\delta}(\xi_{\varepsilon,i})}
u_{\varepsilon}^{p-1-\varepsilon}\big(
\beta_{\varepsilon ,i,0}^{2}u_{\varepsilon}^{2}+h_{\varepsilon ,i}^{2}\big)
\mathrm{d}x
+(1-p)\displaystyle\sum_{i=1}
^{k}\beta_{\varepsilon ,i,0}^{2}
\int_{B_{\delta}(\xi_{\varepsilon,i})}
u_{\varepsilon}^{p+1-\varepsilon}
\mathrm{d}x
+O(\varepsilon).\end{align*}
 By a direct calculation,  we  deduce
\begin{align*}
\Upsilon_{4}&=(p-\varepsilon)
\displaystyle\sum_{i=1}
^{k}\beta_{\varepsilon ,i,N+1}^{2}
\int_{\mathbb{R}^N}
\eta_{i}^{2}u_{\varepsilon}^{p-1-\varepsilon}
\varphi_{\varepsilon ,i}^{2}\mathrm{d}x
+
\displaystyle\sum_{i=1}
^{k}\beta_{\varepsilon ,i,N+1}^{2}
\int_{\mathbb{R}^N}|\nabla\eta_{i}|^{2}
\varphi_{\varepsilon ,i}^{2}\mathrm{d}x
\\&
\quad-\displaystyle\sum_{i=1}
^{k}\beta_{\varepsilon ,i,N+1}^{2}
\int_{\mathbb{R}^N}\eta_{i}^{2}u_{\varepsilon}
\varphi_{\varepsilon ,i}
\Big(\big\langle
(x-\xi_{\varepsilon,i})\cdot\nabla V(x)\big\rangle+2V(x)\Big)\mathrm{d}x,
\\&
=(p-\varepsilon)\displaystyle\sum_{i=1}
^{k}\beta_{\varepsilon ,i,N+1}^{2}
\int_{B_{\delta}(\xi_{\varepsilon,i})}
u_{\varepsilon}^{p-1-\varepsilon}
\varphi_{\varepsilon ,i}^{2}
\mathrm{d}x
-2\displaystyle\sum_{i=1}
^{k}\beta_{\varepsilon ,i,N+1}^{2}
\int_{B_{\delta}(\xi_{\varepsilon,i})}
V(x)u_{\varepsilon}
\varphi_{\varepsilon ,i}\mathrm{d}x
\\&
\quad-\displaystyle\sum_{i=1}
^{k}\beta_{\varepsilon ,i,N+1}^{2}
\int_{B_{\delta}(\xi_{\varepsilon,i})}u_{\varepsilon}
\varphi_{\varepsilon ,i}\big\langle
(x-\xi_{\varepsilon,i})\cdot\nabla V(x)\big\rangle\mathrm{d}x+
O\left(\varepsilon^{\frac{N-2}{2}}\right)
\\&=(p-\varepsilon)\displaystyle\sum_{i=1}
^{k}\beta_{\varepsilon ,i,N+1}^{2}
\int_{B_{\delta}(\xi_{\varepsilon,i})}
u_{\varepsilon}^{p-1-\varepsilon}
\varphi_{\varepsilon ,i}^{2}
\mathrm{d}x+
O(\varepsilon),
\end{align*}
and
\begin{align*}
\Upsilon_{5}&=(p+1-\varepsilon)\displaystyle\sum_{i=1}
^{k}\beta_{\varepsilon ,i,0}\beta_{\varepsilon ,i,N+1}
\int_{\mathbb{R}^N}\eta_{i}^{2}
u_{\varepsilon}^{p-\varepsilon}
\varphi_{\varepsilon ,i}\mathrm{d}x+2\displaystyle\sum_{i=1}
^{k}\beta_{\varepsilon ,i,0}\beta_{\varepsilon ,i,N+1}
\int_{\mathbb{R}^N}|\nabla\eta_{i}|^{2}
u_{\varepsilon}
\varphi_{\varepsilon ,i}\mathrm{d}x
\\&
\quad-2\displaystyle\sum_{i=1}
^{k}\beta_{\varepsilon ,i,0}\beta_{\varepsilon ,i,N+1}
\int_{\mathbb{R}^N}V(x)\eta_{i}^{2}
u_{\varepsilon}^{2}\mathrm{d}x-\displaystyle\sum_{i=1}
^{k}\beta_{\varepsilon ,i,0}\beta_{\varepsilon,i,N+1}
\int_{\mathbb{R}^N}\big\langle
(x-\xi_{\varepsilon,i})\cdot\nabla V(x)\big\rangle\eta_{i}^{2}
u_{\varepsilon}^{2}\mathrm{d}x
\\&=(p+1-\varepsilon)\displaystyle\sum_{i=1}
^{k}\beta_{\varepsilon ,i,0}\beta_{\varepsilon ,i,N+1}
\displaystyle\int_{B_{\delta}
(\xi_{\varepsilon,i})}\Bigg(
\frac{\Big\langle(x-\xi_{\varepsilon,i})\cdot
\nabla  u_{\varepsilon}^{p+1-\varepsilon}
\Big\rangle}{p+1-\varepsilon}
+\frac{2u_{\varepsilon}^{p+1-\varepsilon}}{p-1-\varepsilon}
\Bigg)\mathrm{d}x+O(\varepsilon)\\
&=O(\varepsilon).
\end{align*}
By Proposition \ref{prop2.4} and Lemma \ref{lem4.1}, we get
\begin{align*}
\Upsilon_{6}&=2(p-\varepsilon)\displaystyle\sum_{i=1}
^{k}\beta_{\varepsilon ,i,N+1}\int_{\mathbb{R}^N}
\eta_{i}^{2}u_{\varepsilon}^{p-1-\varepsilon}
\varphi_{\varepsilon ,i}h_{\varepsilon ,i}\mathrm{d}x+
2\displaystyle\sum_{i=1}
^{k}\beta_{\varepsilon ,i,N+1}\int_{\mathbb{R}^N}
|\nabla\eta_{i}|^{2}\varphi_{\varepsilon ,i}
h_{\varepsilon ,i}\mathrm{d}x.
\\&\quad-2\displaystyle\sum_{i=1}
^{k}\beta_{\varepsilon ,i,N+1}\int_{\mathbb{R}^N}
V(x)\eta_{i}^{2}h_{\varepsilon ,i}
u_{\varepsilon}\mathrm{d}x-\displaystyle\sum_{i=1}
^{k}\beta_{\varepsilon ,i,N+1}\int_{\mathbb{R}^N}
\eta_{i}^{2}
u_{\varepsilon}h_{\varepsilon ,i}
\big\langle
(x-\xi_{\varepsilon,i})\cdot\nabla V(x)\big\rangle\mathrm{d}x
\\&\quad-\displaystyle\sum_{i=1}
^{k}\displaystyle\sum_{j=1}
^{N}\beta_{\varepsilon ,i,N+1}\beta_{\varepsilon ,i,j}
\int_{\mathbb{R}^N}\frac{\partial V(x)}
{\partial x_{j}} \eta_{i}^{2} u_{\varepsilon}\varphi_{\varepsilon ,i}
\mathrm{d}x\\&=2(p-\varepsilon)
\displaystyle\sum_{i=1}
^{k}\displaystyle\sum_{j=1}
^{N}\beta_{\varepsilon ,i,N+1}\beta_{\varepsilon ,i,j}
\int_{B_{\delta}(\xi_{\varepsilon,i})}
u_{\varepsilon}^{p-1-\varepsilon}\big\langle
(x-\xi_{\varepsilon,i})\cdot\nabla
u_{\varepsilon}\big\rangle
\frac{\partial u_{\varepsilon}}
{\partial x_{j}}\mathrm{d}x
\\&\quad-\displaystyle\sum_{i=1}
^{k}\displaystyle\sum_{j=1}
^{N}\beta_{\varepsilon ,i,N+1}\beta_{\varepsilon ,i,j}
\frac{\partial V(\xi_{\varepsilon,i})}
{\partial x_{j}}
\int_{B_{\delta}(\xi_{\varepsilon,i})}
 u_{\varepsilon}\varphi_{i,\varepsilon}
\mathrm{d}x
+O(\varepsilon)
\\&=2(p-\varepsilon)\mu_{\varepsilon ,i,
}^{1-\frac{(N-2)\varepsilon}{2} }
\displaystyle\sum_{i=1}
^{k}\displaystyle\sum_{j,m=1}
^{N}\beta_{\varepsilon ,i,N+1}\beta_{\varepsilon ,i,j}
\int_{B_{\delta\mu_{\varepsilon,i}}(
0)}
\tilde{u}_{\varepsilon}^{p-1-\varepsilon}
x_{m}\frac{\partial \tilde{u}_{\varepsilon}}
{\partial x_{m}}
\frac{\partial \tilde{u}_{\varepsilon}}
{\partial x_{j}}\mathrm{d}x
+O(\varepsilon),
%\\&=o
%\Big(\displaystyle\sum_{i=1}
%^{k}\displaystyle\sum_{j=1}
%^{N}|\beta_{\varepsilon ,i,j}||\beta_{\varepsilon,i,N+1}|
%\mu_{\varepsilon,i}\Big)+O\Big(\displaystyle\sum_{i=1}
%^{k}\displaystyle\sum_{j=1}
%^{N}|\beta_{\varepsilon ,i,j}||\beta_{\varepsilon ,i,N+1}|
%\mu_{\varepsilon,i}^{-2 }\Big)
\end{align*}
where $\tilde{u}_{\varepsilon}=\mu_{\varepsilon,i}^{
-\frac{N-2}{2}}u_{\varepsilon}
\big(\xi_{\varepsilon,i}+x/\mu_{\varepsilon,i}\big).$

Thus we  conclude that
\begin{align}\label{6.163}
\|u\|^{2}&=(p-\varepsilon)\displaystyle\sum_{i=1}
^{k}
\displaystyle\int_{B_{\delta}(\xi_{\varepsilon,i})}
u_{\varepsilon}^{p-1-\varepsilon}
\big(\beta_{\varepsilon ,i,0}^{2}
u_{\varepsilon}^{2}+h_{\varepsilon ,i}^{2}+
\beta_{\varepsilon ,i,N+1}^{2}
\varphi_{\varepsilon ,i}^{2}\big)
\mathrm{d}x
\notag\\&\quad
+2(p-\varepsilon)\mu_{\varepsilon ,i,
}^{1-\frac{(N-2)\varepsilon}{2} }
\displaystyle\sum_{i=1}
^{k}\displaystyle\sum_{j,m=1}
^{N}\beta_{\varepsilon ,i,N+1}\beta_{\varepsilon ,i,j}
\int_{B_{\delta\mu_{\varepsilon,i}}(
0)}
\tilde{u}_{\varepsilon}^{p-1-\varepsilon}
x_{m}\frac{\partial \tilde{u}_{\varepsilon}}
{\partial x_{m}}
\frac{\partial \tilde{u}_{\varepsilon}}
{\partial x_{j}}\mathrm{d}x
\notag\\&\quad
+(1-p)\displaystyle\sum_{i=1}
^{k}\beta_{\varepsilon ,i,0}^{2}
\int_{B_{\delta}(\xi_{\varepsilon,i})}
u_{\varepsilon}^{p+1-\varepsilon}
\mathrm{d}x+O(\varepsilon).
\end{align}
On the other hand, we have
\begin{align*}
(p-\varepsilon)\int_{\mathbb{R}^N}
u_{\varepsilon}^{p-1-\varepsilon}u^{2}\mathrm{d}x
&=
(p-\varepsilon)\displaystyle\sum_{i=1}
^{k}\beta_{\varepsilon ,i,N+1}^{2}
\int_{\mathbb{R}^N}
\eta_{i}^{2}u_{\varepsilon}^{p-1-\varepsilon}
\varphi_{\varepsilon ,i}^{2}
\mathrm{d}x+Q_{0}\\&\quad
+2(p-\varepsilon)\displaystyle\sum_{i=1}
^{k}\beta_{\varepsilon ,i,0}\beta_{\varepsilon ,i,N+1}
\int_{B_{2\delta}(\xi_{\varepsilon,i})}
u_{\varepsilon}^{p-\varepsilon}\eta_{i}^{2}
\varphi_{\varepsilon ,i}
\mathrm{d}x\\&\quad
+2(p-\varepsilon)\displaystyle\sum_{i=1}
^{k}\beta_{\varepsilon ,i,N+1}
\int_{B_{2\delta}(\xi_{\varepsilon,i})}
u_{\varepsilon}^{p-1-\varepsilon}\eta_{i}^{2}
h_{\varepsilon ,i}\varphi_{\varepsilon ,i}
\mathrm{d}x
\\:&=(p-\varepsilon)\displaystyle\sum_{i=1}
^{k}\beta_{\varepsilon ,i,N+1}^{2}\bigg(
\int_{B_{\delta}(\xi_{\varepsilon,i})}
u_{\varepsilon}^{p-1-\varepsilon}
\varphi_{\varepsilon ,i}^{2}
\mathrm{d}x+
O\big(
\varepsilon^{\frac{N}{2}}\big)\bigg)+Q_{0}+Q_{2}+Q_{3},
\end{align*}
where (see \eqref{apop6.163})
\begin{align*}
Q_{0}
&=(p-\varepsilon)\displaystyle\sum_{i=1}
^{k}
\int_{B_{\delta}(\xi_{\varepsilon,i})}
u_{\varepsilon}^{p-1-\varepsilon}
\Big(
\beta_{\varepsilon ,i,0}^{2}
u_{\varepsilon}^{2}+h_{\varepsilon ,i}^{2}
\Big)
\mathrm{d}x+O(\varepsilon^{\frac{N}{2}}).\end{align*}
Thanks to Lemma \ref{lem4.1}, we  obtain
\begin{align*}
Q_{2}&=
2(p-\varepsilon)
\displaystyle\sum_{i=1}
^{k}\beta_{\varepsilon ,i,0}\beta_{\varepsilon ,i,N+1}
\int_{B_{\delta}(\xi_{\varepsilon,i})}
\bigg(\frac{\big\langle(x-\xi_{\varepsilon,i})\cdot
\nabla u_{\varepsilon}^{p+1-\varepsilon}\big\rangle}
{p+1-\varepsilon}
+\frac{2u_{\varepsilon}^{p+1-\varepsilon}}
{p-1-\varepsilon}\bigg)
\mathrm{d}x
+O(\varepsilon^{\frac{N}{2}})\\&=
\displaystyle\sum_{i=1}
^{k}\beta_{\varepsilon ,i,0}\beta_{\varepsilon ,i,N+1}
\bigg(\frac{4(p-\varepsilon)}{p-1-\varepsilon}-
\frac{2N(p-\varepsilon)}
{p+1-\varepsilon} \bigg)
\int_{B_{\delta}(\xi_{\varepsilon,i})}
u_{\varepsilon}^{p+1-\varepsilon}\mathrm{d}x
+O(\varepsilon^{\frac{N}{2}})
\\&=O(\varepsilon),
\end{align*}

and
\begin{align*}
Q_{3}&=
2(p-\varepsilon)
\displaystyle\sum_{i=1}
^{k}\displaystyle\sum_{j=1}
^{N}
\beta_{\varepsilon ,i,N+1}\beta_{\varepsilon ,i,j}
\int_{B_{\delta}(\xi_{\varepsilon,i})}
u_{\varepsilon}^{p-1-\varepsilon}
\frac{\partial u_{\varepsilon}}
{\partial x_{j}}\Big(\big\langle
(x-\xi_{\varepsilon,i})\cdot\nabla u_{\varepsilon}\big\rangle
+\frac{2u_{\varepsilon}}
{p-1-\varepsilon}\Big)
\mathrm{d}x
+O(\varepsilon^{\frac{N}{2}})\\&=
2(p-\varepsilon)\mu_{\varepsilon ,i
}^{1-\frac{(N-2)\varepsilon}{2} }
\displaystyle\sum_{i=1}
^{k}\displaystyle\sum_{j,m=1}
^{N}\beta_{\varepsilon ,i,N+1}\beta_{\varepsilon ,i,j}
\int_{B_{\delta\mu_{\varepsilon,i}}(
0)}
\tilde{u}_{\varepsilon}^{p-1-\varepsilon}
x_{m}\frac{\partial \tilde{u}_{\varepsilon}}
{\partial x_{m}}
\frac{\partial \tilde{u}_{\varepsilon}}
{\partial x_{j}}\mathrm{d}x
+O(\varepsilon^{\frac{N}{2}})
%\\&=
%o\Big(\displaystyle\sum_{i=1}
%^{k}\displaystyle\sum_{j=1}
%^{N}|\beta_{\varepsilon ,i,j}||\beta_{\varepsilon ,i,N+1}|
%\mu_{\varepsilon,i}\Big)+O\Big(\displaystyle\sum_{i=1}
%^{k}\displaystyle\sum_{j=1}
%^{N}|\beta_{\varepsilon ,i,N+1}|
%|\beta_{\varepsilon ,i,j}|
%\mu_{\varepsilon,i}^{-N }\Big)
.
\end{align*}
Thus we have
\begin{align}\label{6.162}
(p-\varepsilon)\int_{\mathbb{R}^N}
u_{\varepsilon}^{p-1-\varepsilon}u^{2}\mathrm{d}x
&=(p-\varepsilon)\displaystyle\sum_{i=1}
^{k}
\displaystyle\int_{B_{\delta}(\xi_{\varepsilon,i})}
u_{\varepsilon}^{p-1-\varepsilon}
\big(\beta_{\varepsilon ,i,0}^{2}
u_{\varepsilon}^{2}+h_{\varepsilon ,i}^{2}+
\beta_{\varepsilon ,i,N+1}^{2}
\varphi_{\varepsilon ,i}^{2}\big)
\mathrm{d}x+\displaystyle\sum_{i=1}
^{k}R_{i,1}
,
\end{align}
where
\begin{align*}
R_{i,1}&=
2(p-\varepsilon)\mu_{\varepsilon ,i,
}^{1-\frac{(N-2)\varepsilon}{2} }
\displaystyle\sum_{j,m=1}
^{N}\beta_{\varepsilon ,i,N+1}\beta_{\varepsilon ,i,j}
\int_{B_{\delta\mu_{\varepsilon,i}}(
0)}
\tilde{u}_{\varepsilon}^{p-1-\varepsilon}
x_{m}\frac{\partial \tilde{u}_{\varepsilon}}
{\partial x_{m}}
\frac{\partial \tilde{u}_{\varepsilon}}
{\partial x_{j}}\mathrm{d}x+
O(\varepsilon)\\&=
o(\mu_{\varepsilon ,i
})+O(\varepsilon)
.
\end{align*}
As a consequence, by \eqref{6.161}, \eqref{6.163} and \eqref{6.162}, we derive
\begin{align*}
J(u)&=1+
\frac{\big(1-p\big)
\displaystyle\sum_{i=1}
^{k}\beta_{\varepsilon ,i,0}^{2}
\displaystyle\int_{B_{\delta}(\xi_{\varepsilon,i})}
u_{\varepsilon}^{p+1-\varepsilon}
\mathrm{d}x+O(\varepsilon)
}
{(p-\varepsilon)\displaystyle\sum_{i=1}
^{k}
\displaystyle\int_{B_{\delta}(\xi_{\varepsilon,i})}
u_{\varepsilon}^{p-1-\varepsilon}
f_{\varepsilon ,i}
\mathrm{d}x
+\displaystyle\sum_{i=1}
^{k}R_{i,1}}
\\&\leq1+
\sup_{1\leq i\leq k}
\frac{\big(1-p+\varepsilon\big)\beta_{\varepsilon ,i,0}^{2}
\displaystyle\int_{B_{\delta}(\xi_{\varepsilon,i})}
u_{\varepsilon}^{p+1-\varepsilon}
\mathrm{d}x+O(\varepsilon)}
{(p-\varepsilon)
\displaystyle\int_{B_{\delta}(\xi_{\varepsilon,i})}
u_{\varepsilon}^{p-1-\varepsilon}
f_{\varepsilon ,i}
\mathrm{d}x
+R_{i,1}},
\end{align*}
where $f_{\varepsilon ,i}=
\beta_{\varepsilon ,i,0}^{2}
u_{\varepsilon}^{2}+h_{\varepsilon ,i}^{2}+
\beta_{\varepsilon ,i,N+1}^{2}
\varphi_{\varepsilon ,i}^{2}
$. It follows from similar arguments as in Proposition \ref{Propyu94.2} that \eqref{e2025} holds
and  the proof of Proposition \ref{Prop0094.2} is complete.
\end{proof}

To estimate the eigenvalues $\lambda_{\varepsilon,l}$ precisely, we need to study the asymptotic profile of their corresponding eigenfunctions
$\omega_{\varepsilon,l}$.
\begin{lemma}\label{lem984.3}
For $l=k+1,\cdots,k(N+2)$, there exists at least one $i_{l}\in\{1,\cdots,k\}$ and a non-zero vector  $
 (b_{l,i_{l},0},b_{l,i_{l},1},\cdots,b_{l,i_{l},N})\in\mathbb{R}^{N+1}$ such that
\begin{equation*}
\tilde{\omega}_{\varepsilon,l,i_l}\rightarrow
b_{l,i_{l},0}\Psi_{0}+
\displaystyle\sum_{j=1}^{N}b_{l,i_{l},j}
\Psi_{j}
\ \ \text{in} \ \ C^{2}_{loc}\big(\mathbb{R}^N\big), \end{equation*}
 where $\tilde{\omega}_{\varepsilon,l,i_l}=
\mu_{\varepsilon ,i_{l}}
^{-\frac{N-2}{2}}\omega_{\varepsilon,l}
(\xi_{\varepsilon ,i_{l}}+x/\mu_{\varepsilon ,i_{l}})$,
$\Psi_{0}=\frac{\partial U_{\mu,0}}{\partial \mu}|_{\mu=1}$,
$\Psi_{j}=\frac{\partial U_{1,0}}{\partial x_{j}}$.
%
%\textcolor{red}{For $l=k+1,\cdots,k(N+2)$, there exists at least one $i_{l}\in\{1,\cdots,k\}$ such that
%\begin{equation*}
%\tilde{\omega}_{\varepsilon,l,i_l}\rightarrow
%b_{l,i_{l},0}\Psi_{0}+
%\displaystyle\sum_{j=1}^{N}b_{l,i_{l},j}
%\Psi_{j}
%\ \ \text{in} \ \ C^{2}_{loc}\big(\mathbb{R}^N\big), \end{equation*}
% where the vector  $
% (b_{l,i_{l},0},b_{l,i_{l},1},\cdots,b_{l,i_{l},N})
% \neq\mathbf{0}$, $\tilde{\omega}_{\varepsilon,l,i_l}=
%\mu_{\varepsilon ,i_{l}}
%^{-\frac{N-2}{2}}\omega_{\varepsilon,l}
%(\xi_{\varepsilon ,i_{l}}+x/\mu_{\varepsilon ,i_{l}})$,
%$\Psi_{0}=\frac{\partial U_{\mu,0}}{\partial \mu}|_{\mu=1}$,
%$\Psi_{j}=\frac{\partial U_{1,0}}{\partial x_{j}}$. }

\end{lemma}
\begin{proof}[\bf Proof.]
By virtue of \eqref{aa4.18}, we have 
\begin{align*}
|\omega_{\varepsilon ,l}|
\leq C\|\omega_{\varepsilon ,l}\|_*
\displaystyle\sum_{i=1}^{k}\frac{\mu_{\varepsilon,i}^{\frac{N-2}{2}}}
{(1+\mu_{\varepsilon,i}^{2}|x-\xi_{\varepsilon,i}|^{2})^{\frac{\beta}{2}}},
\end{align*}
where $\beta\in(\frac{N-2}{2}+\sigma,N-2)$.
Then we obtain 
\begin{align*}
\bigg(
\displaystyle\sum_{i=1}^{k}
\frac{\mu_{\varepsilon,i}^{\frac{N-2}{2}}}
{\big(1+\mu_{\varepsilon,i}^{2}|x-\xi_{\varepsilon,i}|^{2}\big)
^{\frac{N-2}{4}+\frac{\sigma}{2}}}\bigg)^{-1}
|\omega_{\varepsilon ,l}(x)|
\leq  C\|\omega_{\varepsilon ,l}\|_*
\displaystyle\sum_{i=1}^{k}\frac{1}
{(1+\mu_{\varepsilon,i}^{2}|x-\xi_{\varepsilon,i}|^{2})
^{\frac{\beta}
{2}-\frac{N-2}{4}-\frac{\sigma}{2}}}
,
\end{align*}
which, together with $\|\omega_{\varepsilon ,l}\|_*=1$, yields that there exists at least one $i_{l}\in\{1,\cdots,k\}$ such that
$\bigg(
\displaystyle\sum_{i=1}^{k}
\frac{\mu_{\varepsilon,i}^{\frac{N-2}{2}}}
{\big(1+\mu_{\varepsilon,i}^{2}|x-\xi_{\varepsilon,i}|^{2}\big)
^{\frac{N-2}{4}+\frac{\sigma}{2}}}\bigg)^{-1}
|\omega_{\varepsilon ,l}(x)|$
 attains its maximum in $ B_{R\mu_{\varepsilon ,i_{l}}^{-1}}
(\xi_{\varepsilon ,i_{l}})$, where $R>0$ is a large constant.
Letting $\tilde{\omega}_{\varepsilon,l,i_{l}}=
\mu_{\varepsilon ,i_{l}}
^{-\frac{N-2}{2}}\omega_{\varepsilon,l}
(\xi_{\varepsilon ,i_{l}}+x/\mu_{\varepsilon ,i_{l}})$, then $\tilde{\omega}_{\varepsilon,l,i_{l}}$ satisfies
\begin{equation*}
-\Delta \tilde{\omega}_{\varepsilon,l,i_l}+
\mu_{i_{l},\varepsilon}
^{-2}V(\xi_{\varepsilon ,i_{l}}+x/\mu_{\varepsilon ,i_{l}})
\tilde{\omega}_{\varepsilon,l,i_l}
=(p-\varepsilon)\lambda_{\varepsilon,l}
\mu_{\varepsilon ,i_{l}}
^{-2} u_\varepsilon^{p-1-\varepsilon}
(\xi_{\varepsilon,i_{l}}+x/\mu_{\varepsilon,i_{l}})
\tilde{\omega}_{\varepsilon,l,i_l}\;\;\text{in}\;\;
\mathbb{R}^N.
\end{equation*}
It follows from the regularity theory of elliptic equation that there is $\omega_{l,i_l}^* \neq 0$ such that
\[
\tilde{\omega}_{\varepsilon,l,i_l} \to \omega_{l,i_l}^* \ \;\text{in} \;\ C^2_{loc}(\mathbb{R}^N),
\]
and
\[
-\Delta \omega_{l,i_l}^*=p U_{1,0}^{p-1}
\omega_{l,i_l}^*\ \; \ \text{in} \ \; \ \mathbb{R}^N.
\]
Then  we infer
$$\omega_{l,i_l}^*=b_{l,i_{l},0}\Psi_{0}+
\displaystyle\sum_{l=1}^{N}b_{l,i_{l},j}
\Psi_{j}, $$
where $\;\; (b_{l,i_{l},0},\cdots,b_{l,i_{l},N})\neq(0,\cdots,0)$. The proof of Lemma \ref{lem984.3} is complete.
\end{proof}

We now use the Pohozaev-type identities in Lemma \ref{lem4.3} to derive  more precise asymptotic estimates for the eigenfunctions $\omega_{\varepsilon,l}$
with $l=1,\cdots,k+1,\cdots,k(N+1)$.
\begin{lemma}\label{lem4.4}
Let $l=k+1,\cdots,k(N+1)$. Then  there exists at least one $i_{l}\in\{1,\cdots,k\}$ and a non-zero vector $\mathbf{b}_{l,i_{l}}=
 (b_{l,i_{l},1},\cdots,b_{l,i_{l},N})\in\mathbb{R}^N$ such that
\begin{align*}
    \tilde{\omega}_{\varepsilon,l,i_l}\rightarrow
\displaystyle\sum_{j=1}^{N}b_{l,i_{l},j}
\Psi_{j}
\ \ \text{in} \ \ C^{2}_{loc}\big(\mathbb{R}^N\big),\end{align*}
where $\tilde{\omega}_{\varepsilon,l,i_l}=
\mu_{\varepsilon ,i_{l}}
^{-\frac{N-2}{2}}\omega_{\varepsilon,l}
(\xi_{\varepsilon ,i_{l}}+x/\mu_{\varepsilon ,i_{l}})$ and
$\mathbf{b}_{l,i_{l}}$ is an eigenvector
of the  matrix $\big(\nabla^{2}V(\xi_{i_l}^{*})\big)$
.
%\textcolor{red}{Let $l=k+1,\cdots,k(N+1)$. Then  there exists at least one $i_{l}\in\{1,\cdots,k\}$ such that
%\begin{align*}
%    \tilde{\omega}_{\varepsilon,l,i_l}\rightarrow
%\displaystyle\sum_{j=1}^{N}b_{l,i_{l},j}
%\Psi_{j}
%\ \ \text{in} \ \ C^{2}_{loc}\big(\mathbb{R}^N\big),\end{align*}
%where the vector $\mathbf{b}_{l,i_{l}}=
% (b_{l,i_{l},1},\cdots,b_{l,i_{l},N})
% \neq\mathbf{0}$, $\tilde{\omega}_{\varepsilon,l,i_l}=
%\mu_{\varepsilon ,i_{l}}
%^{-\frac{N-2}{2}}\omega_{\varepsilon,l}
%(\xi_{\varepsilon ,i_{l}}+x/\mu_{\varepsilon ,i_{l}})$ and
%$\mathbf{b}_{l,i_{l}}$ is an eigenvector
%of the  matrix $\big(\nabla^{2}V(\xi_{i_l}^{*})\big)$
%.}

\end{lemma}
\begin{proof}[\bf Proof.]
In view of Lemma \ref{lem984.3},  for every fixed $l=k+1,\cdots,k(N+1)$, there exists at least one $i_{l}\in\{1,\cdots,k\}$ such that
\begin{align*}
\tilde{\omega}_{\varepsilon,l,i_l}\rightarrow
b_{l,i_{l},0}\Psi_{0}+
\displaystyle\sum_{j=1}^{N}b_{l,i_{l},j}
\Psi_{j}
\ \ \text{in} \ \ C^{2}_{loc}\big(\mathbb{R}^N\big).\end{align*}
Thus we only need to prove that  $b_{l,i_{l},0}=0$. We proceed by contradiction and assume that
$b_{l,i_{l},0}\neq0$. %For simplicity, we omit the subscript of $i_{0}$.
Thanks to  Lemma \ref{lem4.1}, we have
\begin{align*}
 \int_{B_\delta(\xi_{\varepsilon,i_l})}   u_\varepsilon^{p-\varepsilon}
\omega_{\varepsilon,l}
\mathrm{d}x
=&\int_{B_\delta\mu_{\varepsilon,i_l}(0)}
U_{1,0}^{p-\varepsilon}
\tilde{\omega}_{\varepsilon,l,i_l}\mathrm{d}x
+O\big(\|\phi_\varepsilon\|_*\big)
+O\big(\|\phi_\varepsilon\|_*^{p}\big)\\
=&
b_{l,i_l,0}
\int_{\mathbb{R}^N}U_{1,0}^{p}\Psi_0\mathrm{d}x+o(1)=o(1).
\end{align*}
Then we obtain
\begin{align*}\label{a4.4}&
\frac{N-2}{2}(p-\varepsilon)
\big(\lambda_{\varepsilon,l}-1\big)\int_{B_\delta(\xi_{\varepsilon,i_l})} u_\varepsilon^{p-\varepsilon}
\omega_{\varepsilon,l}\mathrm{d}x
-\frac{N-2}{2}\varepsilon\int_{B_\delta(\xi_{\varepsilon,i_l})} u_\varepsilon^{p-\varepsilon}
\omega_{\varepsilon,l}\mathrm{d}x
=o\big(\lambda_{\varepsilon,l}-1\big)+o(\varepsilon).
\end{align*}
A direct calculation leads to
\begin{align*}
\int_{B_\delta(\xi_{\varepsilon,i_l})}\big\langle x-\xi_{\varepsilon,i_l},\nabla V(x)\big\rangle u_\varepsilon\omega_{\varepsilon,l}\mathrm{d}x
\leq C\varepsilon^{\frac{3}{2}}.
\end{align*}
Using Proposition \ref{prop2.4}, we infer

\begin{align*}
&2\int_{B_\delta(\xi_{\varepsilon,i_l})}V(x)u_\varepsilon
\omega_{\varepsilon ,l}\mathrm{d}x\\
=&
2\int_{B_\delta(\xi_{\varepsilon,i_l})} V(x)
\Big(U_{\mu_{\varepsilon,i_l},\xi_{\varepsilon,i_l}}
+\phi_{\varepsilon}\Big)\omega_{\varepsilon,l}\mathrm{d}x\\
=&2V(\xi_{\varepsilon,i_l})\int_{B_\delta(\xi_{\varepsilon,i_l})}
U_{\mu_{\varepsilon,i_l},\xi_{\varepsilon,i_l}}
\omega_{\varepsilon,l}\mathrm{d}x
+O\big(\mu_{\varepsilon,i_l}^{-2}
\|\phi_{\varepsilon}\|_*\big)
+O\left(\varepsilon^2|\log\varepsilon|
\right)
\\
=&2V(\xi_{\varepsilon,i_l})b_{l,i_l,0}\big(1+
o(1)\big)\mu_{\varepsilon,i_l}^{-2}
\int_{\mathbb{R}^N}U_{1,0}\Psi_0
\mathrm{d}x+O\big(\varepsilon^{2-\frac{\sigma}{2}}\big)
\\
=&-2AV(\xi_{i_l}^{*})\mu_{\varepsilon,i_l}^{-2}b_{l,i_l,0}\big(1+
o(1)\big)
+O\big(\varepsilon^{2-\frac{\sigma}{2}}\big).
\end{align*}

On the other hand, it is easy to check that
%%\begin{align*}
%%&\int_{B_\delta(\xi_{\varepsilon,i_l})}
%%\big\langle x-\xi_{\varepsilon,i_l},\nabla u_\varepsilon\big\rangle u_\varepsilon^{p-1-\varepsilon}
%%\omega_{\varepsilon,l}
%%\mathrm{d}x\\
%%%=&\int_{B_\delta(\xi_{\varepsilon,i_l})}
%%%\Big\langle x-\xi_{\varepsilon,i_l},\nabla \big(
%%%U_{\mu_{\varepsilon,i_l},\xi_{\varepsilon,i_l}}
%%%+\phi_{\varepsilon}\big)\Big\rangle u_\varepsilon^{p-1-\varepsilon}
%%%\omega_{\varepsilon,l}
%%%\mathrm{d}x
%%%\\
%%=&\int_{B_\delta(\xi_{\varepsilon,i_l})}
%%\big\langle x-\xi_{\varepsilon,i_l},\nabla
%%U_{\mu_{\varepsilon,i_l},\xi_{\varepsilon,i_l}}
%%\big\rangle
%%U^{p-1-\varepsilon}_{\mu_{\varepsilon,i_l},\xi_{\varepsilon,i_l}}
%%\omega_{\varepsilon,l}
%%\mathrm{d}x
%%+O\big(\|\phi_\varepsilon\|_*\big)\\
%%&\quad
%%+O\Big(\int_{B_\delta(\xi_{\varepsilon,i_l})}
%%| x-\xi_{\varepsilon,i_l}|\big|\nabla
%%U_{\mu_{\varepsilon,i_l},\xi_{\varepsilon,i_l}}
%%\big| |\phi_{\varepsilon}|^{p-1-\varepsilon}
%%|\omega_{\varepsilon,l}|
%%\mathrm{d}x\Big)
%%\\
%%=&\int_{B_\delta(\xi_{\varepsilon,i_l})}
%%\big\langle x-\xi_{\varepsilon,i_l},\nabla
%%U_{\mu_{\varepsilon,i_l},\xi_{\varepsilon,i_l}}
%%\big\rangle U_{\mu_{\varepsilon,i_l},\xi_{\varepsilon,i_l}}^{p-1-\varepsilon}
%%\omega_{\varepsilon ,l}
%%\mathrm{d}x+O\big(\|\phi_\varepsilon\|^{p-1-\varepsilon}_*\big)
%%\\
%%=&D_0b_{l,i_l,0}\big(1+o(1)\big)
%%+O\big(
%%\|\phi_\varepsilon\|^{p-1-\varepsilon}_*\big),
%%\end{align*}
\begin{align*}
&\int_{B_\delta(\xi_{\varepsilon,i_l})}
\big\langle x-\xi_{\varepsilon,i_l},\nabla u_\varepsilon\big\rangle u_\varepsilon^{p-1-\varepsilon}
\omega_{\varepsilon,l}
\mathrm{d}x\\=&
\int_{B_\delta(\xi_{\varepsilon,i_l})}
\big\langle x-\xi_{\varepsilon,i_l},\nabla u_\varepsilon\big\rangle U_{\mu_{\varepsilon,i_l},
\xi_{\varepsilon,i_l}}^{p-1-\varepsilon}
\omega_{\varepsilon,l}
\mathrm{d}x
+O\big(
\|\phi_\varepsilon\|^{p-1-\varepsilon}_*\big)\\
%=&\int_{B_\delta(\xi_{\varepsilon,i_l})}
%\big\langle x-\xi_{\varepsilon,i_l},\nabla
%U_{\mu_{\varepsilon,i_l},
%\xi_{\varepsilon,i_l}}+\nabla
%\phi_{\varepsilon}
%\big\rangle
%U^{p-1-\varepsilon}_{\mu_{\varepsilon,i_l},\xi_{\varepsilon,i_l}}
%\omega_{\varepsilon,l}
%\mathrm{d}x
%+O\big(
%\|\phi_\varepsilon\|^{p-1-\varepsilon}_*\big)\\
%=&\int_{B_\delta(\xi_{\varepsilon,i_l})}
%\big\langle x-\xi_{\varepsilon,i_l},\nabla
%U_{\mu_{\varepsilon,i_l},
%\xi_{\varepsilon,i_l}}
%\big\rangle
%U^{p-1-\varepsilon}_{\mu_{\varepsilon,i_l},\xi_{\varepsilon,i_l}}
%\omega_{\varepsilon,l}
%\mathrm{d}x
%-
%\int_{B_\delta(\xi_{\varepsilon,i_l})}
% (x-\xi_{\varepsilon,i_l})
%\phi_{\varepsilon}
%\nabla\big(
%U^{p-1-\varepsilon}_{\mu_{\varepsilon,i_l},\xi_{\varepsilon,i_l}}
%\omega_{\varepsilon,l}\big)
%\mathrm{d}x
%\\
%&\quad
%-\int_{B_\delta(\xi_{\varepsilon,i_l})}
%\phi_{\varepsilon}
%U^{p-1-\varepsilon}_{\mu_{\varepsilon,i_l},\xi_{\varepsilon,i_l}}
%\omega_{\varepsilon,l}
%\mathrm{d}x
%+
%\int_{\partial B_\delta(\xi_{\varepsilon,i_l})}
%\big\langle x-\xi_{\varepsilon,i_l},
%\gamma_{\varepsilon,i_l}
%\big\rangle\phi_{\varepsilon}
%U^{p-1-\varepsilon}_{\mu_{\varepsilon,i_l},\xi_{\varepsilon,i_l}}
%\omega_{\varepsilon,l}
%\mathrm{d}s
%+O\big(
%\|\phi_\varepsilon\|^{p-1-\varepsilon}_*\big)
%\\
=&\int_{B_\delta(\xi_{\varepsilon,i_l})}
\big\langle x-\xi_{\varepsilon,i_l},\nabla
U_{\mu_{\varepsilon,i_l},
\xi_{\varepsilon,i_l}}
\big\rangle
U^{p-1-\varepsilon}_{\mu_{\varepsilon,i_l},\xi_{\varepsilon,i_l}}
\omega_{\varepsilon,l}
\mathrm{d}x
%+O\big(
%\mu_{\varepsilon,i_l}^{-\frac{N+2+\sigma}{2}}
%\|\phi_\varepsilon\|_*\big)
+O\big(\|
\phi_\varepsilon\|_*\big)+O\big(
\|\phi_\varepsilon\|^{p-1-\varepsilon}_*\big)
\\
=&D_0b_{l,i_l,0}\big(1+o(1)\big)
+O\big(
\|\phi_\varepsilon\|^{p-1-\varepsilon}_*\big),
\end{align*}
where $D_0=\int_{\mathbb{R}^N}\big\langle x,\nabla
U_{1,0}\big\rangle U_{1,0}^{p-1}\Psi_0
\mathrm{d}x$.

Then we obtain
\begin{align*}
\text{LHS\;of\;}\eqref{4.3}=&(p-\varepsilon)
(\lambda_{\varepsilon,l}-1)\bigg(
 D_0b_{l,i_l,0}\big(1+o(1)\big)
+O\big(
\|\phi_\varepsilon\|^{p-1-\varepsilon}_*\big)
\bigg)
 \\
&-2AV(\xi_{i}^{*})b_{l,i_l,0}\big(1+
o(1)\big)\mu_{\varepsilon,i_l}^{-2}
+o(\varepsilon).
\end{align*}
It follows from  Lemma \ref{lem4.1} that
 \begin{equation*}
\begin{split}
\text{RHS\;of\;}\eqref{4.3} =O\big(\varepsilon^{\frac{N-2}{2}}\big).
\end{split}
\end{equation*}
Hence, we deduce
$$
\lambda_{\varepsilon,l}-1=\tilde{M}_{l}
\mu_{\varepsilon,i_l}^{-2}\big(1+
o(1)\big),\ \ l=k+1,\cdots,k(N+1),
$$
where $\tilde{M}_{l}=\frac{2AV(\xi_{i_l}^{*})}
{pD_0}>0$.
As a consequence, we get
\begin{align}\label{an}
\lambda_{\varepsilon,l}=1+\tilde{M}_{l}
\varepsilon
+o\big(\varepsilon\big)>1, \;\; l=k+1,\cdots,k(N+1),\end{align}
which is a  contradiction with  \eqref{po0e2025}. Thus $b_{l,i_{l},0}=0$.

 Without loss of generality, we may assume $b_{l,i_{l},j}\neq0$ for some $j\in\{1,\cdots,N\}$. By virtue of  Lemma \ref{lem4.1}, we have
\begin{align*}
&\int_{ B_{\delta}(\xi_{\varepsilon,i_l})}  u_\varepsilon^{p-1-\varepsilon}
\omega_{\varepsilon,l}
\frac{\partial u_\varepsilon}{\partial x_{j}}\mathrm{d}x
\\
=&-\frac{1}{p-\varepsilon}
\int_{ B_{\delta}(\xi_{\varepsilon,i_l})}
u_\varepsilon^{p-\varepsilon}
\frac{\partial\omega_{\varepsilon,l} }{\partial x_{j}}\mathrm{d}x
+\frac{1}{p-\varepsilon}
\int_{\partial B_{\delta}(\xi_{\varepsilon,i_l})}
u_\varepsilon^{p-\varepsilon}
\omega_{\varepsilon ,l}\nu_{j}
\mathrm{d}s
\\
%=&-\frac{1}{p-\varepsilon}
%\int_{ B_{\delta}(\xi_{\varepsilon,i_l})}
%U_{\mu_{\varepsilon,i_l},\xi_{\varepsilon,i_l}}
%^{p-\varepsilon}
%\frac{\partial\omega_{\varepsilon,l} }{\partial x_{j}}\mathrm{d}x
%+
%O\bigg(
%\int_{ B_{\delta}(\xi_{\varepsilon,i_l})}
%U^{p-1-\varepsilon}
%_{\mu_{\varepsilon,i_l},\xi_{\varepsilon,i_l}}
%\big|\phi_{\varepsilon}\big|\Big|
%\frac{\partial\omega_{\varepsilon,l} }{\partial x_{j}}\Big|\mathrm{d}x
%\bigg)
%\\
%&+
%O\bigg(
%\int_{ B_{\delta}(\xi_{\varepsilon,i_l})}
%\big|\phi_{\varepsilon}\big|^{p-\varepsilon}
%\Big|
%\frac{\partial\omega_{\varepsilon,l} }{\partial x_{j}}\Big|\mathrm{d}x
%\bigg)
%+
%O\bigg(\int_{\partial B_{\delta}(\xi_{\varepsilon,i_l})}
%u_\varepsilon^{p-\varepsilon}
%\big|\omega_{\varepsilon,l}\big|\nu_{j}
%\mathrm{d}x\bigg)
%\\
=&-\frac{1}{p-\varepsilon}
\int_{ B_{\delta}(\xi_{\varepsilon,i_l})}
U_{\mu_{\varepsilon,i_l},\xi_{\varepsilon,i_l}}
^{p-\varepsilon}
\frac{\partial\omega_{\varepsilon,l} }{\partial x_{j}}\mathrm{d}x
+O\Big(\mu_{\varepsilon,i_l}\big(\|\phi_\varepsilon
\|_*+
\|\phi_\varepsilon\|^{p-\varepsilon}_*\big)\Big)
+
O\big(\mu_{\varepsilon,i_l}^{-N}\big)
\\
=&D_j\mu_{\varepsilon,i_l}b_{l,i_l,j}\big(1+o(1)\big)
+
O\big(\varepsilon^{\frac{1-\sigma}{2}}\big),
 \end{align*}
 where $D_j=\int_{\mathbb{R}^N}
U_{1,0}^{p-1}\Psi_j^2
 \mathrm{d}x$.
It follows from Proposition \ref{prop2.4} and Lemma \ref{lem4.1}
that
\begin{align*}
&\int_{B_{\delta}
(\xi_{\varepsilon,i_l})}\frac{\partial V(x)}{\partial x_{j}}u_\varepsilon\omega_{\varepsilon,l}
\mathrm{d}x\\
%=&\int_{B_{\varepsilon^{\tau_{0}}}(\xi_{\varepsilon,i_l})}
%\bigg(
%\frac{\partial V(x)}{\partial x_{j}}-
%\frac{\partial V(\xi_{\varepsilon,i_l})}{\partial x_{j}}+\frac{\partial V(\xi_{\varepsilon,i_l})}{\partial x_{j}}
%\bigg)u_\varepsilon\omega_{\varepsilon ,l}\mathrm{d}x
%+\int_{B_{\delta}(\xi_{\varepsilon,i_l})\backslash
%B_{\varepsilon^{\tau_{0}}}
%(\xi_{\varepsilon,i_l})}
%\frac{\partial V(x)}{\partial x_{j}}u_\varepsilon\omega_{\varepsilon,l}
%\mathrm{d}x
%\\[0.02mm]
=&
\int_{B_{\varepsilon^{\tau_{0}}}(\xi_{\varepsilon,i_l})}
\displaystyle\sum_{m=1}
^{N}\frac{\partial^{2}V(\xi_{\varepsilon,i_l})}{\partial
x_{j} \partial x_{m}}(x_{m}-\xi_{\varepsilon,i_l,m})
u_\varepsilon\omega_{\varepsilon,l}\mathrm{d}x
+
o\big(\mu_{\varepsilon,i_l}^{-3}\big)
\\[0.02mm]
%=&
%\int_{B_{\varepsilon^{\tau_{0}}}
%(\xi_{\varepsilon,i_l})}
%\bigg(\displaystyle\sum_{m=1}
%^{N}\frac{\partial^{2}V(\xi_{\varepsilon,i_l})}{\partial
%x_{j} \partial x_{m}}(x_{m}-\xi_{\varepsilon,i_l,m})
%\bigg)
%\Big(
%U_{\mu_{\varepsilon,i},\xi_{\varepsilon,i_l}}
%+\phi_{\varepsilon}\Big)\omega_{\varepsilon,l_{t}}\mathrm{d}x
%+
%o\big(\mu_{\varepsilon,i_l}^{-3}\big)
%\\[0.02mm]
=&
\int_{B_{\varepsilon^{\tau_{0}}}
(\xi_{\varepsilon,i_l})}
\displaystyle\sum_{m=1}
^{N}\frac{\partial^{2}V(\xi_{\varepsilon,i_l})}{\partial
x_{j} \partial x_{m}} (x_{m}-\xi_{\varepsilon,i_l,m})
U_{\mu_{\varepsilon,i_l},\xi_{\varepsilon,i_l}}
\omega_{\varepsilon,l}\mathrm{d}x+
o\big(\mu_{\varepsilon,i_l}^{-3}\big)
+
O\big(\mu_{\varepsilon,i_l}^{-2}\|\phi_\varepsilon
\|_*\big)
\\[0.02mm]
=&\displaystyle\sum_{m=1}
^{N}\frac{\partial^{2}V(\xi_{i_l}^{\ast})}{\partial
x_{j} \partial x_{m}}b_{l,i_l,m}\big(1+
o(1)\big)\mu_{\varepsilon,i_l}^{-3}
\int_{\mathbb{R}^N}
U_{1,0}\frac{\partial U_{1,0}}{\partial x_{m}}x_{m}\mathrm{d}x+
o\big(\mu_{\varepsilon,i_l}^{-3}\big)\\[0.02mm]
=&-\frac{A}{2}\displaystyle\sum_{m=1}
^{N}\frac{\partial^{2}V(\xi_{i_l}^{\ast})}{\partial
x_{j} \partial x_{m}}b_{l,i_l,m}\big(1+
o(1)\big)\mu_{\varepsilon,i_l}^{-3}
+
o\big(\varepsilon^{\frac{3}{2}}\big),
\end{align*}
where $\tau_{0}>0$ is a small  constant.
Then we conclude
 \begin{align*}
\begin{split}
\text{LHS\;of\;}\eqref{4.003}=
&(p-\varepsilon)(\lambda_{\varepsilon,l}-1)\Big(D_j\mu_{\varepsilon,i_l}b_{l,i_l,j}\big(1+o(1)\big)
+
O\big(\varepsilon^{\frac{1-\sigma}{2}}\big)\Big)\\[0.02mm]
&-\frac{A}{2}\displaystyle\sum_{m=1}
^{N}\frac{\partial^{2}V(\xi_{i_l}^{\ast})}{\partial
x_{j} \partial x_{m}}b_{l,i_l,m}\big(1+
o(1)\big)\mu_{\varepsilon,i_l}^{-3}+
o\big(\varepsilon^{\frac{3}{2}}\big).
\end{split}
 \end{align*}
Similarly, by Lemma \ref{lem4.1}, we see that
 \begin{align*}
\text{RHS\;of\;}\eqref{4.003}=O\big(\varepsilon^{\frac{N-2}{2}}\big).
 \end{align*}
% \begin{align*}
%\text{RHS\;of\;}\eqref{4.003}=&O\Bigg(\int_{\partial B_{\delta}
%(\xi_{\varepsilon,i})}
% \bigg(\displaystyle\sum_{i=1}^{k}
% \frac{\mu_{\varepsilon,i}^{\frac{N}{2}}}
%{\left(1+\mu_{\varepsilon,i}^{2}|x-\xi_{\varepsilon,i}|^{2}\right)^{\frac{N-1}{2}}}
%\bigg)^{2}\mathrm{d}x\Bigg)
% \\
%&
% +O\Bigg(\int_{\partial B_{\delta}
%(\xi_{\varepsilon,i})}\bigg(
% \displaystyle\sum_{i=1}^{k}\frac{\mu_{\varepsilon,i}^{\frac{N-2}{2}}}
%{\left(1+\mu_{\varepsilon,i}^{2}|x-\xi_{\varepsilon,i}|^{2}\right)^{\frac{N-2}{2}}}
% \bigg)^{p+1-\varepsilon}\mathrm{d}x\Bigg)\\
%&+O\Bigg(\int_{\partial B_{\delta}
%(\xi_{\varepsilon,i})}
%\bigg(
% \displaystyle\sum_{i=1}^{k}\frac{\mu_{\varepsilon,i}^{\frac{N-2}{2}}}
%{\left(1+\mu_{\varepsilon,i}^{2}|x-\xi_{\varepsilon,i}|^{2}\right)^{\frac{N-2}{2}}}
% \bigg)^{2}\mathrm{d}x\Bigg)
% \\
%=&O\big(\mu_{\varepsilon,i}^{2-N}\big).
% \end{align*}
Hence, we obtain
\begin{align*}
&(p-\varepsilon)(\lambda_{\varepsilon,l}-1)\Big(D_j\mu_{\varepsilon,i_l}b_{l,i_l,j}\big(1+o(1)\big)
+
O\big(\varepsilon^{\frac{1-\sigma}{2}}\big)\Big)\\
=&\frac{A}{2}\displaystyle\sum_{m=1}
^{N}\frac{\partial^{2}V(\xi_{i_l}^{\ast})}{\partial
x_{j} \partial x_{m}}b_{l,i_l,m}\big(1+
o(1)\big)\mu_{\varepsilon,i_l}^{-3}+
o\big(\varepsilon^{\frac{3}{2}}\big),
\end{align*}
which gives that
\begin{align*}
\tilde{M} b_{l,i_l,j}^{-1}\displaystyle\sum_{m=1}
^{N}\frac{\partial^{2}V(\xi_{i_l}^{\ast})}{\partial
x_{j} \partial x_{m}}b_{l,i_l,m}+
o(1)=(\lambda_{\varepsilon ,l}-1)
\mu_{\varepsilon,i_l}^{4},
 \end{align*}
where $\tilde{M}=\frac{A}{2pD_j}>0$. Therefore, we obtain
\begin{align}\label{aa4.34}
\begin{cases}
(\lambda_{\varepsilon ,l}-1)
\mu_{\varepsilon,i_l}^{4}\rightarrow \tilde{M}\beta_{i_l},
\\[3mm]
\displaystyle\sum_{m=1}
^{N}\frac{\partial^{2}V(\xi_{i_l}^{*})}{\partial
x_{j} \partial x_{m}}b_{l,i_l,m}=\beta_{i_l} b_{l,i_l,j}.
\end{cases}
\end{align}
%It is easy to see that if $b_{l_t,i_0,j}=0$, the above equality holds.
Then  $\beta_{i_l}$ is an eigenvalue of
$\big(\nabla^{2}V(\xi_{i_l}^{*})\big)$ and  $\mathbf{b}_{l,i_l}=\{b_{l,i_l,1},\cdots,
b_{l,i_l,N}\}$ is an eigenvector of $\big(\nabla^{2}V(\xi_{i_l}^{*})\big)$.
This completes the proof of Lemma \ref{lem4.4}.
\end{proof}

\begin{remark}\label{rema4.6}
For $l=k+1,\cdots,k(N+2)$, it follows from the proof of Lemma \ref{lem4.4} that if
$b_{l,i_l,0}\neq0$, then  $b_{l,i_l,j}=0$ for all $j=1,\cdots,N$ and  the converse also  holds.
\end{remark}

\begin{lemma}\label{lem4.p40}
For each $i\in\{1,\cdots,k\}$, there exist
  exactly $N$ orthogonal eigenfunctions $\omega_{\varepsilon,l^s_{i}}$  such that
     the corresponding  rescaled eigenfunctions $\tilde{\omega}_{\varepsilon,l^s_{i},i}$ satisfy
\begin{align}\label{1an2}
 \tilde{\omega}_{\varepsilon, l^s_{i},i}\rightarrow
\displaystyle\sum_{j=1}^{N}b_{l^s_{i},i,j}
\Psi_{j}
\ \ \text{in} \ \ C^{2}_{loc}\big(\mathbb{R}^N\big),
\end{align}
where $s=1,\cdots,N$, $l^s_{i}\in\{k+1,\cdots,k(N+1)\}$, $\mathbf{b}_{l^s_{i},i}$ is as in Lemma \ref{lem4.4}
and $\big\langle\mathbf{b}_{l_{i}^{s},i},
\mathbf{b}_{l_{i}^{t},i}\big\rangle=0$ for $s\neq t$.
%$k+1\leq l_{i}^{a}\neq l_{i}^{b} \leq k(N+1)$, $1\leq a\neq b\leq N$.
\end{lemma}
\begin{proof}
We argue by contradiction and assume that for some $i_{0}\in\{1,\cdots,k\}$, there are at least $s$ orthogonal eigenfunctions $\big\{\omega_{\varepsilon,l_{i_{0}}^{1}},\cdots,
 \omega_{\varepsilon,l_{i_{0}}^{s}}\big\}$ satisfying \eqref{1an2}, $l_{i_{0}}^{t}\in\big\{k+1,\cdots,k(N+1)
 \big\}$, $t=1,\cdots,s$
 and $N+1\leq s\leq kN$.
Using the orthogonality condition and \eqref{1an2}, we  have
\begin{align*}
\displaystyle\sum_{j=1}^{N}b_{l_{i_{0}
}^{t},i_{0},j}
b_{l_{i_{0}}^{m},i_{0},j}
\int_{\mathbb{R}^N} U_{1,0}^{p-1}
\Psi_{j}^{2}\mathrm{d}x
=0\;\ \
\text{for} \ \   t\neq m,
\end{align*}
which gives that
\[
\big\langle\mathbf{b}_{l^t_{i_0},i_{0}},
\mathbf{b}_{l^m_{i_0},i_{0}}\big\rangle=0.
\]
%\begin{align*}
%(b_{l_{i_{0}
%}^{t},i_{0},1},\cdots,
%b_{l_{i_{0}
%}^{t},i_{0},N})\cdot
%(b_{l_{i_{0}
%}^{m},i_{0},1},\cdots,
%b_{l_{i_{0}
%}^{m},i_{0},N})=0.
%\end{align*}
By Lemma \ref{lem4.4}, we know  that  $\mathbf{b}_{l_{i_{0}
}^{t},i_0}=\{b_{l_{i_{0}
}^{t},i_0,1},\cdots,b_{l_{i_{0}
}^{t},i_0,N}\}$ is an eigenvector  of  $\big(\nabla^{2}V(\xi_{i_0}^{*})\big)$.
Then we  conclude that the  matrix
$\big(\nabla^{2}V(\xi_{i_{0}}^{*})\big)$ admits at least $N+1$ orthogonal eigenvectors, which is a contradiction.
\end{proof}
\begin{lemma}\label{lem4.p4}
Let $l=k+1,\cdots,k(N+1)$. Then there exists a unique $i_{l}\in\{1,\cdots,k\}$ such that
 \begin{align}\label{1an1}
 \tilde{\omega}_{\varepsilon, l,i_l}\rightarrow
\displaystyle\sum_{j=1}^{N}b_{l,i_{l},j}
\Psi_{j}
\ \ \text{in} \ \ C^{2}_{loc}\big(\mathbb{R}^N\big),
\end{align}
where $\tilde{\omega}_{\varepsilon,l,i_l}$ and
$\mathbf{b}_{l,i_{l}}=
 (b_{l,i_{l},1},\cdots,b_{l,i_{l},N})$ are as in Lemma \ref{lem4.4}.
\end{lemma}
\begin{proof}[\bf Proof.]
On the contrary, we assume that there exists another  $i_l^\prime\in\{1,\cdots,k\}$ such that
 the corresponding rescaled  eigenfunction $\tilde{\omega}_{\varepsilon,l,i_l^\prime}$ satisfying \eqref{1an1}.
 Then  for some $i_{l}^*\in\{1,\cdots,k\}$, there exist at least $s$ orthogonal eigenfunctions $\big\{\omega_{\varepsilon,l_{1}},\cdots,
 \omega_{\varepsilon,l_{s}}\big\}$ such that
 \begin{align*}%\label{1an}
 \tilde{\omega}_{\varepsilon, l_{t},i_{l}^*}\rightarrow
\displaystyle\sum_{j=1}^{N}b_{l_{t},i_{l}^*,j}
\Psi_{j}
\ \ \text{in} \ \ C^{2}_{loc}\big(\mathbb{R}^N\big),
\end{align*}
where $(b_{l_{t},i_{l}^*,1},\cdots,b_{l_{t},i_{l}^*,N})
\neq(0,\cdots,0),$ $l_{t}\in\big\{k+1,\cdots,k(N+1)\big\}$, $t=1,\cdots,s$ and $N+1\leq s\leq kN$.
Thus we arrive at   a contradiction with Lemma \ref{lem4.p40}.
\end{proof}
Based on the above results, we can investigate the precise asymptotic estimates of $\lambda_{\varepsilon,l}$ for $l=k+1,\cdots,k(N+1)$.
Let $\beta_{i,m}$ be the eigenvalue of the matrix $\big(\nabla^{2}V(\xi_{i}^{*})\big)$, where $i=1,\cdots,k$, $m=1,\cdots,N$. Without loss of generality,
we can assume that $\beta_{i,1}\leq \beta_{i,2}\leq\cdots\leq\beta_{i,N}.$

\begin{lemma}\label{lem4.04}
For each $l=k+1,\cdots,k(N+1)$, there exists a unique $i_l\in\{1,2,\cdots,k\}$ such that
\begin{align*}
\lambda_{\varepsilon,l}=1+
\tilde{M}\beta_{i_l,m_l}\varepsilon^2+
o\big(\varepsilon^2\big),
\end{align*}
where $\tilde{M}>0$ is a constant independent of $\varepsilon$ and $\beta_{i_l,m_l}$ is an eigenvalue
of  matrix $\big(\nabla^{2}V(\xi_{i_l}^{*})\big)$ corresponding to the eigenvector
 $\mathbf{b}_{l,i_l}$.
\end{lemma}
\begin{proof}[\bf Proof.]
It follows from  \eqref{aa4.34} and \eqref{1an1} that for every $l=k+1,\cdots,k(N+1)$,
there exists only one $i_{l}=1,\cdots,k$ and $\mathbf{b}_{l,i_{l}}\neq0$ such that
\begin{equation*}
\begin{cases}
(\lambda_{\varepsilon ,l}-1)
\mu_{\varepsilon ,i_{l}}^{4}\rightarrow \tilde{M}\beta_{i_l,m_l},
\\[3mm]
\displaystyle\sum_{h=1}
^{N}\frac{\partial^{2}V(\xi_{i_{l}}^{*})}{\partial
x_{j} \partial x_{h}}b_{l,i_{l},h}=\beta_{i_l,m_l} b_{l,i_{l},j}.
\end{cases}
\end{equation*}
Hence  $\beta_{i_l,m_l}$ is an eigenvalue of the matrix
$\big(\nabla^{2}V(\xi_{i_{l}}^{*})\big)$ corresponding to the eigenvector
$\mathbf{b}_{l,i_{l}}$ and
\begin{equation*}
\lambda_{\varepsilon,l}=1+
\frac{\tilde{M}\beta_{i_l,m_l}}{\mu_{\varepsilon ,i_{l}}^{4}}+
o\big(\mu_{\varepsilon ,i_{l}}^{-4}\big).
\end{equation*}
%where $\mathbf{b}_{l,i_{l}}=\big(b_{l,i_{l},1},
%\cdots,b_{l,i_{l},N}\big)
%\neq(0,\cdots,0)$.
\end{proof}
Finally, we establish the    asymptotic profiles of the  eigenfunctions $\omega_{\varepsilon,l}$ of \eqref{aa4.1} for  $l=k(N+1)+1,\cdots,k(N+2)$.
\begin{lemma}\label{lem4.6e789} For each  $l=k(N+1)+1,\cdots,k(N+2)$,   there is a unique  $i_{l}\in\{1,\cdots,k\}$ and $b_{l,i_{l},0}\neq0$ such that
 \begin{align}\label{898.6e7}
\tilde{\omega}_{\varepsilon,l,i_{l}}\rightarrow
b_{l,i_{l},0}\Psi_{0}
\ \ \text{in} \ \ C^{2}_{loc}\big(\mathbb{R}^N\big),
\end{align}
where $\tilde{\omega}_{\varepsilon,l,i_{l}}=\mu_{\varepsilon,i_l}^{-\frac{N-2}{2}}\omega_{\varepsilon,l}
(\xi_{\varepsilon,i_l}+x/\mu_{\varepsilon,i_l})$.
\end{lemma}
\begin{proof}[\bf Proof.]
From Lemma \ref{lem984.3}, we  see that for every fixed $l=k(N+1)+1,\cdots,k(N+2)$, there exists  $i_{l}\in\{1,\cdots,k\}$ such that
$$\tilde{\omega}_{\varepsilon,l,i_{l}}\rightarrow
b_{l,i_{l},0}\Psi_{0}+
\displaystyle\sum_{j=1}^{N}b_{l,i_{l},j}
\Psi_{j}
\ \ \text{in} \ \ C^{2}_{loc}\big(\mathbb{R}^N\big),$$
where $(b_{l,i_{l},0},\cdots,b_{l,i_{l},N})
\neq(0,\cdots,0).$

We first show that for every fixed $l=k(N+1)+1,\cdots,k(N+2)$, we have
$$\mathbf{b}_{l,i_{l}}=\big(b_{l,i_{l},1},
\cdots,b_{l,i_{l},N}\big)
=(0,\cdots,0).$$
 Assume on the contrary that there exists some $j=1,\cdots,N$ such that $b_{l,i_{l},j}\neq0$. Then,  by
Remark \ref{rema4.6}, we find that  $b_{l,i_{l},0}=0$ and
$\mathbf{b}_{l,i_{l}}$ is an eigenvector
of the matrix $\big(\nabla^{2}V(\xi_{i_{l}}^{*})\big)$.  Thus
\begin{align}\label{lemkui}
\int_{\mathbb{R}^N} u_{\varepsilon}^{p-1-\varepsilon}(x)
\omega_{\varepsilon,l}(x)
\omega_{\varepsilon,l_{1}}(x)\mathrm{d}x=0\;\;\;
\text{for}\;\;
\  l_1\in\{k+1,\cdots,k(N+1)\}.
\end{align}
Moreover, it follows from Lemma \ref{lem4.p40}  that there exist  only $N$ eigenfunctions
$\omega_{\varepsilon,l_{s}}(x)$, $l_{s}\in\big\{k+1,\cdots,k(N+1)\big\}$, $s=1,\cdots,N$ such that
\begin{align*}
\begin{cases}
\tilde{\omega}_{\varepsilon,l_{s},i_{l}}\rightarrow
\displaystyle\sum_{j=1}^{N}b_{l_{s},i_{l},j}
\Psi_{j}\neq0
\ \ \text{in} \ \ C^{2}_{loc}\big(\mathbb{R}^N\big),
\\[3mm]
\big\langle\mathbf{b}_{l_s,i_{l}},
\mathbf{b}_{l_t,i_{l}}\big\rangle=0,  \ \ k+1\leq l_{s}, l_{t} \leq k(N+1), \ \  s\neq t.
\end{cases}
\end{align*}
which, combining with \eqref{lemkui}, implies that $\big\langle\mathbf{b}_{l,i_{l}},
\mathbf{b}_{l_s,i_{l}}\big\rangle=0.$
Thus we reach a contradiction,  since every matrix $\big(\nabla^{2}V(\xi_{i}^{*})\big)$ has $N$ orthogonal eigenvectors with $i=1,\cdots,k$.
Hence $\mathbf{b}_{l,i_{l}}=\mathbf{0}$.\\

Next we  claim that
 for each $i\in\{1,\cdots,k\}$, there exists
  only one eigenfunction $\omega_{\varepsilon,l_{i}}$, $l_{i}\in\{k(N+1)+1,\cdots,k(N+2)\}$
  such that  $\tilde{\omega}_{\varepsilon,l_{i},i}$ satisfies
  \begin{align}\label{898.6e2347}
\tilde{\omega}_{\varepsilon,l_{i},i}\rightarrow
b_{l_{i},i,0}\Psi_{0}
\; \; \text{in} \; \; C^{2}_{loc}\big(\mathbb{R}^N\big), \ \ \text{where}  \ \ b_{l_{i},i,0}\neq0 .\end{align}

We prove it by contradiction. If the claim were false, then  for some $i\in\{1,\cdots,k\}$, there exist at least two orthogonal eigenfunctions
$\omega_{\varepsilon,l_{i}^{\prime}}$ and  $\omega_{\varepsilon,l_{i}^{\prime\prime}}$ such that their corresponding  rescaled eigenfunctions
satisfy \eqref{898.6e2347}, where $
l_{i}^{\prime}\neq l_{i}^{\prime\prime}$. Since
\begin{align*}
\int_{\mathbb{R}^N} u_{\varepsilon}^{p-1-\varepsilon}(x)
\omega_{\varepsilon,l_{i}^{\prime}}(x)
\omega_{\varepsilon,l_{i}^{\prime\prime}}(x)\mathrm{d}x=0.
\end{align*}
Letting $\varepsilon\rightarrow0$, then we have
$$
b_{l_{i}^{\prime},i,0}b_{l_{i}^{\prime\prime},i,0}\int_{\mathbb{R}^N} U_{1,0}^{p-1}
\Psi_{0}^{2}=0,$$
which is  a contradiction. Thus our claim holds.
\vskip 0.1cm

If for some $l\in\big\{k(N+1)+1,\cdots,k(N+2)\big\}$, there exists another $i^\prime_{l}\in\{1,\cdots,k\}$ such that
\begin{align*}
\tilde{\omega}_{\varepsilon,l,i^\prime_{l}}\rightarrow
b_{l,i^\prime_{l},0}\Psi_{0}
\ \ \text{in} \ \ C^{2}_{loc}\big(\mathbb{R}^N\big).
\end{align*}
Then for some $i^*_{l}\in\{1,\cdots,k\}$, there exist  at least two orthogonal eigenfunctions verifying  \eqref{898.6e2347}, which is a contradiction and the proof of Lemma \ref{lem4.6e789} is finished.
\end{proof}
\begin{lemma}\label{lem4.p0}
For $ l=k(N+1)+1,\cdots,k(N+2)$, then we have
\begin{align*}%\label{y987}
\lambda_{\epsilon,l}=1+\tilde{M}_{l}
\varepsilon
+o(\varepsilon),
\end{align*}
where $\tilde{M}_{l}>0$ is a constant independent of  $\varepsilon$.
\end{lemma}
\begin{proof}[\bf Proof.]
By the similar arguments of \eqref{an} in Lemma \ref{lem4.4}, we readily obtain  the desired estimates.
\end{proof}
\begin{proof}[\bf Proof of Theorem \ref{athm1.2}]
It follows from Proposition \ref{Prop4.2}, Lemma \ref{lem4.04} and Lemma \ref{lem4.p0} that
\begin{equation}\label{l2}
\begin{cases}
\lambda_{\varepsilon,l}<1,
\ \ \ \ l=1,\cdots,k,
\\[3mm]
\lambda_{\varepsilon,l}=1+
\tilde{M}\beta_{i_l,m_l}\varepsilon^2+
o\big(\varepsilon^2\big), \ \ \ \ l=k+1,\cdots,k(N+1),
\\[3mm]
\lambda_{\epsilon,l}=1+\tilde{M}_{l}
\varepsilon
+o(\varepsilon),
\ \ l=k(N+1)+1,\cdots,k(N+2),
\end{cases}
\end{equation}
where $\tilde{M}>0$ and $\tilde{M}_{l}>0$ are  constants independent of $\varepsilon$ and $\beta_{i_l,m_l}$ is
an eigenvalue of  matrix $\big(\nabla^{2}V(\xi_{i_l}^{*})\big)$ as in Lemma \ref{lem4.04}.
%$\beta_{l-k}$ is
%an eigenvalue of some non-degenerate matrix $\big(\nabla^{2}V(\xi_{i}^{*})\big)$  which corresponds to the eigenvector $\mathbf{b}_{l,i}=\big(b_{l,i,1},\cdots,b_{l,i,N}\big)
%\neq(0,\cdots,0)$,
%$i=1,\cdots,k$, $l=k+1,\cdots,k(N+1)$.
Taking advantage of \eqref{l2}, we deduce that  the Morse index of $u_{\varepsilon}$ is
$$m(u_{\varepsilon})=\displaystyle\sum_{i=1}^{k}m
\big(\nabla^{2}V(\xi_{i}^{*})\big)+k.$$
Thus we complete the proof of Theorem \ref{athm1.2}.
\end{proof}
\begin{proof}[\bf Proof of Theorem \ref{athm1.3}]
By virtue of \eqref{l2}, we find  that $1$ is not the eigenvalue of problem \eqref{aa4.1}. Then the solution $u_\varepsilon$ obtained in Theorem A is non-degenerate.
\end{proof}

\noindent{\bf Author Contributions} All authors contributed equally to the writing of this paper and have reviewed it for submission.

 \smallskip
 
\noindent{\bf Funding}\, Z. Liu  was supported by the  Natural Science Foundation of Henan (Grant No. 252300421051) and
									the National Natural Science Foundation of China (Grant Nos. 12371111, 12471104). Z. Liu was also  supported by the Open Research Fund of Key Laboratory of Nonlinear Analysis
									\& Applications (Central China Normal University), Ministry of Education,  China (Grant No. NAA2024ORG004).
S. Tian was supported by the Fundamental Research Funds for the Central Universities (104972025KFYjc0115).

\smallskip
\noindent{\bf Data availability}\\
No datasets were generated or analysed during the current study.\\
%\begin{small}
%a
%\end{small}

\noindent{\bf Declarations }

\smallskip

\noindent{\bf Ethical Approval}  Not applicable.

\noindent{\bf Conflict of interest}\\
\smallskip
The authors declare that they have no conflicts of interest.
 %The authors declare that they have no conflict of interest to this work.\\
%The authors state that there is no conflict of interest. \\

\noindent{\bf Competing interests}\\
\smallskip
Competing interests The authors declare no competing interests.

%%%%%%%%%%%%%%%%%%%%%%%%%%%%%%%%%%%%%%%%%%%%%%%%%%%%%%%%%%%%%%%%

%%%%%%%%%%%%%%%%%%%%%%%%%%%%%%%%%%%%%%%%%%%%%%%%%%%%%%%%%%%%%%%%

\end{document}